\documentclass[generic]{imsart}

\RequirePackage[OT1]{fontenc}
\RequirePackage{amsthm,amsmath,amsfonts,amssymb,mathtools,thmtools}
\RequirePackage{natbib}
\RequirePackage[colorlinks,citecolor=blue,urlcolor=blue]{hyperref}
\RequirePackage{enumitem,graphicx,MnSymbol,mathrsfs}

\pubyear{2017}
\volume{X}
\issue{X}
\firstpage{1}
\lastpage{59}

\startlocaldefs

\numberwithin{equation}{section}
\theoremstyle{plain}
\newtheorem{stelling}{Theorem}[section]
\newtheorem{voorbeeld}[stelling]{Example}
\newtheorem{lemma}[stelling]{Lemma}
\theoremstyle{definition}                     
\newtheorem{definitie}{Definition}[section]   

\newcommand{\field}[1]{\mathbb{#1}} 
\newcommand{\R}{\field{R}} 

\newcommand{\refFig}[1]{Fig.~{\ref{#1}}}
\newcommand{\refEq}[1]{Eq.~{\eqref{#1}}}
\newcommand{\refTab}[1]{Table~{\ref{#1}}}
\newcommand{\refEx}[1]{Example~{\ref{#1}}}

\newcommand{\bibOrder}[1]{}

\newcommand{\Lc}{\mathcal{L}}

\newcommand{\Mc}{\mathcal{M}}
\newcommand{\Nc}{\mathcal{N}}

\newcommand{\Pc}{\mathcal{P}}

\newcommand{\Ps}{\mathscr{P}}

\newcommand{\Rc}{\mathcal{R}}

\newcommand{\Xc}{\mathcal{X}}
\newcommand{\Yc}{\mathcal{Y}}

\newcommand{\der}{\textnormal{d}}
\newcommand{\obs}{\textnormal{obs}}

\newcommand{\Var}{\text{Var}}
\newcommand{\Cov}{\text{Cov}}

\makeatletter
\newcommand{\xdashrightarrow}[2][]{\ext@arrow 0359\rightarrowfill@@{#1}{#2}}
\newcommand{\xdashleftarrow}[2][]{\ext@arrow 3095\leftarrowfill@@{#1}{#2}}
\newcommand{\xdashleftrightarrow}[2][]{\ext@arrow 3359\leftrightarrowfill@@{#1}{#2}}
\def\rightarrowfill@@{\arrowfill@@\relax\relbar\rightarrow}
\def\leftarrowfill@@{\arrowfill@@\leftarrow\relbar\relax}
\def\leftrightarrowfill@@{\arrowfill@@\leftarrow\relbar\rightarrow}
\def\downarrowfill@@{\arrowfill@@\leftarrow\relbar\rightarrow}
\def\arrowfill@@#1#2#3#4{%
  $\m@th\thickmuskip0mu\medmuskip\thickmuskip\thinmuskip\thickmuskip
   \relax#4#1
   \xleaders\hbox{$#4#2$}\hfill
   #3$%
}
\makeatother

\usepackage{tikz}
\usepackage{xparse}

\NewDocumentCommand\DownArrow{O{2.0ex} O{black}}{%
   \mathrel{\tikz[baseline] \draw [<-, line width=0.5pt, #2] (0,0) -- ++(0,#1);}
}

\endlocaldefs

\begin{document}

\begin{frontmatter}
\title{A Tutorial on Fisher Information\thanksref{t1}}
\runtitle{Fisher information tutorial}

\begin{aug}
\author{\fnms{Alexander} \snm{Ly}\ead[label=e1]{a.ly@uva.nl}},
\author{\fnms{Maarten} \snm{Marsman}\ead[label=e2]{m.marsman@uva.nl}}, 
\author{\fnms{Josine} \snm{Verhagen}\ead[label=e3]{josineverhagen@gmail.com}},
\author{\fnms{Raoul} \snm{Grasman}\ead[label=e4]{r.p.p.p.grasman@uva.nl}}
\and
\author{\fnms{Eric-Jan} \snm{Wagenmakers} \thanksref{t1} \ead[label=e5]{ej.wagenmakers@gmail.com}}
\ead[label=u1, url]{www.alexander-ly.com/}
\ead[label=u2, url]{https://jasp-stats.org/}

\thankstext{t1}{This work was supported by the starting grant ``Bayes or Bust'' awarded by the European Research Council (283876). Correspondence concerning this article may be addressed to Alexander Ly, email address: a.ly@uva.nl. The authors would like to thank Jay Myung, Trisha Van Zandt, and three anonymous reviewers for their comments on an earlier version of this paper. The discussions with Helen Steingroever, Jean-Bernard Salomond, Fabian Dablander, Nishant Mehta, Alexander Etz, Quentin Gronau and Sacha Epskamp led to great improvements of the manuscript. Moreover, the first author is grateful to Chris Klaassen, Bas Kleijn and Henk Pijls for their patience and enthusiasm with which they taught, and answered questions from a not very docile student.}
\runauthor{Ly, et. al.}

\affiliation{University of Amsterdam}

\address{University of Amsterdam \\
Department of Psychological Methods \\
PO Box 15906 \\
Nieuwe Achtergracht 129-B \\
1001 NK  Amsterdam \\
The Netherlands\\
e-mail: \printead*{e1} \\
url: \printead*{u1}; \printead*{u2}}
\end{aug}


\begin{abstract}
In many statistical applications that concern mathematical psychologists, the concept of Fisher information plays an important role. In this tutorial we clarify the concept of Fisher information as it manifests itself across three different statistical paradigms. First, in the frequentist paradigm, Fisher information is used to construct hypothesis tests and confidence intervals using maximum likelihood estimators; second, in the Bayesian paradigm, Fisher information is used to define a default prior; lastly, in the minimum description length paradigm, Fisher information is used to measure model complexity.
\end{abstract}

\begin{keyword}[class=MSC]
\kwd[Primary ]{62-01} 
\kwd{62B10} 
\kwd[; secondary ]{62F03} 
\kwd{62F12} 
\kwd{62F15} 
\kwd{62B10} 
\end{keyword}

\begin{keyword}
\kwd{Confidence intervals}
\kwd{hypothesis testing} 
\kwd{Jeffreys's prior} 
\kwd{minimum description length}
\kwd{model complexity}
\kwd{model selection} 
\kwd{statistical modeling}
\end{keyword}

\setcounter{tocdepth}{1}
%
\tableofcontents
\end{frontmatter}

\section{Introduction}

Mathematical psychologists develop and apply quantitative models in order to describe human behavior and understand latent psychological processes. Examples of such models include Stevens' law of psychophysics that describes the relation between the objective physical intensity of a stimulus and its subjectively experienced intensity (\citealp{stevens1957psychophysical}); Ratcliff's diffusion model of decision making that measures the various processes that drive behavior in speeded response time tasks (\citealp{ratcliff1978theory}); and multinomial processing tree models that decompose performance in memory tasks into the contribution of separate latent mechanisms (\citealp{batchelderRiefer1980}; \citealp{chechile1973}).

When applying their models to data, mathematical psychologists may operate from within different statistical paradigms and focus on different substantive questions. For instance, working within the classical or frequentist paradigm a researcher may wish to test certain hypotheses or decide upon the number of trials to be presented to participants in order to estimate their latent abilities. Working within the Bayesian paradigm a researcher may wish to know how to determine a suitable default prior on the parameters of a model. Working within the minimum description length (MDL) paradigm a researcher may wish to compare rival models and quantify their complexity. Despite the diversity of these paradigms and purposes, they are connected through the concept of Fisher information.

Fisher information plays a pivotal role throughout statistical modeling, but an accessible introduction for mathematical psychologists is lacking. The goal of this tutorial is to fill this gap and illustrate the use of Fisher information in the three statistical paradigms mentioned above: frequentist, Bayesian, and MDL. This work builds directly upon the \emph{Journal of Mathematical Psychology} tutorial article by \citet{myung2003tutorial} on maximum likelihood estimation. The intended target group for this tutorial are graduate students and researchers with an affinity for cognitive modeling and mathematical statistics.

To keep this tutorial self-contained we start by describing our notation and key concepts. We then provide the definition of Fisher information and show how it can be calculated. The ensuing sections exemplify the use of Fisher information for different purposes. Section~{\ref{fiInFreq}} shows how Fisher information can be used in frequentist statistics to construct confidence intervals and hypothesis tests from maximum likelihood estimators (MLEs). Section~{\ref{fiInBayes}} shows how Fisher information can be used in Bayesian statistics to define a default prior on model parameters. In Section~{\ref{fiInMdl}} we clarify how Fisher information can be used to measure model complexity within the MDL framework of inference.

\subsection{Notation and key concepts}
Before defining Fisher information it is necessary to discuss a series of fundamental concepts such as the nature of statistical models, probability mass functions, and statistical independence. Readers familiar with these concepts may safely skip to the next section.

A \emph{statistical model} is typically defined through a function \( f(x_{i} \, | \, \theta) \) that represents how a parameter \( \theta \) is functionally related to potential outcomes \( x_{i} \) of a random variable \( X_{i} \). For ease of exposition, we take \( \theta \) to be one-dimensional throughout this text. The generalization to vector-valued \( \theta \) can be found in Appendix~{\ref{appendixFiMatrix}}, see also \citet{myung2005information}. 

As a concrete example, \( \theta \) may represent a participant's intelligence, \( X_{i} \) a participant's (future) performance on the \( i \)th item of an IQ test, \( x_{i}=1 \) the potential outcome of a correct response, and \( x_{i}=0 \) the potential outcome of an incorrect response on the \( i \)th item. Similarly, \( X_{i} \) is the \( i \)th trial in a coin flip experiment with two potential outcomes: heads, \( x_{i}=1 \), or tails, \( x_{i}=0 \). Thus, we have the binary outcome space \( \Xc=\{ 0, 1 \} \). The coin flip model is also known as the Bernoulli distribution \( f(x_{i} \, | \, \theta) \) that relates the coin's propensity \( \theta \in (0, 1) \) to land heads to the potential outcomes as
\begin{align}
\label{BernoulliPdf}
f(x_{i} \, | \, \theta) = \theta^{x_{i}} (1- \theta)^{1-x_{i}}, \text{ where } x_{i} \in \Xc=\{ 0, 1 \}.
\end{align}
Formally, if \( \theta \) is known, fixing it in the functional relationship \( f \) yields a function \( p_{\theta}(x_{i})=f(x_{i} \, | \, \theta) \) of the potential outcomes \( x_{i} \). This \( p_{\theta}(x_{i}) \) is referred to as a \emph{probability density function} (pdf) when \( X_{i} \) has outcomes in a continuous interval, whereas it is known as a \emph{probability mass function} (pmf) when \( X_{i} \) has discrete outcomes. The pmf \( p_{\theta}(x_{i})=P(X_{i}=x_{i}\, | \, \theta) \) can be thought of as a data generative device as it specifies how \( \theta \) defines the chance with which \( X_{i} \) takes on a potential outcome \( x_{i} \). As this holds for any outcome \( x_{i} \) of \( X_{i} \), we say that \( X_{i} \) is distributed according to \( p_{\theta}(x_{i}) \). For brevity, we do not further distinguish the continuous from the discrete case, and refer to \( p_{\theta}(x_{i}) \) simply as a pmf.

For example, when the coin's true propensity is \( \theta^{*}=0.3 \), replacing \( \theta \) by \( \theta^{*} \) in the Bernoulli distribution yields the pmf \( p_{0.3}(x_{i})=0.3^{x_{i}}0.7^{1-x_{i}} \), a function of all possible outcomes of \( X_{i} \). A subsequent replacement \( x_{i}=0 \) in the pmf \( p_{0.3}(0)=0.7 \) tells us that this coin generates the outcome \( 0 \) with 70\% chance.

In general, experiments consist of \( n \) trials yielding a potential set of outcomes \( x^{n}=(x_{1}, \ldots, x_{n}) \) of the random vector \( X^{n}=(X_{1}, \ldots, X_{n}) \). These \( n \) random variables are typically assumed to be \emph{independent and identically distributed} (iid). Identically distributed implies that each of these \( n \) random variables is governed by one and the same \( \theta \), while independence implies that the joint distribution of all these \( n \) random variables simultaneously is given by a product, that is, %
\begin{align}
\label{iid}
f(x^{n} \, | \, \theta)=f(x_{1} \, | \, \theta) \times \ldots \times f(x_{n} \, | \, \theta)= \prod_{i=1}^{n} f(x_{i} \, | \, \theta).
\end{align}
As before, when \( \theta \) is known, fixing it in this relationship \( f(x^{n} \, | \, \theta) \) yields the (joint) pmf of \( X^{n} \) as \( p_{\theta}(x^{n})=p_{\theta}(x_{1}) \times \ldots \times p_{\theta}(x_{n}) = \prod_{i=1}^{n} p_{\theta}(x_{i}) \).

In psychology the iid assumption is typically evoked when experimental data are analyzed in which participants have been confronted with a sequence of \( n \) items of roughly equal difficulty. When the participant can be either correct or incorrect on each trial, the participant's performance \( X^{n} \) can then be related to an \( n \)-trial coin flip experiment governed by one single \( \theta \) over all \( n \) trials. The random vector \( X^{n} \) has \( 2^{n} \) potential outcomes \( x^{n} \). For instance, when \( n=10 \), we have \( 2^{n}=1{,}024 \) possible outcomes and we write \( \Xc^{n} \) for the collection of all these potential outcomes. The chance of observing a potential outcome \( x^{n} \) is determined by the coin's propensity \( \theta \) as follows
\begin{align}
\label{nBernoulliPdf}
f(x^{n} \, | \, \theta) =f(x_{1} \, | \, \theta) \times \ldots \times f(x_{n} \, | \, \theta)= \theta^{\sum_{i=1}^{n} x_{i}} (1 - \theta)^{n-\sum_{i=1}^{n} x_{i}}, \text{ where } x^{n} \in \Xc^{n}.
\end{align}
When the coin's true propensity \( \theta \) is \( \theta^{*}=0.6 \), replacing \( \theta \) by \( \theta^{*} \) in \refEq{nBernoulliPdf}
yields the joint pmf \( p_{0.6}(x^{n})=f(x^{n} \, | \, \theta=0.6) = 0.6^{\sum_{i=1}^{n} x_{i}} 0.4^{n- \sum_{i=1}^{n}x_{i}} \). The pmf with a particular outcome entered, say, \( x^{n}=(1,1,1,1,1,1,1,0,0,0) \) reveals that the coin with \( \theta^{*}=0.6 \) generates this particular outcome with 0.18\% chance. 

\subsection{Definition of Fisher information}
In practice, the true value of \( \theta \) is not known and has to be inferred from the observed data. The first step typically entails the creation of a data summary. For example, suppose once more that \( X^{n} \) refers to an \( n \)-trial coin flip experiment and suppose that we observed \( x^{n}_{\obs}=(1, 0, 0, 1, 1, 1, 1,0, 1, 1) \). To simplify matters, we only record the number of heads as \( Y=\sum_{i=1}^{n} X_{i} \), which is a function of the data. Applying our function to the specific observations yields the realization \( y_{\obs}=Y(x^{n}_{\obs})=7 \). Since the coin flips \( X^{n} \) are governed by \( \theta \), so is a function of \( X^{n} \); indeed, \( \theta \) relates to the potential outcomes \( y \) of \( Y \) as follows
\begin{align}
\label{binomialPdf}
f(y \, | \, \theta) = {n \choose y} \theta^{y} (1-\theta)^{n-y}, \text{ where } y \in \Yc = \{ 0, 1, \ldots, n \},
\end{align}
where \( {n \choose y} = {n! \over y! (n-y)!} \) enumerates the possible sequences of length \( n \) that consist of \( y \) heads and \(n - y\) tails. For instance, when flipping a coin \( n=10 \) times, there are 120 possible sequences of zeroes and ones that contain \( y=7 \) heads and \(n-y=3\) tails. The distribution \( f(y \, | \, \theta) \) is known as the binomial distribution.

The summary statistic \( Y \) has \( n+1 \) possible outcomes, whereas \( X^{n} \) has \( 2^{n} \). For instance, when \( n=10 \) the statistic \( Y \) has only \( 11 \) possible outcomes, whereas \( X^{n} \) has \( 1{,}024 \). This reduction results from the fact that the statistic \( Y \) ignores the order with which the data are collected. Observe that the conditional probability of the raw data given \( Y=y \) is equal to \( P(X^{n} \, | \, Y=y, \theta)=1 / {n \choose y} \) and that it does not depend on \( \theta \). This means that after we observe \( Y=y \) the conditional probability of \( X^{n} \) is independent of \( \theta \), even though each of the distributions of \( X^{n} \) and \( Y \) separately do depend on \( \theta \). We, therefore, conclude that there is no information about \( \theta \) left in \( X^{n} \) after observing \( Y=y \) (\citealp{fisher1920mathematical}; \citealp{stigler1973studies}). 

More generally, we call a function of the data, say, \( T=t(X^{n}) \) a \emph{statistic}. A statistic is referred to as \emph{sufficient} for the parameter \( \theta \), if the expression \( P(X^{n} \, | \, T=t, \theta) \) does not depend on \( \theta \) itself. To quantify the amount of information about the parameter \( \theta \) in a sufficient statistic \( T \) and the raw data, Fisher introduced the following measure.

\begin{definitie}[Fisher information]
The \emph{Fisher information} \( I_{X}(\theta) \) of a random variable \( X \) about \( \theta \) is defined as%
\footnote{Under mild regularity conditions Fisher information is equivalently defined as %
\begin{align}
\label{FIOne}
I_{X}(\theta) = -E \Big ( \tfrac{\der^2 }{\der \theta ^2} \log f (X \, | \, \theta) \Big ) = \begin{cases} -\sum_{x \in \Xc} \Big ( \tfrac{\der^2 }{\der \theta^2} \log f (x \, | \, \theta) \Big ) p_{\theta}(x) & \text{if } X \text{ is discrete,} \\
-\int_{\Xc} \Big ( \tfrac{\der^2 }{\der \theta ^2} \log f (x \, | \, \theta) \Big ) p_{\theta}(x) \der x & \text{if } X \text{ is continuous.}
\end{cases}
\end{align}
where \( \tfrac{\der^2 }{\der \theta ^2} \log f (x \, | \, \theta) \) denotes the second derivate of the logarithm of \( f \) with respect \( \theta \).} %
\begin{align}
\label{FIOne0}
I_{X}(\theta) = 
\begin{cases} \sum_{x \in \Xc} \Big ( \tfrac{\der }{\der \theta} \log f (x \, | \, \theta) \Big )^{2} p_{\theta} (x) & \text{if } X \text{ is discrete,} \\
\int_{\Xc} \Big ( \tfrac{\der }{\der \theta} \log f (x \, | \, \theta) \Big )^{2} p_{\theta} (x)\der x & \text{if } X \text{ is continuous.}
\end{cases}
\end{align}
The derivative \( \tfrac{\der }{\der \theta} \log f(x \, | \, \theta) \) is known as the \emph{score function}, a function of \( x \), and describes how sensitive the model (i.e., the functional form \( f \)) is to changes in \( \theta \) at a particular \( \theta \). The Fisher information measures the overall sensitivity of the functional relationship \( f \) to changes of \( \theta \) by weighting the sensitivity at each potential outcome \( x \) with respect to the chance defined by \( p_{\theta}(x)=f(x \, | \, \theta) \). The weighting with respect to \( p_{\theta}(x) \) implies that the Fisher information about \( \theta \) is an expectation.

Similarly, Fisher information \( I_{X^{n}}(\theta) \) within the random vector \( X^{n} \) about \( \theta \) is calculated by replacing \( f(x \, | \, \theta) \) with \( f(x^{n} \, | \, \theta) \), thus, \( p_{\theta}(x) \) with \( p_{\theta}(x^{n}) \) in the definition. Moreover, under the assumption that the random vector \( X^{n} \) consists of \( n \) iid trials of \( X \) it can be shown that \( I_{X^{n}}(\theta)=n I_{X}(\theta) \), which is why \( I_{X}(\theta) \) is also known as the unit Fisher information.%
\footnote{Note the abuse of notation -- we dropped the subscript \( i \) for the \( i \)th random variable \( X_{i} \) and denote it simply by \( X \) instead.} %
Intuitively, an experiment consisting of \( n=10 \) trials is expected to be twice as informative about \( \theta \) compared to an experiment consisting of only \( n=5 \) trials. \( \hfill \diamond \)
\end{definitie}

Intuitively, we cannot expect an arbitrary summary statistic \( T \) to extract more information about \( \theta \) than what is already provided by the raw data. Fisher information adheres to this rule, as it can be shown that
\begin{align}
\label{informationInequality}
I_{X^{n}}(\theta) \geq I_{T}(\theta),
\end{align}
with equality if and only if \( T \) is a sufficient statistic for \( \theta \). 

\begin{voorbeeld}[name=The information about $\theta$ within the raw data and a summary statistic]
A direct calculation with a Bernoulli distributed random vector \( X^{n} \) shows that the Fisher information about \( \theta \) within an \( n \)-trial coin flip experiment is given by
\begin{align}
\label{FIBernoulli}
I_{X^{n}}(\theta)=n I_{X}(\theta) = n {1 \over \theta (1 - \theta)},
\end{align}
where \( I_{X}(\theta)=\tfrac{1}{\theta(1-\theta)} \) is the Fisher information of \( \theta \) within a single trial. As shown in \refFig{fig:FIBernoulli}, the unit Fisher information \( I_{X}(\theta) \) depends on \( \theta \). %
\begin{figure}
\centering
\includegraphics[width=9cm]{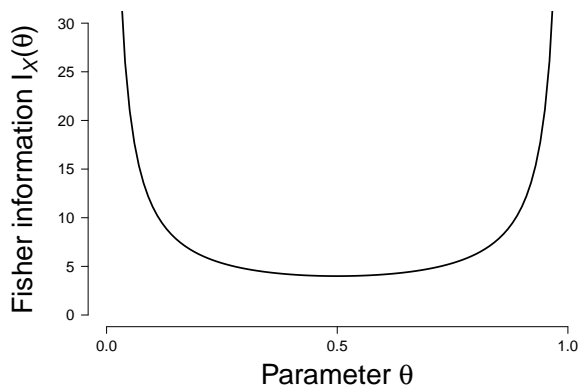}
\caption{The unit Fisher information \( I_{X}(\theta)={1 \over \theta (1- \theta)} \) as a function of \( \theta \) within the Bernoulli model. As \( \theta \) reaches zero or one the expected information goes to infinity.}
\label{fig:FIBernoulli}
\end{figure}
Similarly, we can calculate the Fisher information about \( \theta \) within the summary statistic \( Y \) by using the binomial model instead. This yields \( I_{Y}(\theta)={n \over \theta (1-\theta)} \). Hence, \( I_{X^{n}}(\theta)=I_{Y}(\theta) \) for any value of \( \theta \). In other words, the expected information in \( Y \) about \( \theta \) is the same as the expected information about \( \theta \) in \( X^{n} \), regardless of the value of \( \theta \).
\( \hfill \diamond \)
\end{voorbeeld}
Observe that the information in the raw data \( X^{n} \) and the statistic \( Y \) are equal for every \( \theta \), and specifically also for its unknown true value \( \theta^{*} \). That is, there is no statistical information about \( \theta \) lost when we use a sufficient statistic \( Y \) instead of the raw data \( X^{n} \). This is particular useful when the data set \( X^{n} \) is large and can be replaced by single number \( Y \).

\section{The Role of Fisher Information in Frequentist Statistics}
\label{fiInFreq}
Recall that \( \theta \) is unknown in practice and to infer its value we might: (1) provide a best guess in terms of a point estimate; (2) postulate its value and test whether this value aligns with the data, or (3) derive a confidence interval. In the frequentist framework, each of these inferential tools is related to the Fisher information and exploits the data generative interpretation of a pmf. Recall that given a model \( f(x^{n} \, | \, \theta) \) and a known \( \theta \), we can view the resulting pmf \( p_{\theta}(x^{n}) \) as a recipe that reveals how \( \theta \) defines the chances with which \( X^{n} \) takes on the potential outcomes \( x^{n} \). 

\begin{figure}
\centering
\includegraphics[width = 9cm]{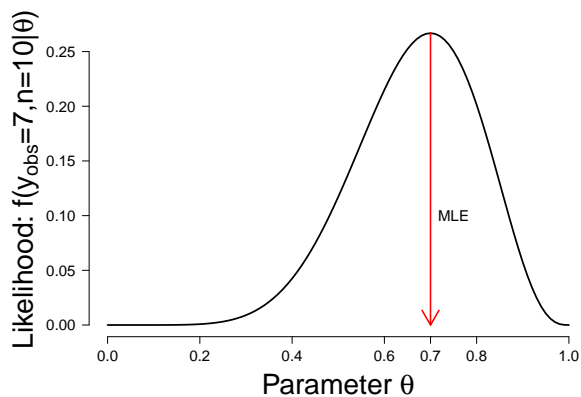}
\caption{The likelihood function based on observing \( y_{\obs}=7 \) heads in \( n=10 \) trials. For these data, the MLE is equal to \( \hat{\theta}_{\obs}=0.7 \), see the main text for the interpretation of this function.}
\label{figLikelihoodBin}
\end{figure}

This data generative view is central to Fisher's conceptualization of the \emph{maximum likelihood estimator} (MLE; \citealp{fisher1912absolute}; \citealp{fisher1922mathematical}; \citealp{fisher1925theory}; \citealp{le1990maximum}; \citealp{myung2003tutorial}). For instance, the binomial model implies that a coin with a hypothetical propensity \( \theta=0.5 \) will generate the outcome \( y=7 \) heads out of \( n=10 \) trials with 11.7\% chance, whereas a hypothetical propensity of \( \theta=0.7 \) will generate the same outcome \( y=7 \) with 26.7\% chance. Fisher concluded that an actual observation \( y_{\obs} = 7 \) out of \( n=10 \) is therefore more likely to be generated from a coin with a hypothetical propensity of \( \theta=0.7 \) than from a coin with a hypothetical propensity of \( \theta=0.5 \). \refFig{figLikelihoodBin} shows that for this specific observation \( y_{\obs} = 7 \), the hypothetical value \( \theta=0.7 \) is the maximum likelihood \emph{estimate}; the number \( \hat{\theta}_{\obs}=0.7 \). %
This estimate is a realization of the maximum likelihood \emph{estimator} (MLE); in this case, the MLE is the function \( \hat{\theta}=\tfrac{1}{n} \sum_{i=1}^{n} X_{i} = \tfrac{1}{n} Y \), i.e., the sample mean. Note that the MLE is a statistic, that is, a function of the data.

\subsection{Using Fisher information to design an experiment}
Since \( X^{n} \) depends on \( \theta \) so will a function of \( X^{n} \), in particular, the MLE \( \hat{\theta} \). The distribution of the potential outcomes of the MLE \( \hat{\theta} \) is known as the \emph{sampling distribution} of the estimator and denoted as \( f(\hat{\theta}_{\obs} \, | \, \theta) \). As before, when \( \theta^{*} \) is assumed to be known, fixing it in \( f(\hat{\theta}_{\obs} \, | \, \theta) \) yields the pmf \( p_{\theta^{*}}(\hat{\theta}_{\obs}) \), a function of the potential outcomes of \( \hat{\theta} \). This function \( f \) between the parameter \( \theta \) and the potential outcomes of the MLE \( \hat{\theta} \) is typically hard to describe, but for \( n \) large enough it can be characterized by the Fisher information.

For iid data and under general conditions,%
\footnote{Basically, when the Fisher information exists for all parameter values. For details see the advanced accounts provided by \citet{bickel1993efficient}, \citet{hajek1970characterization}, \citet{inagaki1970limiting}, \citet{lecam1970assumptions} and Appendix~{\ref{appendixReg}}.} %
the difference between the true \( \theta^{*} \) and the MLE converges in distribution to a normal distribution, that is, %
\begin{align}
\label{mleCLT}
\sqrt{n} ( \hat{\theta} - \theta^{*}) \overset{D}{\rightarrow} \Nc(0, I_{X}^{-1}(\theta^{*})), \text{ as n } \rightarrow \infty.
\end{align}
Hence, for large enough \( n \), the ``error'' is approximately normally distributed%
\footnote{Note that \( \hat{\theta} \) is random, while the true value \( \theta^{*} \) is fixed. As such, the error \( \hat{\theta}-\theta^{*} \) and the rescaled error \( \sqrt{n} ( \hat{\theta} - \theta^{*}) \) are also random. We used \( \overset{D}{\rightarrow} \) in \refEq{mleCLT} to convey that the distribution of the left-hand side goes to the distribution on the right-hand side. Similarly, \( \overset{D}{\approx} \) in \refEq{mleCLT2} implies that the distribution of the left-hand side is \emph{approximately} equal to the distribution given on the right-hand side. Hence, for finite \( n \) there will be an error due to using the normal distribution as an approximation to the true sampling distribution. This approximation error is ignored in the constructions given below, see Appendix~{\ref{appendixFreq}} for a more thorough discussion.} %
\begin{align}
\label{mleCLT2}
(\hat{\theta} - \theta^{*}) \overset{D}{\approx} \Nc \Big (0, 1/ \big (n I_{X}(\theta^{*}) \big ) \Big ).
\end{align}
This means that the MLE \( \hat{\theta} \) generates potential estimates \( \hat{\theta}_{\obs} \) around the true value \( \theta^{*} \) with a standard error given by the inverse of the square root of the Fisher information at the true value \( \theta^{*} \), i.e., \( 1/\sqrt{n I_{X}(\theta^{*})} \), whenever \( n \) is large enough. Note that the chances with which the estimates of \( \hat{\theta} \) are generated depend on the true value \( \theta^{*} \) and the sample size \( n \). Observe that the standard error decreases when the unit information \( I_{X}(\theta^{*}) \) is high or when \( n \) is large. As experimenters we do not have control over the true value \( \theta^{*} \), but we can affect the data generating process by choosing the number of trials \( n \). Larger values of \( n \) increase the amount of information in \( X^{n} \), heightening the chances of the MLE producing an estimate \( \hat{\theta}_{\obs} \) that is close to the true value \( \theta^{*} \). The following example shows how this can be made precise.

\begin{voorbeeld}[name=Designing a binomial experiment with the Fisher information, label=designBin]
Recall that the potential outcomes of a normal distribution fall within one standard error of the population mean with 68\% chance. Hence, when we choose \( n \) such that \( 1/\sqrt{nI_{X}(\theta^{*})}=0.1 \) we design an experiment that allows the MLE to generate estimates within \( 0.1 \) distance of the true value with 68\% chance. To overcome the problem that \( \theta^{*} \) is not known, we solve the problem for the worst case scenario. For the Bernoulli model this is given by \( \theta=1/2 \), the least informative case, see \refFig{fig:FIBernoulli}. As such, we have \( {1 / \sqrt{n I_{X}(\theta^{*})}} \leq {1 / \sqrt{n I_{X}(1/2)}} = {1 / (2 \sqrt{n})} = 0.1 \), where the last equality is the target requirement and is solved by \( n=25 \).

This leads to the following interpretation. After simulating \( k=100 \) data sets \( x^{n}_{\obs, 1}, \ldots, x^{n}_{\obs, k} \) each with \( n=25 \) trials, we can apply to each of these data sets the MLE yielding \( k \) estimates \( \hat{\theta}_{\obs, 1}, \ldots, \hat{\theta}_{\obs, k} \). The sampling distribution implies that at least \( 68 \) of these \( k=100 \) estimate are expected to be at most \( 0.1 \) distance away from the true \( \theta^{*} \). \( \hfill \diamond \)
\end{voorbeeld}

\subsection{Using Fisher information to construct a null hypothesis test}
The (asymptotic) normal approximation to the sampling distribution of the MLE can also be used to construct a null hypothesis test. When we postulate that the true value equals some hypothesized value of interest, say, \( \theta^{*}=\theta_{0} \), a simple plugin then allows us to construct a prediction interval based on our knowledge of the normal distribution. More precisely, the potential outcomes \( x^{n} \) with \( n \) large enough and generated according to \( p_{\theta^{*}}(x^{n}) \) leads to potential estimates \( \hat{\theta}_{\obs} \) that fall within the range
\begin{align}
\label{nonRejectionInterval}
\left (\theta^{*}-1.96 \sqrt{ \tfrac{1}{n} I_{X}^{-1}(\theta^{*})}, \theta^{*} +1.96 \sqrt{ \tfrac{1}{n} I_{X}^{-1}(\theta^{*})} \right ),
\end{align}
with (approximately) 95\% chance. This 95\%-prediction interval \refEq{nonRejectionInterval} allows us to construct a point null hypothesis test based on a pre-experimental postulate \( \theta^{*}=\theta_{0} \).

\begin{voorbeeld}[name=A null hypothesis test for a binomial experiment, label=bernoulliH0Ex]
Under the null hypothesis \( H_{0}: \theta^{*}=\theta_{0} = 0.5\), we predict that an outcome of the MLE based on \( n=10 \) trials will lie between \( (0.19, 0.81) \) with 95\% chance. This interval follows from replacing \( \theta^{*} \) by \( \theta_{0} \) in the 95\%-prediction interval \refEq{nonRejectionInterval}. The data generative view implies that if we simulate \( k=100 \) data sets each with the same \( \theta^{*}=0.5 \) and \( n=10 \), we would then have \( k \) estimates \( \hat{\theta}_{\obs, 1}, \ldots, \hat{\theta}_{\obs, k} \) of which five are expected to be outside this 95\% interval \( (0.19, 0.81) \). Fisher, therefore, classified an outcome of the MLE that is smaller than 0.19 or larger than 0.81 as extreme under the null and would then reject the postulate \( H_{0}: \theta_{0}=0.5 \) at a significance level of \( .05 \). \( \hfill \diamond \)
\end{voorbeeld}
The normal approximation to the sampling distribution of the MLE and the resulting null hypothesis test is particularly useful when the exact sampling distribution of the MLE is unavailable or hard to compute. %
\begin{voorbeeld}[name=An MLE null hypothesis test for the Laplace model, label=laplaceEx]
Suppose that we have \( n \) iid samples from the Laplace distribution %
\begin{align}
f(x_{i} \, | \, \theta) = \tfrac{1}{2 b} \exp \big ( - \tfrac{| x_{i} - \theta |}{b} \big ),
\end{align}
where \( \theta \) denotes the population mean and the population variance is given by \( 2b^2\). It can be shown that the MLE for this model is the sample median, \( \hat{\theta}=\hat{M} \), and the unit Fisher information is \( I_{X}(\theta)=b^{-2} \). The exact sampling distribution of the MLE is unwieldy \citep{kotz2001laplace} and not presented here. Asymptotic normality of the MLE is practical, as it allows us to discard the unwieldy exact sampling distribution and, instead, base our inference on a more tractable (approximate) normal distribution with a mean equal to the true value \(\theta^{*} \) and a variance equal to \( {b^{2} / n} \). For \( n=100 \), \( b=1 \) and repeated sampling under the hypothesis \( H_{0}:\theta^{*}=\theta_{0} \), approximately 95\% of the estimates (the observed sample medians) are expected to fall in the range \( (\theta_{0} - 0.196, \theta_{0}+0.196) \). \( \hfill \diamond \)
\end{voorbeeld}

\subsection{Using Fisher information to compute confidence intervals}
An alternative to both point estimation and null hypothesis testing is interval estimation. In particular, a 95\%-confidence interval can be obtained by replacing in the prediction interval \refEq{nonRejectionInterval} the unknown true value \( \theta^{*} \) by an estimate \( \hat{\theta}_{\obs} \). Recall that a simulation with \( k=100 \) data sets each with \( n \) trials leads to \( \hat{\theta}_{\obs, 1}, \ldots, \hat{\theta}_{\obs, k} \) estimates, and each estimate leads to a different 95\%-confidence interval. 
It is then expected that \( 95 \) of these \( k=100 \) intervals encapsulate the true value \( \theta^{*} \).%
\footnote{But see \citet{brown2001interval}.} %
Note that these intervals are centred around different points whenever the estimates differ and that their lengths differ, as the Fisher information depends on \( \theta \). %

\begin{voorbeeld}[name=An MLE confidence interval for the Bernoulli model, label=bernoulliCIEx]
When we observe \( y_{\obs, 1}=7 \) heads in \( n=10 \) trials, the MLE then produces the estimate \( \hat{\theta}_{\obs, 1}=0.7 \). Replacing $\theta^\ast$ in the prediction interval \refEq{nonRejectionInterval} with \( \theta^{*}=\hat{\theta}_{\obs, 1} \) yields an approximate 95\%-confidence interval \( (0.42, 0.98) \) of length \( 0.57 \). On the other hand, had we instead observed \(y_{\obs, 2}=6 \) heads, the MLE would then yield \(\hat{\theta}_{\obs, 2} = 0.6 \) resulting in the interval \( (0.29, 0.90) \) of length \( 0.61 \). \( \hfill \diamond \)
\end{voorbeeld}

In sum, Fisher information can be used to approximate the sampling distribution of the MLE when \( n \) is large enough. Knowledge of the Fisher information can be used to choose \( n \) such that the MLE produces an estimate close to the true value, construct a null hypothesis test, and compute confidence intervals.

\section{The Role of Fisher Information in Bayesian Statistics}
\label{fiInBayes}
This section outlines how Fisher information can be used to define the Jeffreys's prior, a default prior commonly used for estimation problems and for nuisance parameters in a Bayesian hypothesis test (e.g., \citealp{bayarri2012criteria}; \citealp{dawid2011posterior}; \citealp{gronau2017informed}; \citealp{jeffreys1961theory}; \citealp{liang2008mixtures}; \citealp{li2015mixtures}; \citealp{ly2016harold, ly2016evaluation}; \citealp{lyInpressanalytic}; \citealp{ly2017bayesian}; \citealp{robert2016expected}). To illustrate the desirability of the Jeffreys's prior we first show how the naive use of a uniform prior may have undesirable consequences, as the uniform prior depends on the representation of the inference problem, that is, on how the model is parameterized. This dependence is commonly referred to as lack of invariance: different parameterizations of the same model result in different posteriors and, hence, different conclusions. We visualize the representation problem using simple geometry and show how the geometrical interpretation of Fisher information leads to the Jeffreys's prior that is parameterization-invariant. 

\subsection{Bayesian updating}
Bayesian analysis centers on the observations \( x^{n}_{\obs} \) for which a generative model \( f \) is proposed that functionally relates the observed data to an unobserved parameter \( \theta \). Given the observations \( x^{n}_{\obs} \), the functional relationship \( f \) is inverted using Bayes' rule to infer the relative plausibility of the values of \( \theta \). This is done by replacing the potential outcome part \( x^{n} \) in \( f \) by the actual observations yielding a \emph{likelihood function} \( f(x^{n}_{\obs} \, | \, \theta) \), which is a function of \( \theta \). %
In other words, \( x^{n}_{\obs} \) is known, thus, fixed, and the true \( \theta \) is unknown, therefore, free to vary. 
The candidate set of possible values for the true \( \theta \) is denoted by \( \Theta \) and referred to as the parameter space. Our knowledge about \( \theta \) is formalized by a distribution \( g(\theta) \) over the parameter space \( \Theta \). This distribution is known as the prior on \( \theta \), as it is set before any datum is observed. We can use Bayes' theorem to calculate the posterior distribution over the parameter space \( \Theta \) given the data that were actually observed as follows
\begin{align}
\label{BayesTheorem}
g(\theta \, | \, X^{n}=x^{n}_{\obs}) = \frac{f(x^{n}_{\obs} \, | \, \theta) g(\theta)}{\int_{\Theta} f(x^{n}_{\obs} \, | \, \theta) g(\theta) \, \der \theta}.
\end{align}
This expression is often verbalized as %
\begin{align}
\label{BayesTheoremVerbalised}
\text{posterior} = \frac{\text{likelihood} \times \text{prior}}{\text{marginal likelihood}}.
\end{align}%
The posterior distribution is a combination of what we knew before we saw the data (i.e., the information in the prior), and what we have learned from the observations in terms of the likelihood (e.g., \citealp{lee2013bayesian}). Note that the integral is now over \( \theta \) and not over the potential outcomes.

\subsection{Failure of the uniform distribution on the parameter as a noninformative prior}
When little is known about the parameter \( \theta \) that governs the outcomes of \( X^{n} \), it may seem reasonable to express this ignorance with a uniform prior distribution \( g(\theta) \), as no parameter value of \( \theta \) is then favored over another. This leads to the following type of inference:

\begin{voorbeeld}[Uniform prior on $\theta$]
\label{bayesBer}
Before data collection, \( \theta \) is assigned a uniform prior, that is, \( g(\theta)= {1 / V_{\Theta}} \) with a normalizing constant of \( V_{\Theta}=1 \) as shown in the left panel of \refFig{posteriorTheta}. %
\begin{figure}
    \centering
\begin{tabular}{ccc}
    \begin{minipage}{.375 \textwidth}
    \centering
    \qquad \scalebox{1.25}{Uniform prior on \( \theta \)}
    \includegraphics[width= \linewidth]{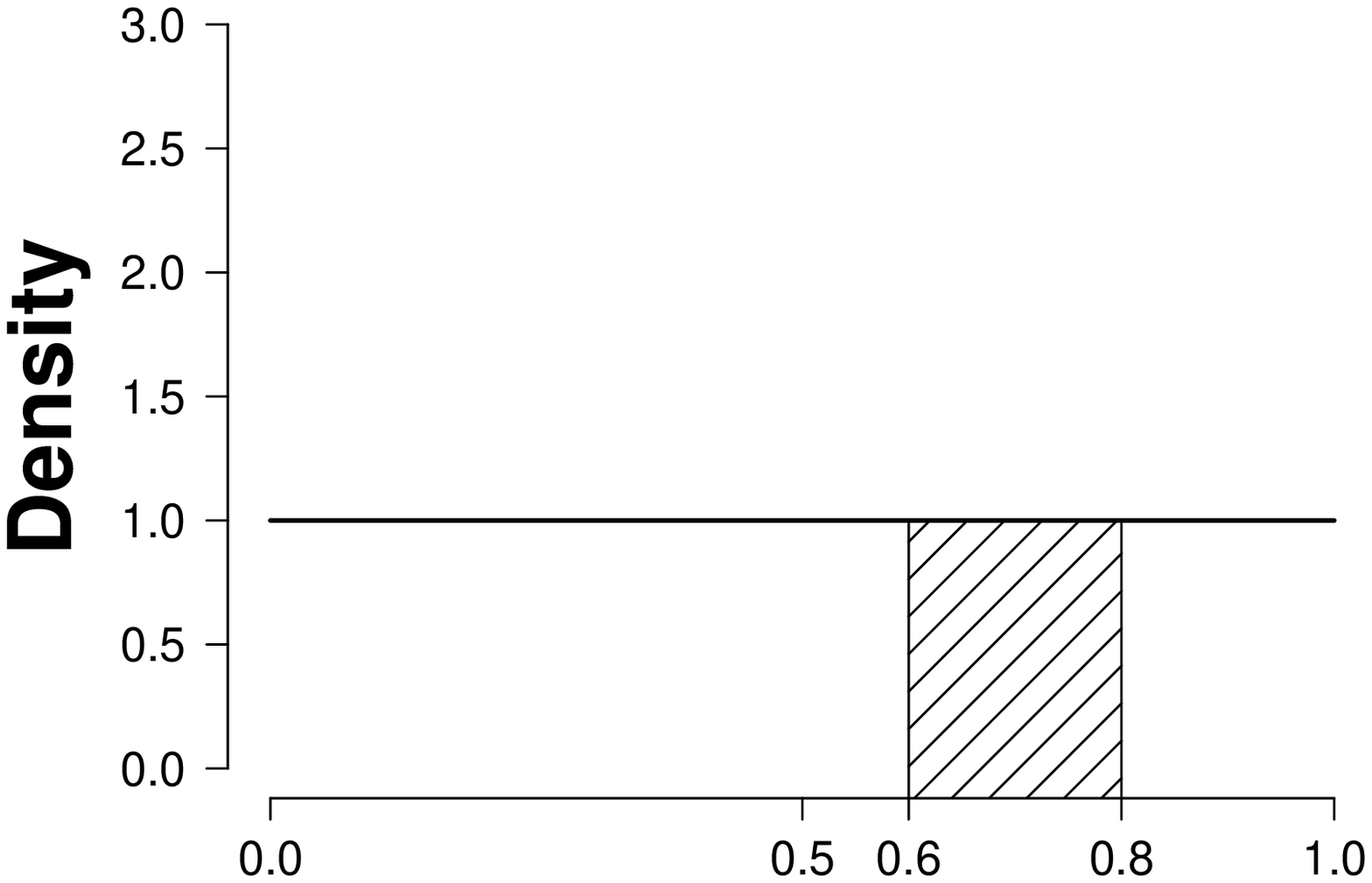}
    
    \quad \qquad \scalebox{1.25}{Propensity \( \theta \)}
    \end{minipage}  & %
    \scalebox{1.25}{\( \xrightarrow[n \, = \, 10]{y_{\obs} \, = \, 7} \)} & %
    \begin{minipage}{.375 \textwidth}
    \centering
    \scalebox{1.25}{Posterior \( \theta \) from \( \theta \sim U[0,1] \)}
    \includegraphics[width=\linewidth]{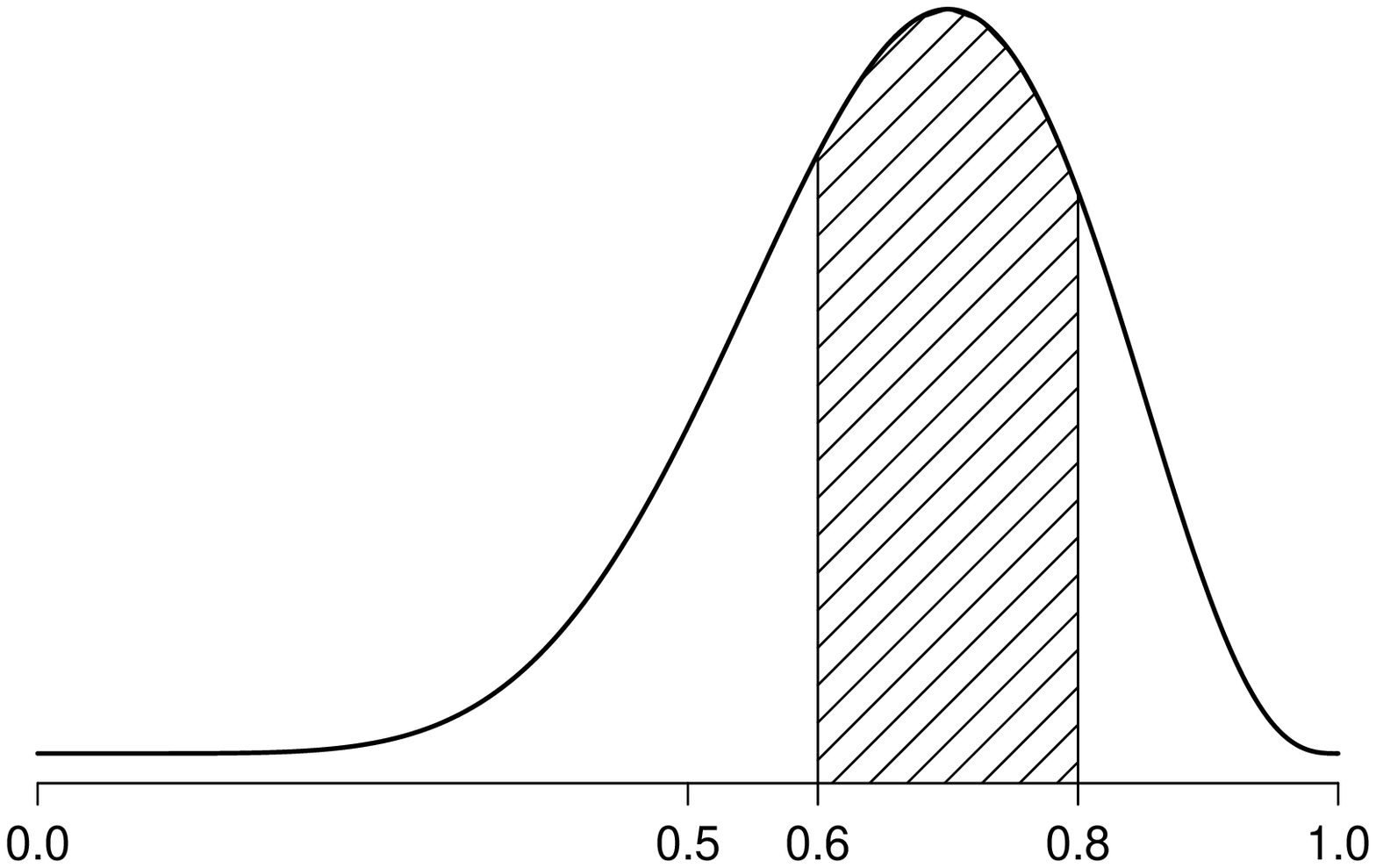}
    \scalebox{1.25}{Propensity \( \theta \)}
    \end{minipage}  %
\end{tabular}    
\caption{Bayesian updating based on observations \( x^{n}_{\obs} \) with \( y_{\obs}=7 \) heads out of \( n=10 \) tosses. In the left panel, the uniform prior distribution assigns equal probability to every possible value of the coin's propensity \( \theta \). In the right panel, the posterior distribution is a compromise between the prior and the observed data.}  
\label{posteriorTheta}
\end{figure}  
Suppose that we observe coin flip data \( x^{n}_{\obs} \) with \( y_{\obs}=7 \) heads out of \( n=10 \) trials. To relate these observations to the coin's propensity \( \theta \) we use the Bernoulli distribution as our \( f(x^{n} \, | \, \theta) \). A replacement of \( x^{n} \) by the data actually observed yields the likelihood function \( f(x^{n}_{\obs} \, | \, \theta)= \theta^{7} (1- \theta)^{3} \), which is a function of \( \theta \). Bayes' theorem now allows us to update our prior to the posterior that is plotted in the right panel of \refFig{posteriorTheta}. \( \hfill \diamond \)
 \end{voorbeeld}
Note that a uniform prior on \( \theta \) has the length, more generally, volume, of the parameter space as the normalizing constant; in this case, \( V_{\Theta}=1 \), which equals the length of the interval \( \Theta=(0,1) \). Furthermore, a uniform prior can be characterized as the prior that gives equal probability to all sub-intervals of equal length. Thus, the probability of finding the true value \( \theta^{*} \) within a sub-interval \( J_{\theta}=(\theta_{a}, \theta_{b}) \subset \Theta=(0, 1) \) is given by the relative length of \( J_{\theta} \) with respect to the length of the parameter space, that is,
\begin{align}
\label{uniformGeo}
 P \Big (\theta^{*} \in J_{\theta} \Big )= \int_{J_{\theta}} g(\theta) \der \theta = \frac{1}{V_{\Theta}} \int_{\theta_{a}}^{\theta_{b}}1 \der \theta = \frac{\theta_{b}-\theta_{a}}{V_{\Theta}}.
\end{align}
Hence, before any datum is observed, the uniform prior expresses the belief \( P (\theta^{*} \in J_{\theta})=0.20 \) of finding the true value \( \theta^{*} \) within the interval \( J_{\theta}=(0.6, 0.8) \). After observing \( x^{n}_{\obs} \) with \( y_{\obs}=7 \) out of \( n=10 \), this prior is updated to the posterior belief of \( P(\theta^{*} \in J_{\theta} \, | \, x^{n}_{\obs})=0.54 \), see the shaded areas in \refFig{posteriorTheta}.

Although intuitively appealing, it can be unwise to choose the uniform distribution by default, as the results are highly dependent on how the model is parameterized. In what follows, we show how a different parameterization leads to different posteriors and, consequently, different conclusions.

\begin{voorbeeld}[Different representations, different conclusions]
\label{coinBernoulliEx}
The propensity of a coin landing heads up is related to the angle \( \phi \) with which that coin is bent. Suppose that the relation between the angle \( \phi \) and the propensity \( \theta \) is given by the function \( \theta=h(\phi)=\frac{1}{2} + \frac{1}{2} \big ( \tfrac{\phi}{\pi} \big )^{3} \), chosen here for mathematical convenience.%
\footnote{Another example involves the logit formulation of the Bernoulli model, that is, in terms of \( \phi = \log (\tfrac{\theta}{1- \theta} ) \), where \( \Phi = \R \). This logit formulation is the basic building block in item response theory. We did not discuss this example as the uniform prior on the logit cannot be normalized and, therefore, not easily represented in the plots.} %
When \( \phi \) is positive the tail side of the coin is bent inwards, which increases the coin's chances to land heads. As the function \( \theta=h(\phi) \) also admits an inverse function \( h^{-1}(\theta)=\phi \), we have an equivalent formulation of the problem in \refEx{bayesBer}, but now described in terms of the angle \( \phi \) instead of the propensity \( \theta \).

As before, in order to obtain a posterior distribution, Bayes' theorem requires that we specify a prior distribution. As the problem is formulated in terms of \( \phi \), one may believe that a noninformative choice is to assign a uniform prior \( \tilde{g}(\phi) \) on \( \phi \), as this means that no value of \( \phi \) is favored over another. A uniform prior on \( \phi \) is in this case given by \( \tilde{g}(\phi)= 1/ V_{\Phi} \) with a normalizing constant \( V_{\Phi}=2 \pi \), because the parameter \( \phi \) takes on values in the interval \( \Phi=(-\pi, \pi) \). This uniform distribution expresses the belief that the true \( \phi^{*} \) can be found in any of the intervals \( (-1.0 \pi, -0.8 \pi), (-0.8 \pi, -0.6 \pi), \ldots, (0.8 \pi, 1.0 \pi) \) with 10\% probability, because each of these intervals is 10\% of the total length, see the top-left panel of \refFig{inductionFromUniformPhi}. %
\begin{figure}[h]
\begin{tabular}{ccc}
    \begin{minipage}{.375 \textwidth}
    \centering
    \qquad \scalebox{1.25}{Uniform prior on \( \phi \)}
    \includegraphics[width= \linewidth]{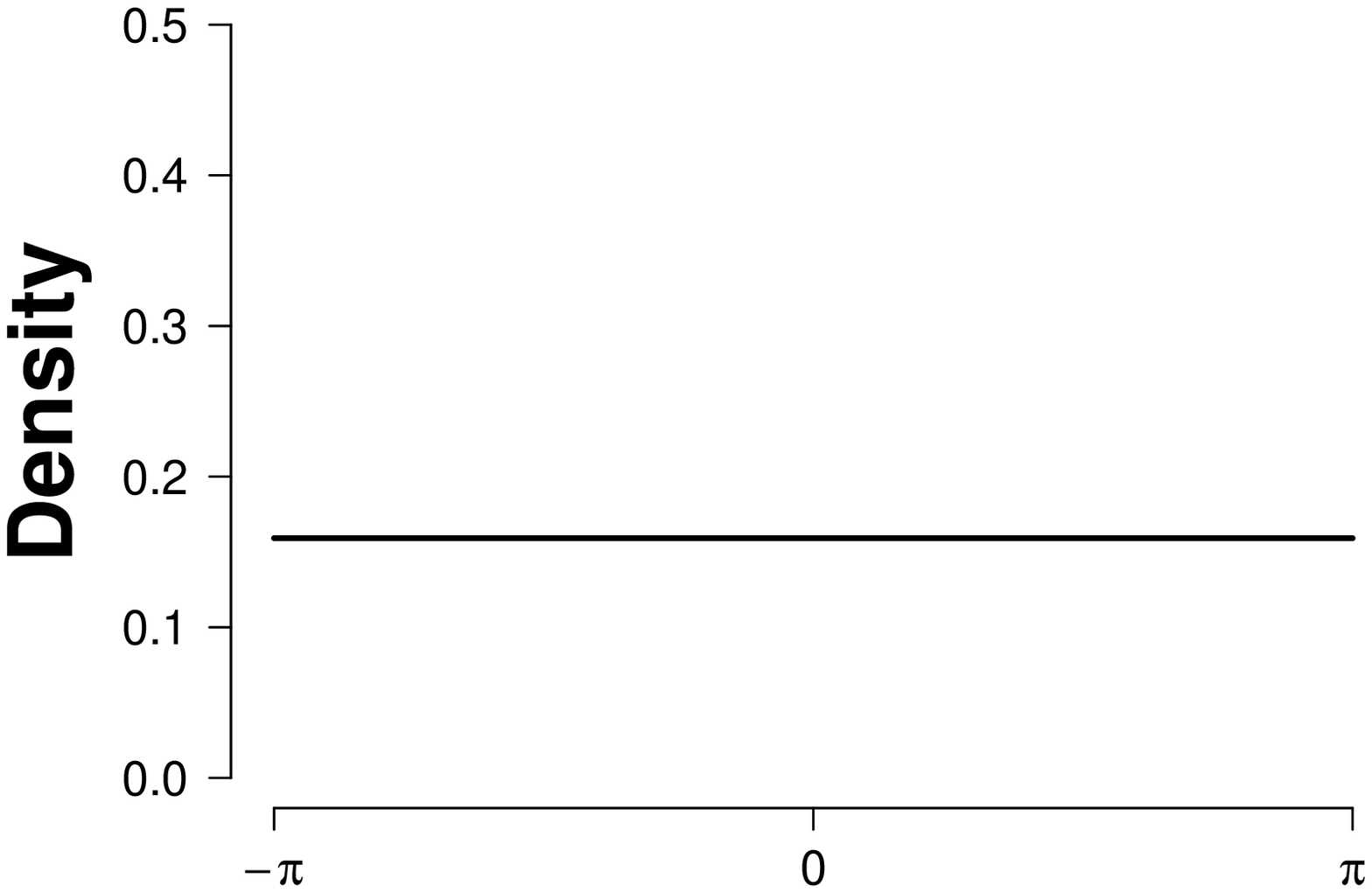}
    \end{minipage}  & %
    \scalebox{1.25}{\( \xrightarrow[n \,= \, 10]{y_{\obs} \, = \, 7} \)} & %
    \begin{minipage}{.375 \textwidth}
    \centering
    \scalebox{1.25}{Posterior \( \phi \) from \( \phi \sim U[-\pi, \pi] \)}
    \includegraphics[width=\linewidth]{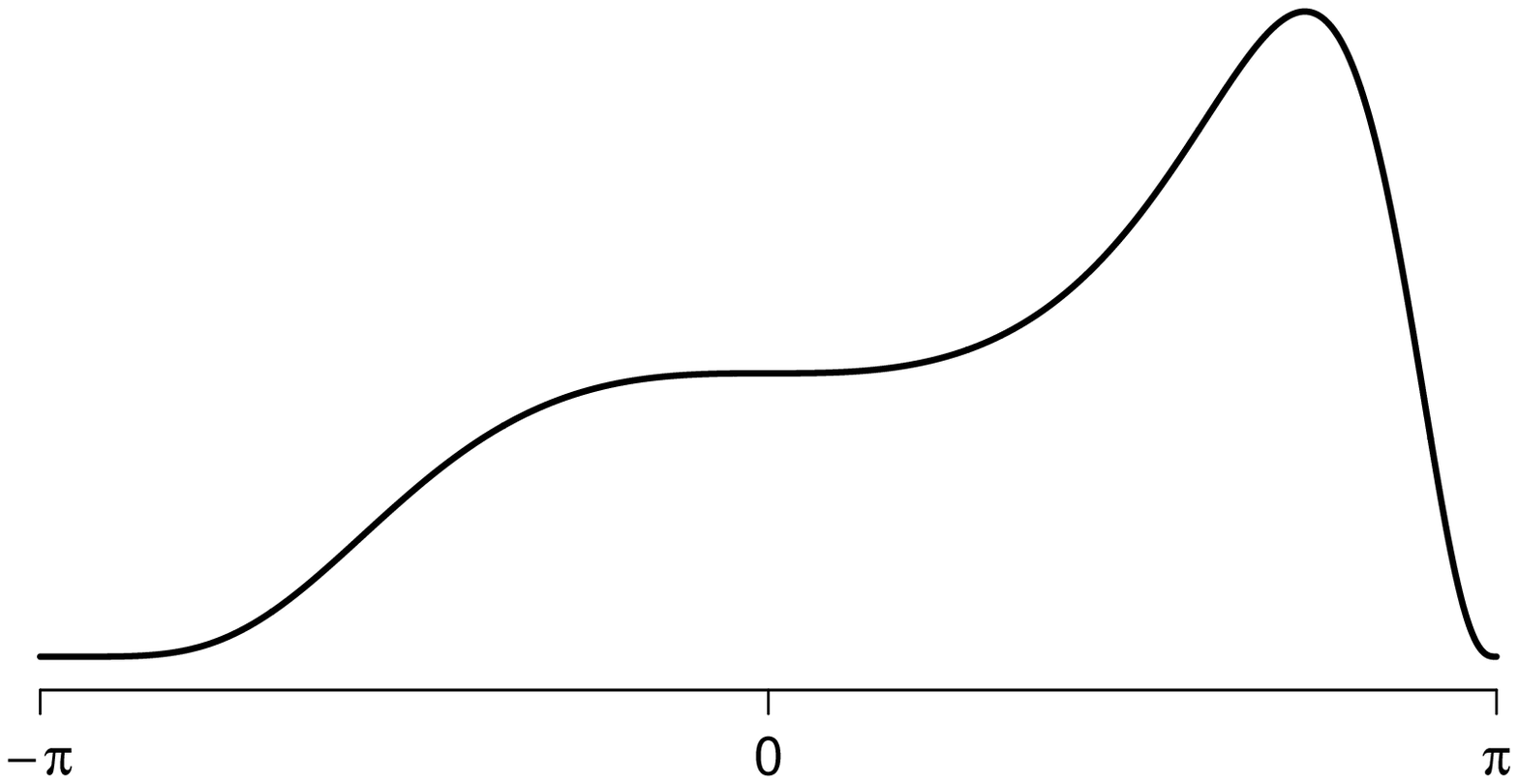}
    \end{minipage}  %
    \\
    \begin{minipage}{.375 \textwidth} %
    \begin{tabular}{cc}
    & \scalebox{1.25}{Angle \( \phi \)} \\
    & \scalebox{1.1}{\( h \)} \scalebox{1.5}{\( \dasheddownarrow \)} \\
    & \scalebox{1.25}{Prior \( \theta \) from \( \phi \sim U[-\pi, \pi] \)}
    \end{tabular}
    \end{minipage}  & %
      & %
    \begin{minipage}{.375 \textwidth}
    \begin{tabular}{c}
    \scalebox{1.25}{Angle \( \phi \)} \\
    \scalebox{1.1}{\( h \)} \scalebox{1.5}{\( \downarrow \)} \\
    \scalebox{1.25}{Posterior \( \theta \) from \( \phi \sim U[-\pi, \pi] \)}
    \end{tabular}
    \end{minipage}  %
    \\
    \begin{minipage}{.375 \textwidth}
    \centering
    \includegraphics[width= \linewidth]{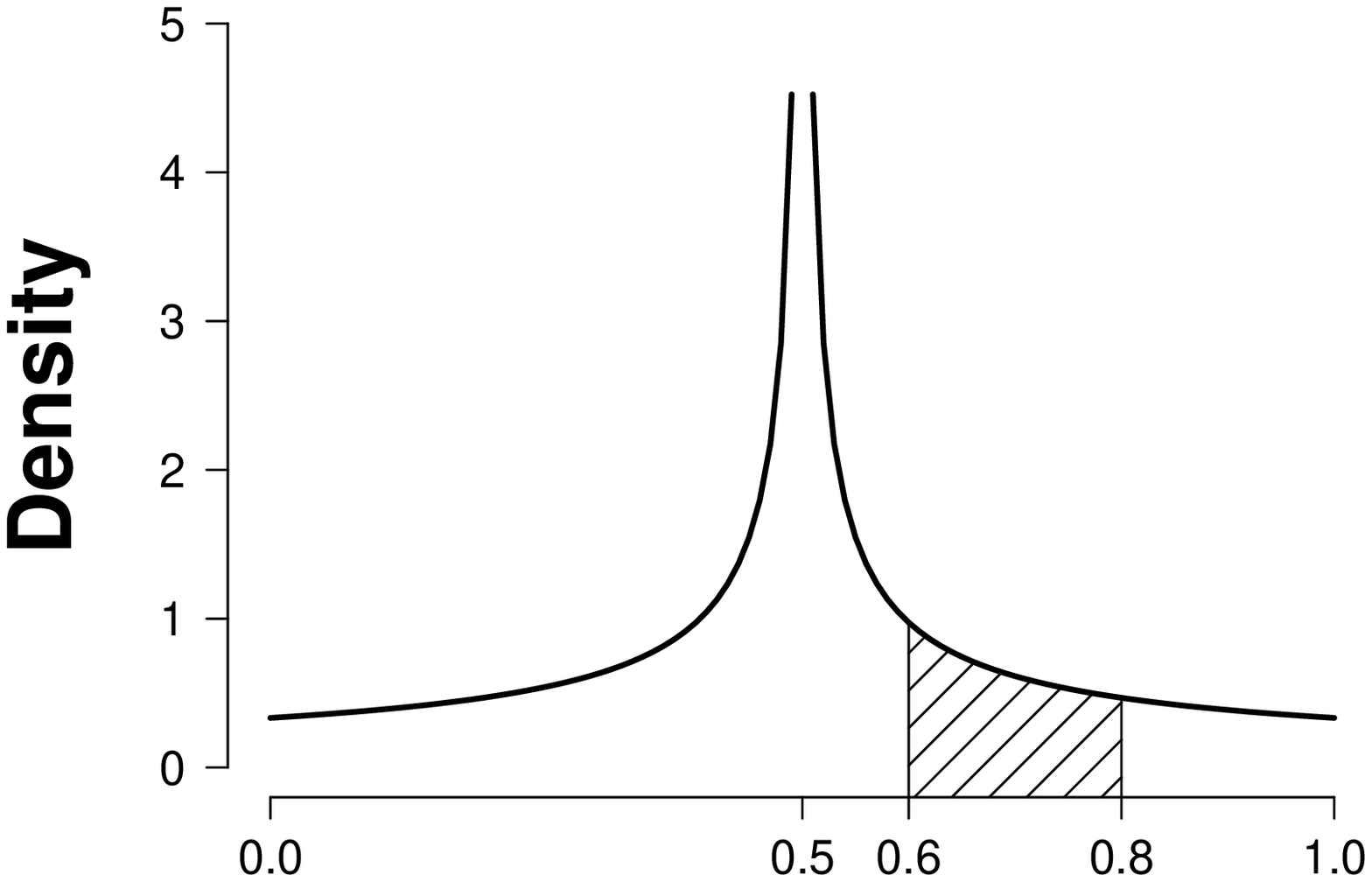} %
    \begin{tabular}{cc}
    & \scalebox{1.25}{Propensity \( \theta \)}
    \end{tabular}
    \end{minipage}  & %
    \scalebox{1.25}{\( \xdashrightarrow[n \,= \, 10]{y_{\obs} \, = \, 7} \)} & %
    \begin{minipage}{.375 \textwidth}
    \centering
    \includegraphics[width=\linewidth]{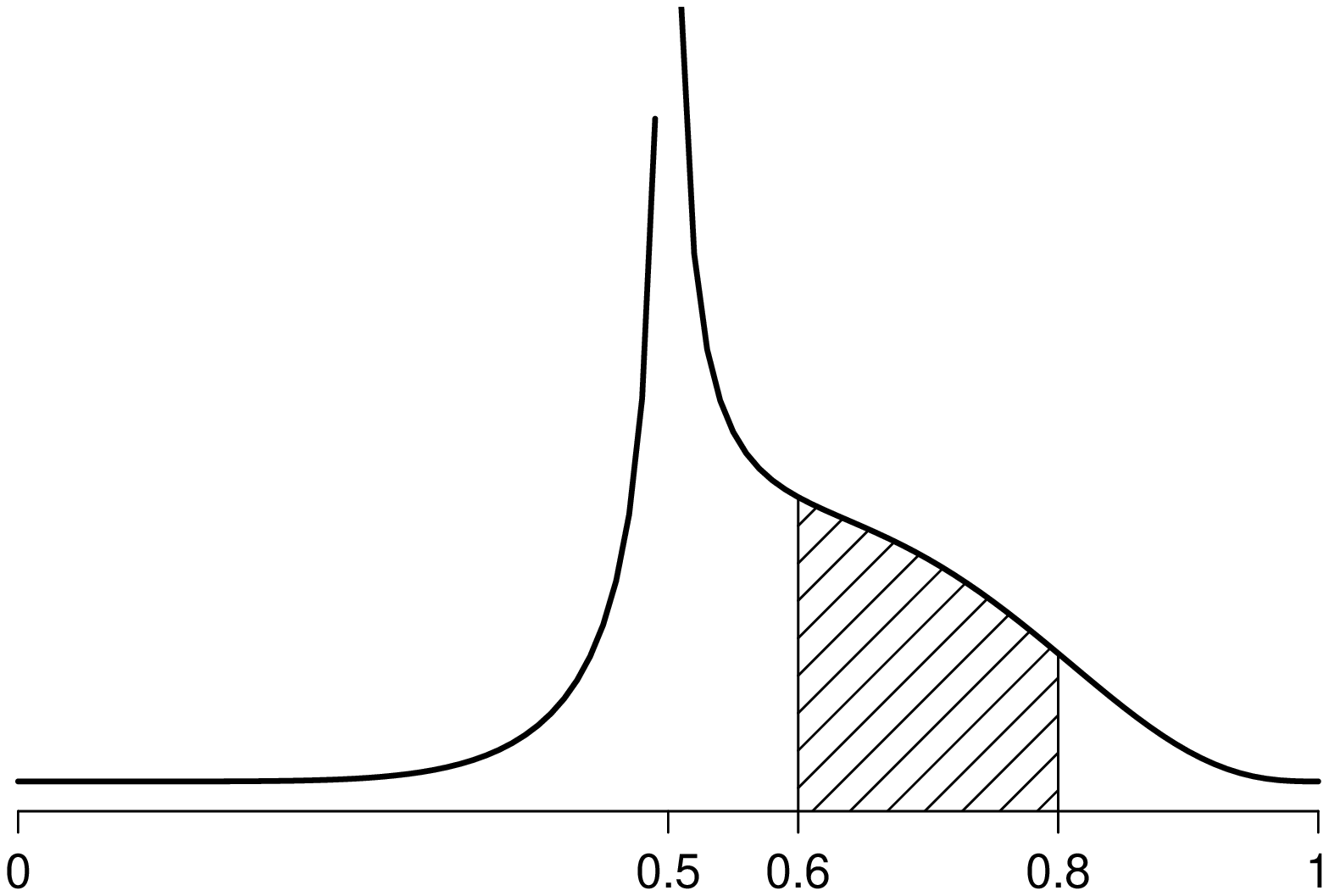}
    \scalebox{1.25}{Propensity \( \theta \)}
    \end{minipage}
\end{tabular}
\caption{Bayesian updating based on observations \( x^{n}_{\obs} \) with \( y_{\obs}=7 \) heads out of \( n=10 \) tosses when a uniform prior distribution is assigned to the the coin's angle \( \phi \). The uniform distribution is shown in the top-left panel. Bayes' theorem results in a posterior distribution for \( \phi \) that is shown in the top-right panel. This posterior \( \tilde{g}(\phi \, | \, x^{n}_{\obs}) \) is transformed into a posterior on \( \theta \) (bottom-right panel) using \( \theta=h(\phi) \). The same posterior on \( \theta \) is obtained if we proceed via an alternative route in which we first transform the uniform prior on \( \phi \) to the corresponding prior on \( \theta \) and then apply Bayes' theorem with the induced prior on \( \theta \). A comparison to the results from \protect \refFig{posteriorTheta} reveals that posterior inference differs notably depending on whether a uniform distribution is assigned to the angle \( \phi \) or to the propensity \( \theta \).}
\label{inductionFromUniformPhi}
\end{figure}
For the same data as before, the posterior calculated from Bayes' theorem is given in top-right panel of \refFig{inductionFromUniformPhi}. As the problem in terms of the angle \( \phi \) is equivalent to that of \( \theta=h(\phi) \) we can use the function \( h \) to translate the posterior in terms of \( \phi \) to a posterior on \( \theta \), see the bottom-right panel of \refFig{inductionFromUniformPhi}. This posterior on \( \theta \) is noticeably different from the posterior on \( \theta \) shown in Figure~\ref{posteriorTheta}.

Specifically, the uniform prior on \( \phi \) corresponds to the prior belief \( \tilde{P}(\theta^{*} \in J_{\theta})=0.13 \) of finding the true value \( \theta^{*} \) within the interval \( J_{\theta}=(0.6, 0.8) \). After observing \( x^{n}_{\obs} \) with \( y_{\obs}=7 \) out of \( n=10 \), this prior is updated to the posterior belief of \( \tilde{P}(\theta^{*} \in J_{\theta} \, | \, x^{n}_{\obs})=0.29 \),%
\footnote{The tilde makes explicit that the prior and posterior are derived from the uniform prior \( \tilde{g}(\phi) \) on \( \phi \).} %
see the shaded areas in \refFig{inductionFromUniformPhi}. Crucially, the earlier analysis that assigned a uniform prior to the propensity \( \theta \) yielded a posterior probability \( P(\theta^{*} \in J_{\theta} \, | \, x^{n}_{\obs})=0.54 \), which is markedly different from the current analysis that assigns a uniform prior to the angle \( \phi \).

The same posterior on \( \theta \) is obtained when the prior on \( \phi \) is first translated into a prior on \( \theta \) (bottom-left panel) and then updated to a posterior with Bayes' theorem. Regardless of the stage at which the transformation is applied, the resulting posterior on \( \theta \) differs substantially from the result plotted in the right panel of \refFig{posteriorTheta}. \( \hfill \diamond \)
\end{voorbeeld}

Thus, the uniform prior distribution is not a panacea for the quantification of prior ignorance, as the conclusions depend on how the problem is parameterized. In particular, a uniform prior on the coin's angle \( \tilde{g}(\phi)={1 / V_{\Phi}} \) yields a highly informative prior in terms of the coin's propensity \( \theta \). This lack of invariance caused Karl Pearson, Ronald Fisher and Jerzy Neyman to reject 19th century Bayesian statistics that was based on the uniform prior championed by Pierre-Simon Laplace. This rejection resulted in, what is now known as, frequentist statistics, see also \cite{hald2008history}, \cite{lehmann2011fisher}, and \cite{stigler1986history}. 

\subsection{A default prior by Jeffreys's rule}
Unlike the other fathers of modern statistical thoughts, Harold Jeffreys continued to study Bayesian statistics based on formal logic and his philosophical convictions of scientific inference (see, e.g., \citealp{aldrich2005statistical}; \citealp{etz2017haldane}; \citealp{jeffreys1961theory}; \citealp{ly2016harold, ly2016evaluation}; \citealp{robert2009harold}; \citealp{wrinch1919some, wrinch1921certain, wrinch1923certain}). Jeffreys concluded that the uniform prior is unsuitable as a default prior due to its dependence on the parameterization. As an alternative, \cite{jeffreys1946invariant} proposed the following prior based on Fisher information %
\begin{align}
\label{JeffreysRule}
g_{J}(\theta) = {1 \over V} \sqrt{I_{X}(\theta)} , \text{ where } V=\int_{\Theta} \sqrt{I_{X}(\theta)} \der \theta,
\end{align}%
which is known as the prior derived from Jeffreys's rule or the \emph{Jeffreys's prior} in short. The Jeffreys's prior is parameterization-invariant, which implies that it leads to the same posteriors regardless of how the model is represented. %
\begin{voorbeeld}[Jeffreys's prior]
\label{jeffreysPriorEx}
The Jeffreys's prior of the Bernoulli model in terms of \( \phi \) is %
\begin{align}
\label{jeffreysPriorPhi}
g_{J}(\phi) = {3 \phi^{2} \over V \sqrt{\pi^{6} - \phi^{6}}}, \text{ where } V=\pi,
\end{align}
which is plotted in the top-left panel of \refFig{jeffreysPrior}. %
\begin{figure}
\begin{tabular}{ccc}
    \begin{minipage}{.375 \textwidth}
    \centering
    \qquad \scalebox{1.25}{Jeffreys's prior on \( \phi \)}
    \includegraphics[width= \linewidth]{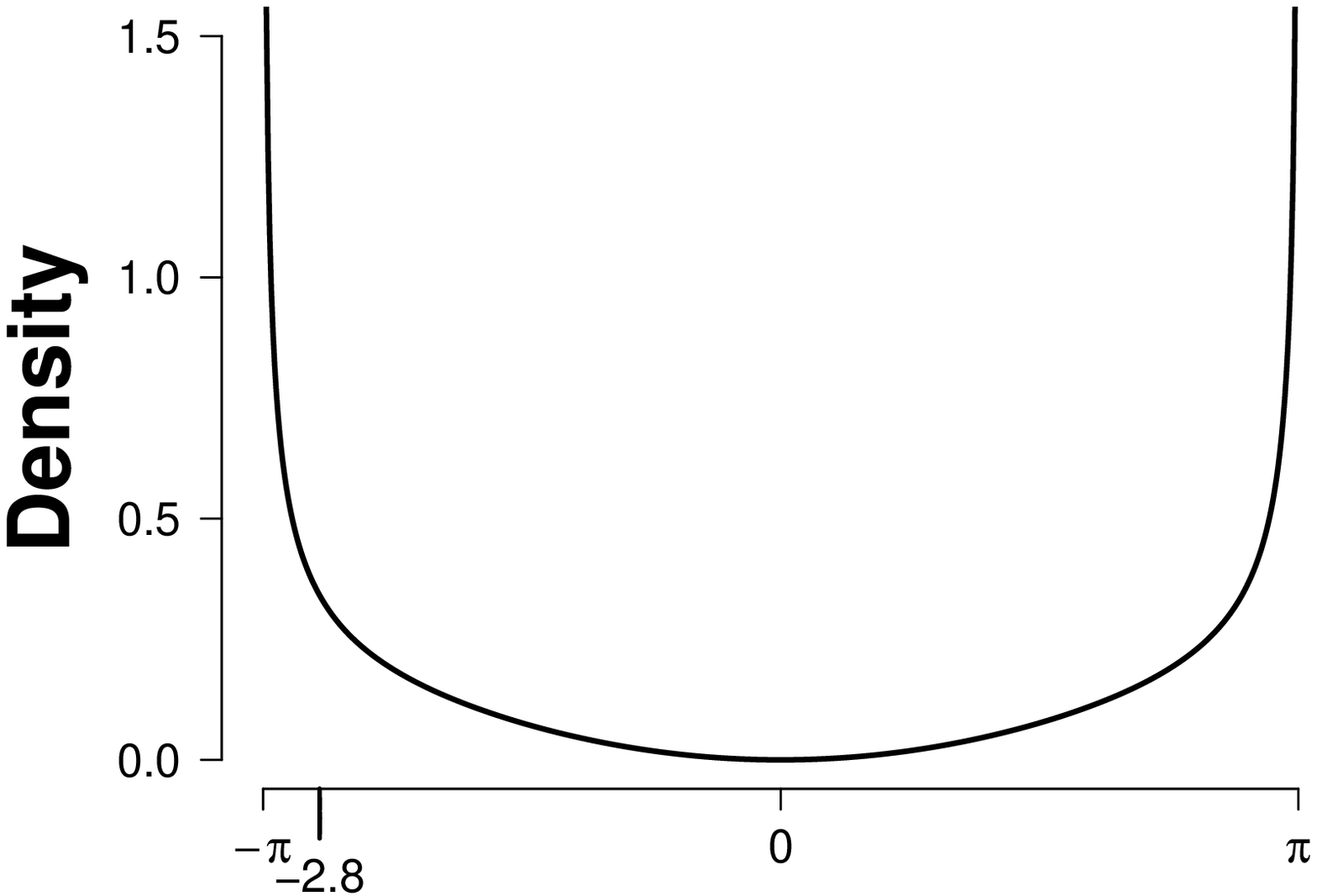}
    \end{minipage}  & %
    \scalebox{1.25}{\( \xrightarrow[n \,= \, 10]{y_{\obs} \, = \, 7} \)} & %
    \begin{minipage}{.375 \textwidth}
    \centering
    \scalebox{1.25}{Jeffreys's posterior on \( \phi \)}
    \includegraphics[width=\linewidth]{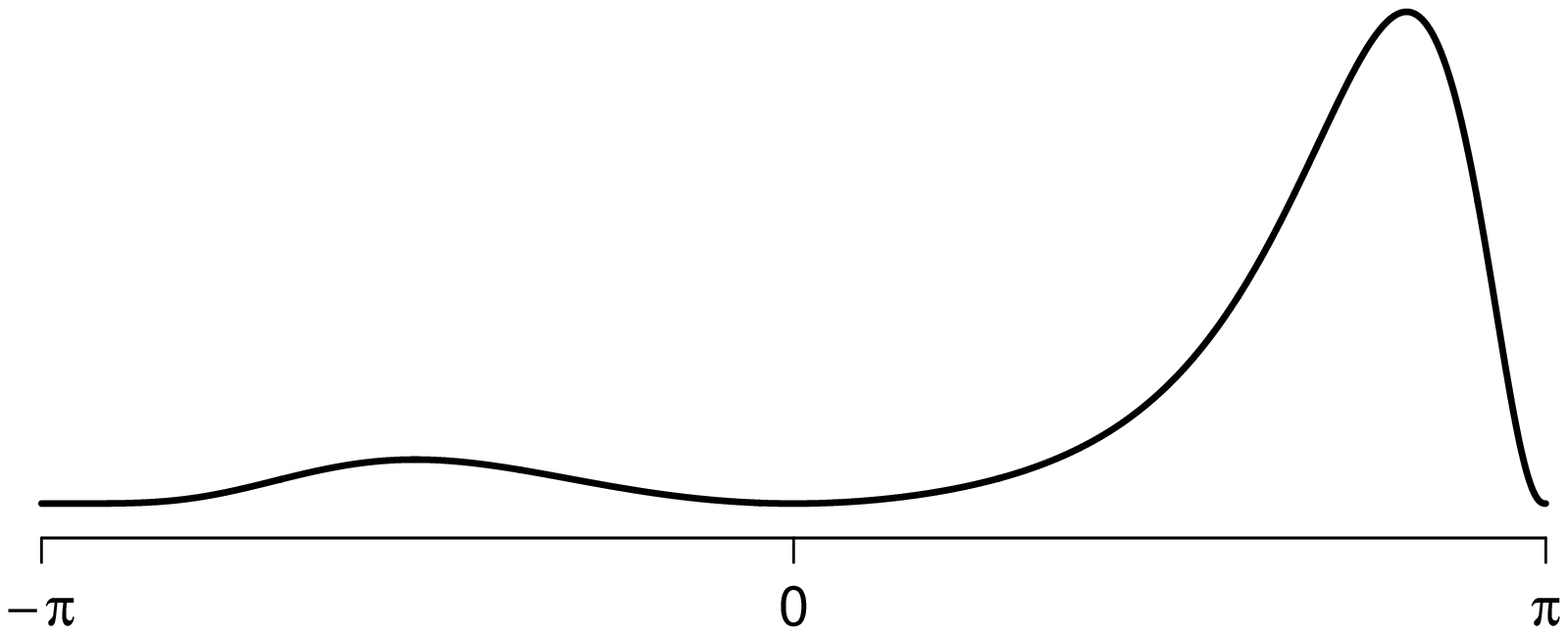}
    \end{minipage}  %
    \\
    \begin{minipage}{.375 \textwidth}
    \centering
    \begin{tabular}{cc}
    & \scalebox{1.25}{Angle \( \phi \)} \\
    & \, \scalebox{1.1}{\( h \)} \scalebox{1.5}{\( \downarrow \uparrow \)} \scalebox{1.1}{\( h^{-1} \)}\\
    & \scalebox{1.25}{Jeffreys's prior on \( \theta \)}
    \end{tabular}
    \end{minipage}  & %
      & %
    \begin{minipage}{.375 \textwidth}
    \begin{tabular}{c}
    \scalebox{1.25}{Angle \( \phi \)} \\
    \scalebox{1.1}{\( h \)} \scalebox{1.5}{\( \downarrow \uparrow \)} \scalebox{1.1}{\( h^{-1} \)} \\
    \scalebox{1.25}{Jeffreys's posterior on \( \theta \)}
    \end{tabular}
    \end{minipage}  %
    \\ %
    \begin{minipage}{.375 \textwidth}
    \centering
    \includegraphics[width= \linewidth]{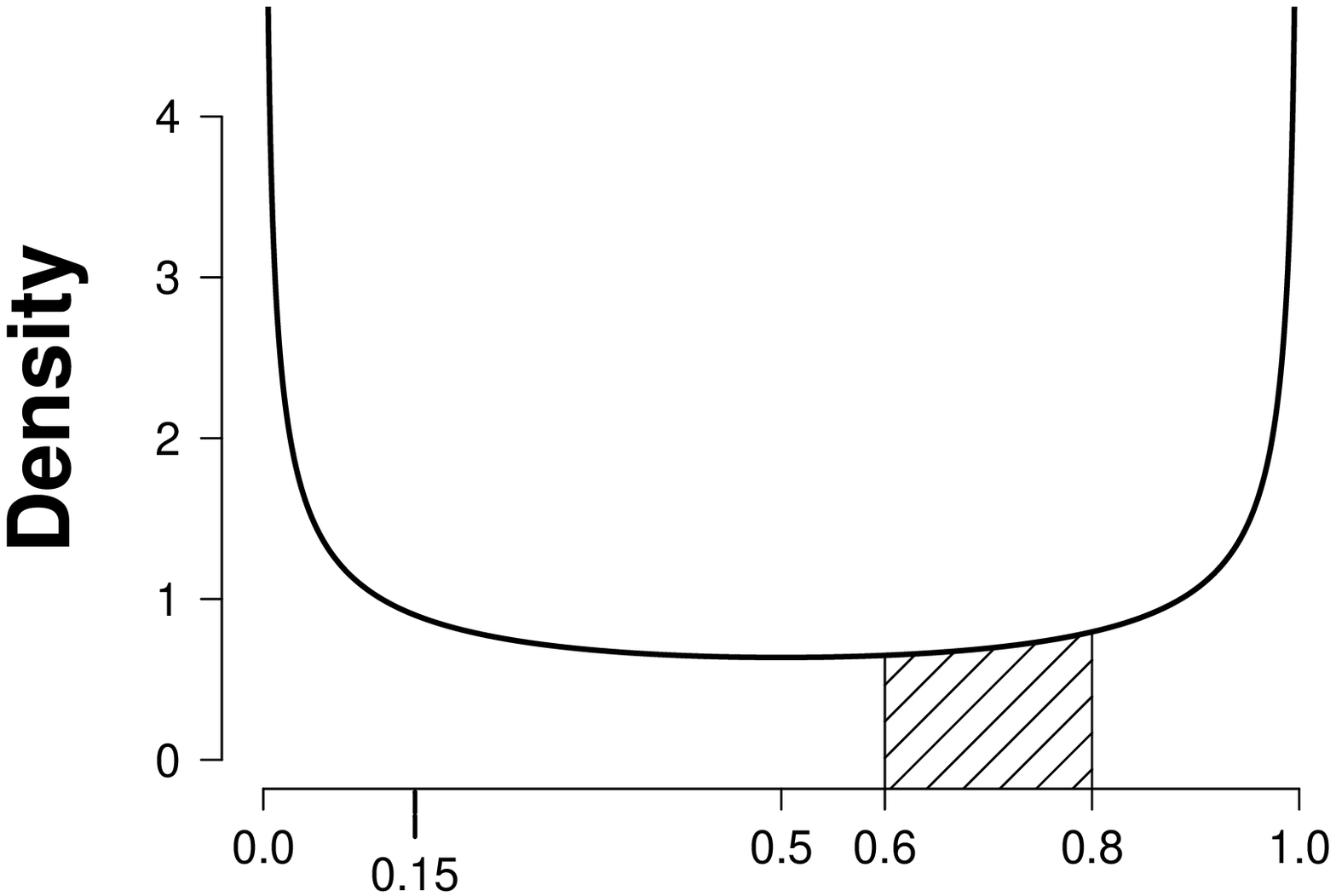}
    \begin{tabular}{cc}
    & \scalebox{1.25}{Propensity \( \theta \)}
    \end{tabular}
    \end{minipage}  & %
    \scalebox{1.25}{\( \xrightarrow[n \,= \, 10]{y_{\obs} \, = \, 7} \)} & %
    \begin{minipage}{.375 \textwidth}
    \centering
    \includegraphics[width=\linewidth]{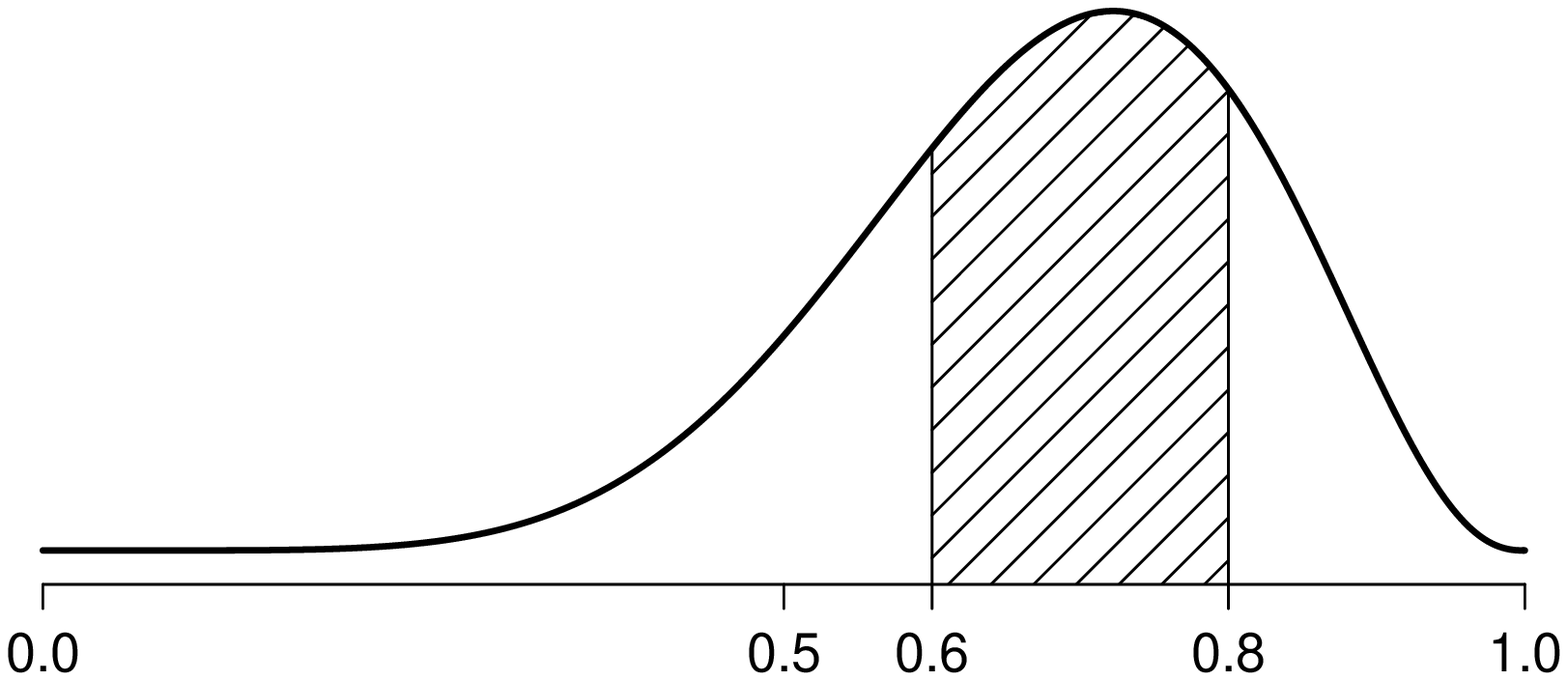}
    \scalebox{1.25}{Propensity \( \theta \)}
    \end{minipage}
\end{tabular}
\caption{For priors constructed through Jeffreys's rule it does not matter whether the problem is represented in terms of the angles \( \phi \) or its propensity \( \theta \). Thus, not only is the problem equivalent due to the transformations \( \theta=h(\phi) \) and its backwards transformation \( \phi=h^{-1}(\theta) \), the prior information is the same in both representations. This also holds for the posteriors.} 
\label{jeffreysPrior}
\end{figure}
The corresponding posterior is plotted in the top-right panel, which we transformed into a posterior in terms of \( \theta \) using the function \( \theta=h(\phi) \) shown in the bottom-right panel.%
\footnote{The subscript \( J \) makes explicit that the prior and posterior are based on the prior derived from Jeffreys's rule, i.e., \( g_{J}(\theta) \) on \( \theta \), or equivalently, \( g_{J}(\phi) \) on \( \phi \).} %

Similarly, we could have started with the Jeffreys's prior in terms of \( \theta \) instead, that is, %
\begin{align}
\label{jeffreysPriorTheta}
g_{J}(\theta) = {1 \over V \sqrt{\theta(1-\theta)}}, \text{ where } V=\pi.
\end{align}%
The Jeffreys's prior and posterior on \( \theta \) are plotted in the bottom-left and the bottom-right panel of \refFig{jeffreysPrior}, respectively. The Jeffreys's prior on \( \theta \) corresponds to the prior belief \( P_{J}(\theta^{*} \in J_{\theta})=0.14 \) of finding the true value \( \theta^{*} \) within the interval \( J_{\theta}=(0.6, 0.8) \). After observing \( x^{n}_{\obs} \) with \( y_{\obs}=7 \) out of \( n=10 \), this prior is updated to the posterior belief of \( P_{J}(\theta^{*} \in J_{\theta} \, | \, x^{n}_{\obs})=0.53 \), see the shaded areas in \refFig{jeffreysPrior}. The posterior is identical to the one obtained from the previously described updating procedure that starts with the Jeffreys's prior on \( \phi \) instead of on \( \theta \). \( \hfill \diamond \)
\end{voorbeeld}%
This example shows that the Jeffreys's prior leads to the same posterior knowledge regardless of how we as researcher represent the problem. Hence, the same conclusions about \( \theta \) are drawn regardless of whether we (1) use Jeffreys's rule to construct a prior on \( \theta \) and update with the observed data, or (2) use Jeffreys's rule to construct a prior on \( \phi \), update to a posterior distribution on \( \phi \), which is then transformed to a posterior on \( \theta \). 

\subsection{Geometrical properties of Fisher information}
In the remainder of this section we make intuitive that the Jeffreys's prior is in fact uniform in the model space. We elaborate on what is meant by model space and how this can be viewed geometrically. This geometric approach illustrates (1) the role of Fisher information in the definition of the Jeffreys's prior, (2) the interpretation of the shaded area, and (3) why the normalizing constant is \( V=\pi \), regardless of the chosen parameterization. 

\subsubsection{The model space \( \Mc \)}
\label{secCompleteSet2d}
Before we describe the geometry of statistical models, recall that at a pmf can be thought of as a data generating device of \( X \), as the pmf specifies the chances with which \( X \) takes on the potential outcomes \( 0 \) and \( 1 \). Each such pmf has to fulfil two conditions: (i) the chances have to be non-negative, that is, \( 0 \leq p(x)= P(X=x) \) for every possible outcome \( x \) of \( X \), and (ii) to explicitly convey that there are \( w=2 \) outcomes, and none more, the chances have to sum to one, that is, \( p(0)+p(1)=1 \). We call the largest set of functions that adhere to conditions (i) and (ii) the complete set of pmfs \( \Pc \). 

\begin{figure}[h]
    \centering
\begin{tabular}{cc}
    \begin{minipage}{.4 \textwidth}
    \centering
    \includegraphics[width= \linewidth]{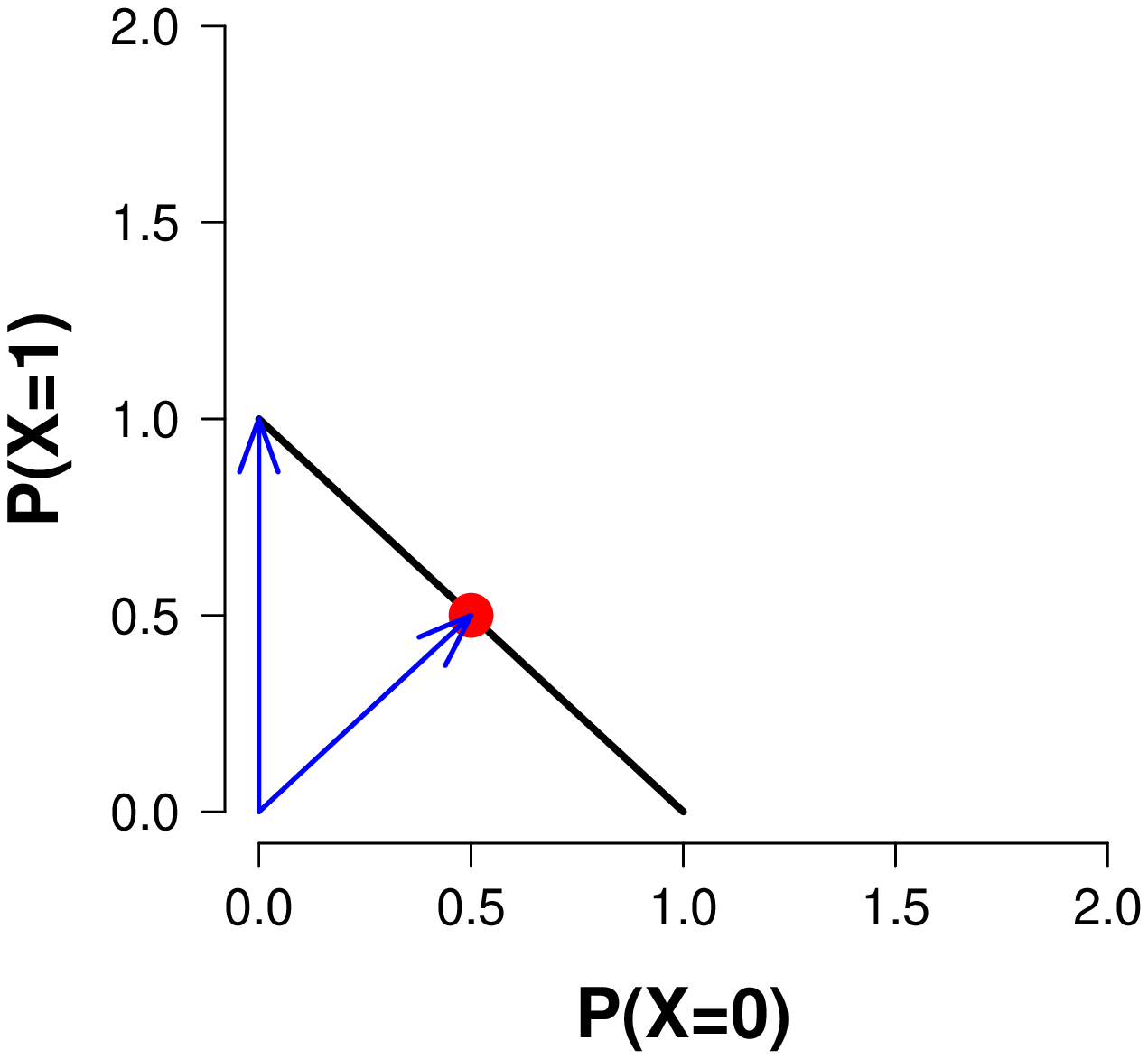}
    \end{minipage}  & %
    \begin{minipage}{.4 \textwidth}
    \centering
    \includegraphics[width=\linewidth]{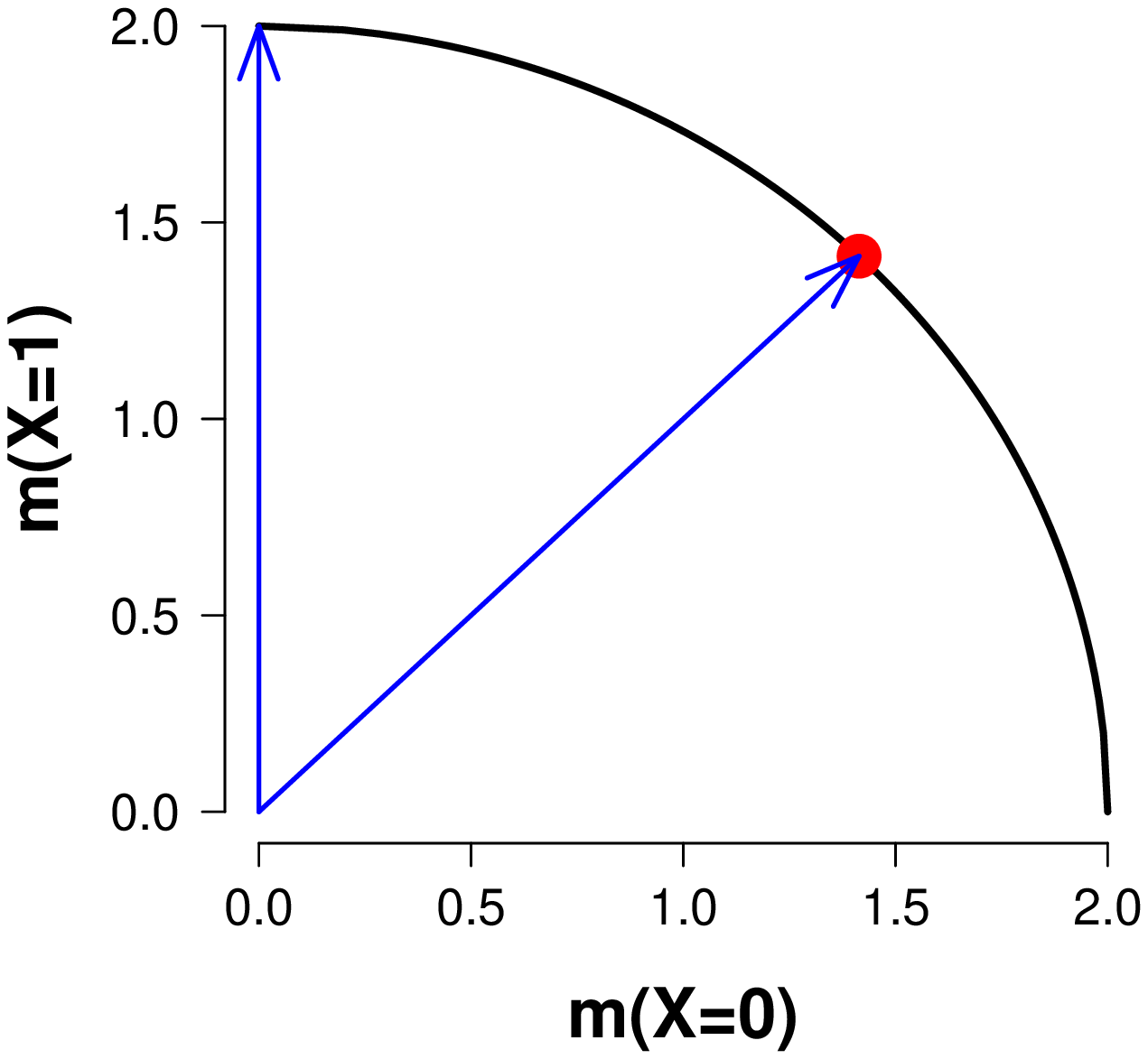}
    \end{minipage}  %
\end{tabular}  
\caption{The true pmf of \( X \) with the two outcomes \( \{0, 1 \} \) has to lie on the line (left panel) or more naturally on the positive part of the circle (right panel). The dot represents the pmf \( p_{e}(X) \).}
\label{2dGeoA}
\end{figure}

As any pmf from \( \Pc \) defines \( w=2 \) chances, we can represent such a pmf as a vector in \( w \) dimensions. To simplify notation, we write \( p(X) \) for all \( w \) chances simultaneously, hence, \( p(X) \) is the vector \( p(X)=[p(0), p(1)] \) when \( w=2 \). %
The two chances with which a pmf \( p(X) \) generates outcomes of \( X \) can be simultaneously represented in the plane with \( p(0)=P(X=0) \) on the horizontal axis and \( p(1)=P(X=1) \) on the vertical axis. In the most extreme case, we have the pmf \( p(X)=[1, 0] \) or \( p(X)=[0,1] \). These two extremes are linked by a straight line in the left panel of \refFig{2dGeoA}. 
\emph{Any} pmf --and the true pmf \( p^{*}(X) \) of \( X \) in particular-- can be uniquely identified with a vector on the line and vice versa. For instance, the pmf \( p_{e}(X)=[1/2, 1/2] \) (i.e., the two outcomes are generated with the same chance) is depicted as the dot on the line. 

This vector representation allows us to associate to each pmf of \( X \) a norm, that is, a length. Our intuitive notion of length is based on the \emph{Euclidean norm} and entails taking the root of the sums of squares. For instance, we can associate to the pmf \( p_{e}(X) \) the length \( \| p_{e}(X) \|_{2} = \sqrt{ (1/2)^{2} + (1/2)^{2}} = 1/ \sqrt{2} \approx 0.71 \). On the other hand, the length of the pmf that states that \( X=1 \) is generated with 100\% chance has length one. Note that by eye, we conclude that \( p_{e}(X) \), the arrow pointing to the dot in the left panel in \refFig{2dGeoA} is indeed much shorter than the arrow pointing to extreme pmf \( p(X)=[0,1] \). 

This mismatch in lengths can be avoided when we represent each pmf \( p(X) \) by two times its square root instead \citep{kass1989geometry}, that is, by \( m(X)=2 \sqrt{p(X)}=[2 \sqrt{p(0)}, 2 \sqrt{p(1)}] \).%
\footnote{The factor two is used to avoid a scaling of a quarter, though, its precise value is not essential for the ideas conveyed here. To simplify matters, we also call \( m(X) \) a pmf.} %
A pmf that is identified as the vector \( m(X) \) is now two units away from the origin, that is, \( \| m(X) \|_{2} = \sqrt{ m(0)^{2} + m(1)^{2}} = \sqrt{4(p(0)+p(1))}=2 \). For instance, the pmf \( p_{e}(X) \) is now represented as \( m_{e}(X) \approx [1.41, 1.41] \). The model space \( \Mc \) is collection of all transformed pmfs and represented as the surface of (the positive part of) a circle, see the right panel of \refFig{2dGeoA}.%
\footnote{Hence, the model space \( \Mc \) is the collection of all functions on \( \Xc \) such that (i) \( m(x) \geq 0  \) for every outcome \( x \) of \( X \), and (ii) \( \sqrt{m(0)^{2} + m(1)^{2}}=2 \). This vector representation of all the pmfs on \( X \) has the advantage that it also induces an inner product, which allows one to project one vector onto another, see \citet[p. 4]{rudin1991functional}, \citet[p. 94]{van1998asymptotic} and Appendix~{\ref{appendixReg}}.} %
By representing the set of all possible pmfs of \( X \) as vectors \( m(X)=2 \sqrt{p(X)} \) that reside on the sphere \( \Mc \), we adopted our intuitive notion of distance. As a result, we can now, by simply looking at the figures, clarify that a uniform prior on the parameter space may lead to a very informative prior in the model space \( \Mc \). 

\subsubsection{Uniform on the parameter space versus uniform on the model space}
\begin{figure}[h]
   \centering
\begin{tabular}{cc}
    \begin{minipage}{.4 \textwidth}
    \centering
    \includegraphics[width= \linewidth]{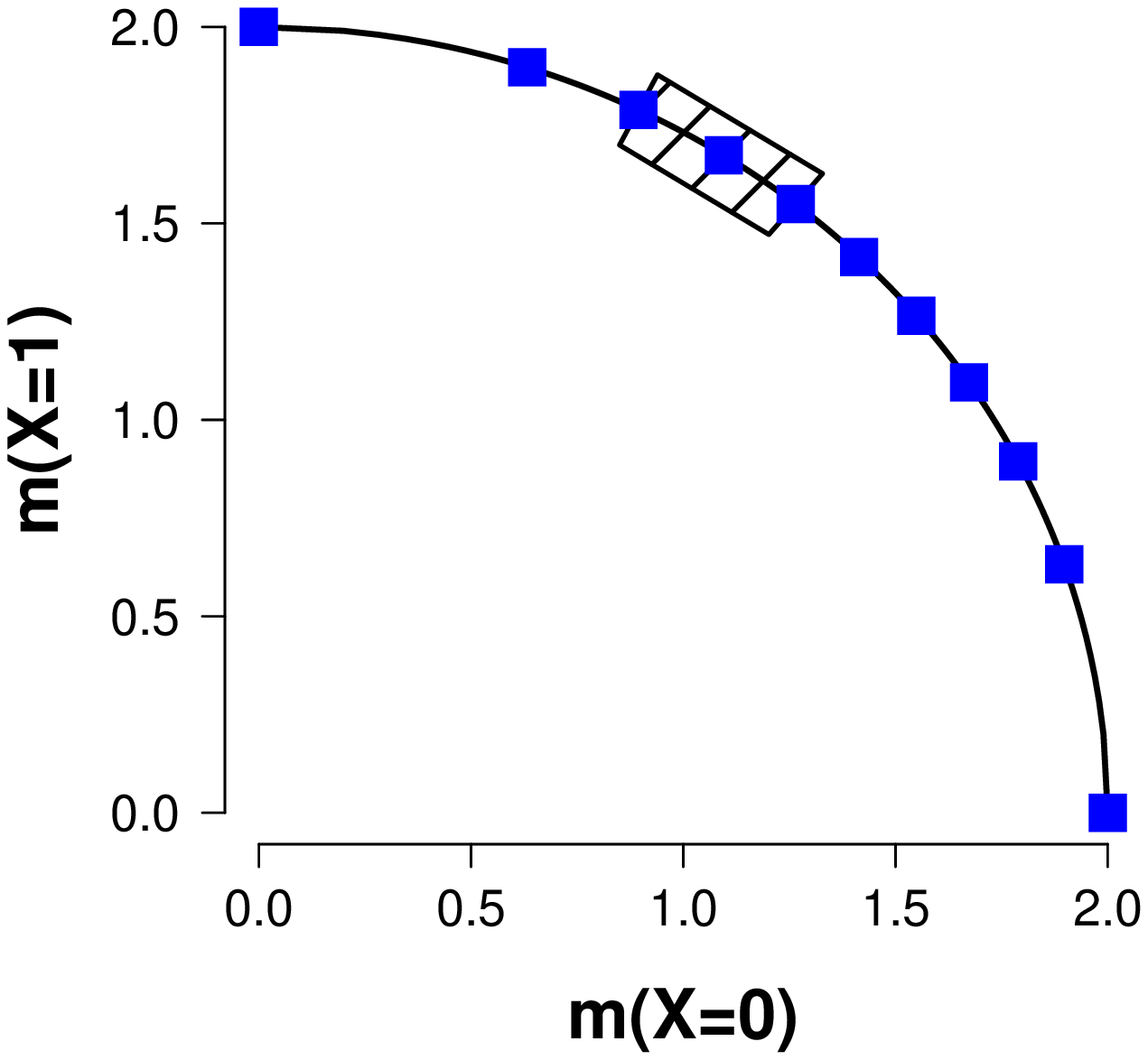}
    \end{minipage}  & %
    \begin{minipage}{.4 \textwidth}
    \centering
    \includegraphics[width= \linewidth]{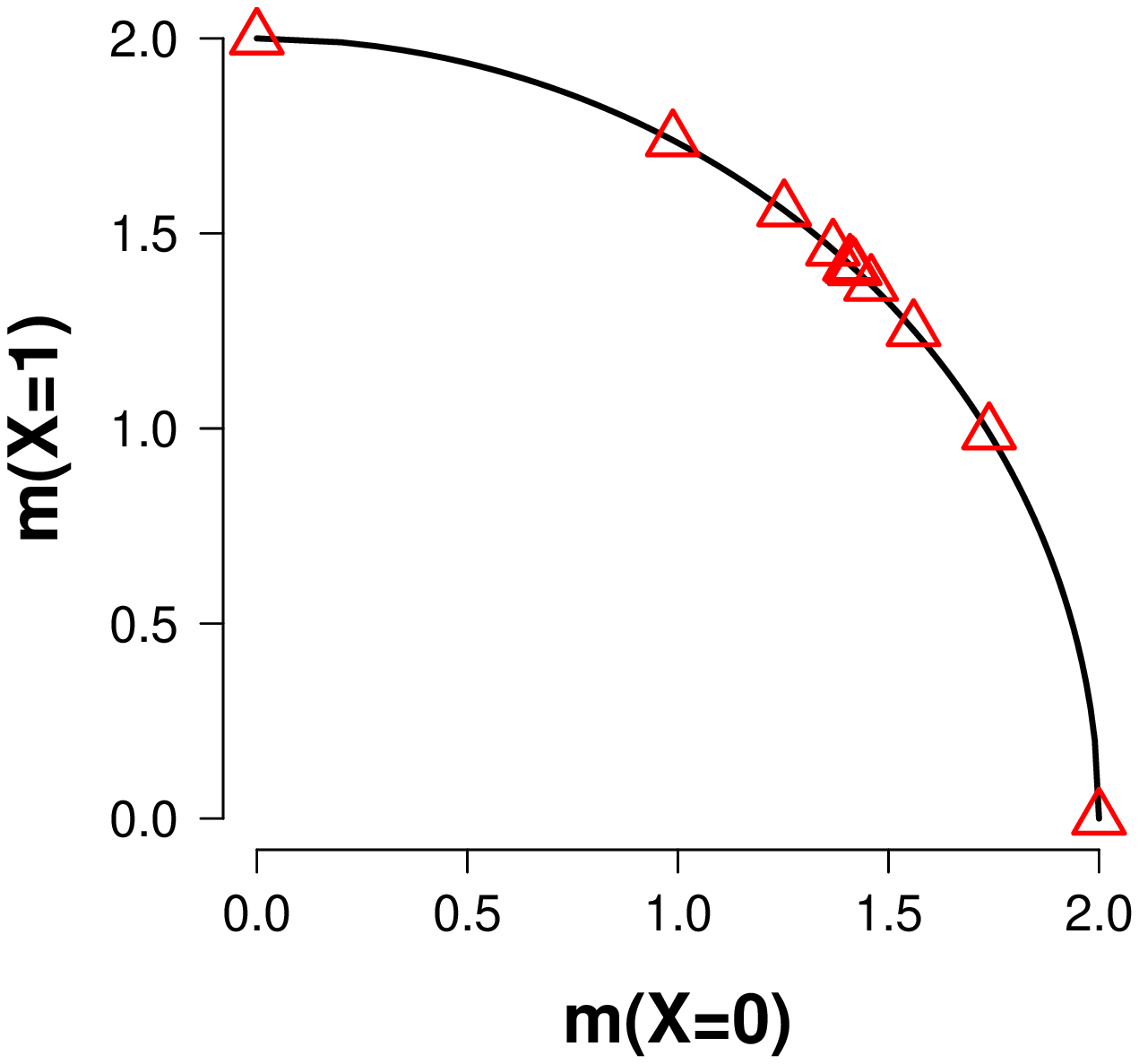}
    \end{minipage}  %
\end{tabular}    
\caption{The parameterization in terms of propensity \( \theta \) (left panel) and angle \( \phi \) (right panel) differ from each other substantially, and from a uniform prior in the model space. Left panel: The eleven squares (starting from the right bottom going anti-clockwise) represent pmfs that correspond to \( \theta = 0.0, 0.1, 0.2, \ldots, 0.9, 1.0 \). The shaded area corresponds to the shaded area in the bottom-left panel of \refFig{jeffreysPrior} and accounts for 14\% of the model's length. Right panel: Similarly, the eleven triangles (starting from the right bottom going anti-clockwise) represent pmfs that correspond to \( \phi = -1.0 \pi, -0.8 \pi, -\ldots 0.8 \pi, 1.0 \pi \).}
\label{2dGeoB}
\end{figure}

As \( \Mc \) represents the largest set of pmfs, any model defines a subset of \( \Mc \). Recall that the function \( f(x \, | \, \theta) \) represents how we believe a parameter \( \theta \) is functionally related to an outcome \( x \) of \( X \). For each \( \theta \) this parameterization yelds a pmf \( p_{\theta}(X) \) and, thus, also \( m_{\theta}(X) = 2 \sqrt{p_{\theta}(X)} \). We denote the resulting set of vectors \( m_{\theta}(X) \) so created by \( \Mc_{\Theta} \). For instance, the Bernoulli model \( f(x \, | \, \theta)=\theta^{x}(1-\theta)^{1-x} \) consists of pmfs given by \( p_{\theta}(X)=[f(0\, | \, \theta), f(1\, | \, \theta)]=[1-\theta, \theta] \), which we represent as the vectors \( m_{\theta}(X)=[2 \sqrt{1-\theta} , 2 \sqrt{\theta}] \). Doing this for every \( \theta \) in the parameter space \( \Theta \) yields the candidate set of pmfs \( \Mc_{\Theta} \). In this case, we obtain a saturated model, since \( \Mc_{\Theta}=\Mc \), see the left panel in \refFig{2dGeoB}, where the right most square on the curve corresponds to \( m_{0}(X)=[2, 0] \). By following the curve in an anti-clockwise manner we encounter squares that represent the pmfs \( m_{\theta}(X) \) corresponding to \( \theta = 0.1, 0.2, \ldots, 1.0 \) respectively. %
In the right panel of \refFig{2dGeoB} the same procedure is repeated, but this time in terms of \( \phi \) at \( \phi=-1.0 \pi, -0.8 \pi, \ldots, 1.0 \pi \). Indeed, filling in the gaps shows that the Bernoulli model in terms of \( \theta \) and \( \phi \) fully overlap with the largest set of possible pmfs, thus, \( \Mc_{\Theta}=\Mc=\Mc_{\Phi} \). \refFig{2dGeoB} makes precise what is meant when we say that the models \( \Mc_{\Theta} \) and \( \Mc_{\Phi} \) are equivalent; the two models define the same candidate set of pmfs that we believe to be viable data generating devices for \( X \).

However, \( \theta \) and \( \phi \) represent \( \Mc \) in a substantially different manner. As the representation \( m(X)=2 \sqrt{p(X)} \) respects our natural notion of distance, we conclude, by eye, that a uniform division of \( \theta \)s with distance, say, \( \der \theta=0.1 \) does not lead to a uniform partition of the model. More extremely, a uniform division of \( \phi \) with distance \( \der \phi = 0.2 \pi \) (10\% of the length of the parameter space) also does not lead to a uniform partition of the model. In particular, even though the intervals \( (-\pi, -0.8 \pi) \) and \( (-0.2 \pi, 0) \) are of equal length in the parameter space \( \Phi \), they do not have an equal displacement in the model \( \Mc_{\Phi} \). In effect, the right panel of \refFig{2dGeoB} shows that the 10\% probability that the uniform prior on \( \phi \) assigns to \( \phi^{*} \in (-\pi, -0.8 \pi) \) in parameter space is redistributed over a larger arc length of the model \( \Mc_{\Phi} \) compared to the 10\% assigned to \( \phi^{*} \in (-0.2 \pi, 0) \). Thus, a uniform distribution on \( \phi \) favors the pmfs \( m_{\phi}(X) \) with \( \phi \) close to zero. Note that this effect is cancelled by the Jeffreys's prior, as it puts more mass on the end points compared to \( \phi=0 \), see the top-left panel of \refFig{jeffreysPrior}. Similarly, the left panel of \refFig{2dGeoB} shows that the uniform prior \( g(\theta) \) also fails to yield an equiprobable assessment of the pmfs in model space. Again, the Jeffreys's prior in terms of \( \theta \) compensates for the fact that the interval \( (0, 0.1) \) as compared to \( (0.5, 0.6) \) in \( \Theta \) is more spread out in model space. However, it does so less severely compared to the Jeffreys's prior on \( \phi \). To illustrate, we added additional tick marks on the horizontal axis of the priors in the left panels of \refFig{jeffreysPrior}. The tick mark at \( \phi=-2.8 \) and \( \theta=0.15 \) both indicate the 25\% quantiles of their respective Jeffreys's priors. Hence, the Jeffreys's prior allocates more mass to the boundaries of \( \phi \) than to the boundaries of \( \theta \) to compensate for the difference in geometry, see \refFig{2dGeoB}. More generally, the Jeffreys's prior uses Fisher information to convert the geometry of the model to the parameter space. 

Note that because the Jeffreys's prior is specified using the Fisher information, it takes the functional relationship \( f(x \, | \, \theta) \) into account. The functional relationship makes precise how the parameter is linked to the data and, thus, gives meaning and context to the parameter. On the other hand, a prior on \( \phi \) specified without taking the functional relationship \( f(x \, | \, \phi) \) into account is a prior that neglects the context of the problem. For instance, the right panel of \refFig{2dGeoB} shows that this neglect with a uniform prior on \( \phi \) results in having the geometry of \( \Phi=(-\pi, \pi) \) forced onto the model \( \Mc_{\Phi} \). 

\subsection{Uniform prior on the model}
\refFig{2dGeoB} shows that neither a uniform prior on \( \theta \), nor a uniform prior on \( \phi \) yields a uniform prior on the model. Alternatively, we can begin with a uniform prior on the model \( \Mc \) and convert this into priors on the parameter spaces \( \Theta \) and \( \Phi \). This uniform prior on the model translated to the parameters is exactly the Jeffreys's prior. 

Recall that a prior on a space \( S \) is uniform, if it has the following two defining features: (i) the prior is proportional to one, and (ii) a normalizing constant given by \( V_{S} = \int_{S} 1 \der s \) that equals the length, more generally, volume of \( S \). For instance, a replacement of \( s \) by \( \phi \) and \( S \) by \( \Phi=(- \pi, \pi) \) yields the uniform prior on the angles with the normalizing constant \( V_{\Phi}= \int_{\Phi} 1 \der \phi = 2 \pi \). Similarly, a replacement of \( s \) by the pmf \( m_{\theta}(X) \) and \( S \) by the function space \( \Mc_{\Theta} \) yields a uniform prior on the model \( \Mc_{\Theta} \). The normalizing constant then becomes a daunting looking integral in terms of displacements \( \der m_{\theta}(X) \) between functions in model space \( \Mc_{\Theta} \). Fortunately, it can be shown, see Appendix~{\ref{appendixBayes}}, that \( V \) simplifies to %
\begin{align}
 V=\int_{\Mc_{\Theta}} 1 \der m_{\theta}(X) = \int_{\Theta} \sqrt{I_{X}(\theta)} \der \theta.
\end{align}
Thus, \( V \) can be computed in terms of \( \theta \) by multiplying the distances \( \der \theta \) in \( \Theta \) by the root of the Fisher information. Heuristically, this means that the root of the Fisher information translates displacements \( \der m_{\theta}(X) \) in the model \( \Mc_{\Theta} \) to distances \( \sqrt{I_{X}(\theta)} \der \theta \) in the parameter space \( \Theta \). %

Recall from \refEx{jeffreysPriorEx} that regardless of the parameterization, the normalizing constant of the Jeffreys's prior was \( \pi \). To verify that this is indeed the length of the model, we use the fact that the circumference of a quarter circle with radius \( r=2 \) can also be calculated as \( V=(2 \pi r)/4=\pi \). 

Given that the Jeffreys's prior corresponds to a uniform prior on the model, we deduce that the shaded area in the bottom-left panel of \refFig{jeffreysPrior} with \( P_{J}(\theta^{*} \in J_{\theta})=0.14 \), implies that the model interval \( J_{m} =\Big (m_{0.6}(X), m_{0.8}(X) \Big ) \), the shaded area in the left panel of \refFig{2dGeoB}, accounts for 14\% of the model's length. After updating the Jeffreys's prior with the observations \( x^{n}_{\obs} \) consisting of \( y_{\obs}=7 \) out of \( n=10 \) the probability of finding the true data generating pmf \( m^{*}(X) \) in this interval of pmfs \( J_{m} \) is increased to 53\%. 

In conclusion, we verified that the Jeffreys's prior is a prior that leads to the same conclusion regardless of how we parameterize the problem. This parameterization-invariance property is a direct result of shifting our focus from finding the true parameter value within the parameter space to the proper formulation of the estimation problem --as discovering the true data generating pmf \( m_{\theta^{*}}(X)=2 \sqrt{p_{\theta^{*}}(X)} \) in \( \Mc_{\Theta} \) and by expressing our prior ignorance as a uniform prior on the model \( \Mc_{\Theta} \). 

\section{The Role of Fisher Information in Minimum Description Length}
\label{fiInMdl}
In this section we graphically show how Fisher information is used as a measure of model complexity and its role in model selection within the minimum description length framework (MDL; \citealp{deRooij2010luckiness}; \citealp{grunwald2005advances}; \citealp{grunwald2007minimum}; \citealp{myung2000model}; \citealp{myung2006model}; \citealp{pitt2002toward}).

The primary aim of a model selection procedure is to select a single model from a set of competing models, say, models \( \mathcal{M}_{1} \) and \( \mathcal{M}_{2} \), that best suits the observed data \( x^{n}_{\obs} \). Many model selection procedures have been proposed in the literature, but the most popular methods are those based on penalized maximum likelihood criteria, such as the Akaike information criterion (AIC; \citealp{akaike1974new}; \citealp{burnham2002model}), the Bayesian information criterion (BIC; \citealp{raftery1995bayesian}; \citealp{schwarz1978estimating}), and the Fisher information approximation (FIA; \citealp{grunwald2007minimum}; \citealp{rissanen1996fisher}). These criteria are defined as follows
\begin{align}
\label{AIC}
\text{AIC} = & - 2 \log f_{j} \big (x^{n}_{\obs} \, | \, \hat{\theta}_{j}(x^{n}_{\obs}) \big ) & & + \quad 2 d_{j}, \\
\label{BIC}
\text{BIC} = & - 2 \log f_{j} \big (x^{n}_{\obs} \, | \, \hat{\theta}_{j}(x^{n}_{\obs}) \big ) & & + \quad d_{j} \log(n), \\
\label{FIA}
 \text{FIA} = & \underbrace{- \log f_{j} \big (x^{n}_{\obs} \, | \, \hat{\theta}_{j}(x^{n}_{\obs}) \big )}_{\text{Goodness-of-fit}} & & + \underbrace{\frac{d_{j}}{2} \log \frac{n}{2 \pi}}_{\text{Dimensionality}} + \underbrace{\log \left ( \int_{\Theta} \sqrt{\text{det } I_{\Mc_{j}}(\theta_{j})} \, \der \theta_{j} \right )}_{\text{Geometric complexity}},
\end{align}
%
%
%
where \( n \) denotes the sample size, \( d_{j} \) the number of free parameters, \( \hat{\theta}_{j} \) the MLE, \( I_{\Mc_{j}}(\theta_{j}) \) the unit Fisher information, and \( f_{j} \) the functional relationship between the potential outcome \( x^{n} \) and the parameters \( \theta_{j} \) within model \( \Mc_{j} \).%
\footnote{For vector-valued parameters \( \theta_{j} \), we have a Fisher information matrix and \( \text{det } I_{\Mc_{j}}(\theta_{j}) \) refers to the determinant of this matrix. This determinant is always non-negative, because the Fisher information matrix is always a positive semidefinite symmetric matrix. Intuitively, volumes and areas cannot be negative (Appendix~{\ref{appendixOrtho}}).} %
Hence, except for the observations \( x^{n}_{\obs} \), all quantities in the formulas depend on the model \( \Mc_{j} \). We made this explicit using a subscript \( j \) to indicate that the quantity, say, \( \theta_{j} \) belongs to model \( \Mc_{j} \).%
\footnote{For the sake of clarity, we will use different notations for the parameters within the different models. We introduce two models in this section: the model \( \Mc_{1} \) with parameter \( \theta_{1}=\vartheta \) which we pit against the model \( \Mc_{2} \) with parameter \( \theta_{2}=\alpha \).} %
For all three criteria, the model yielding the lowest criterion value is perceived as the model that generalizes best \citep{myungInpressModel}. 

Each of the three model selection criteria tries to strike a balance between model fit and model complexity. Model fit is expressed by the goodness-of-fit terms, which involves replacing the potential outcomes \( x^{n} \) and the unknown parameter \( \theta_{j} \) of the functional relationships \( f_{j} \) by the actually observed data \( x^{n}_{\obs} \), as in the Bayesian setting, and the maximum likelihood estimate \( \hat{\theta}_{j}(x^{n}_{\obs}) \), as in the frequentist setting.

The positive terms in the criteria account for model complexity. A penalization of model complexity is necessary, because the support in the data cannot be assessed by solely considering goodness-of-fit, as the ability to fit observations increases with model complexity (e.g., \citealp{roberts2000persuasive}). As a result, the more complex model necessarily leads to better fits but may in fact overfit the data. The overly complex model then captures idiosyncratic noise rather than general structure, resulting in poor model generalizability (\citealp{myung2000model}; \citealp{wagenmakers2006model}).

The focus in this section is to make intuitive how FIA acknowledges the trade-off between goodness-of-fit and model complexity in a principled manner by graphically illustrating this model selection procedure, see also \citet{balasubramanian1996geometric}, \citet{kass1989geometry}, \citet{myung2000counting}, and \citet{rissanen1996fisher}. We exemplify the concepts with simple multinomial processing tree (MPT) models (e.g., \citealp{batchelder1999theoretical}; \citealp{klauer2011flexibility}; \citealp{wu2010minimum}). For a more detailed treatment of the subject we refer to Appendix~{\ref{appendixMdl}}, \citet{deRooij2010luckiness}, \citet{grunwald2007minimum}, \citet{myung2006model}, and the references therein.

\subsubsection{The description length of a model}
Recall that each model specifies a functional relationship \( f_{j} \) between the potential outcomes of \( X \) and the parameters \( \theta_{j} \). This \( f_{j} \) is used to define a so-called \emph{normalized maximum likelihood} (NML) code. For the \( j \)th model its NML code is defined as%
\begin{align}
\label{NML}
p_{\text{NML}}(x^{n}_{\obs} \, | \, \Mc_{j}) = {f_{j}(x^{n}_{\obs} \, | \, \hat{\theta}_{j} ( x^{n}_{\obs})) \over \sum_{x^{n} \in \Xc^{n}} f_{j}(x^{n} \, | \, \hat{\theta}_{j}(x^{n}))},
\end{align}
where the sum in the denominator is over all possible outcomes \( x^{n} \) in \( \Xc^{n} \), and where \( \hat{\theta}_{j} \) refers to the MLE within model \( \Mc_{j} \). The NML code is a relative goodness-of-fit measure, as it compares the observed goodness-of-fit term against the sum of all possible goodness-of-fit terms. Note that the actual observations \( x^{n}_{\obs} \) only affect the numerator, by a plugin of \( x^{n}_{\obs} \) and its associated maximum likelihood estimate \( \hat{\theta}(x^{n}_{\obs}) \) into the functional relationship \( f_{j} \) belonging to model \( \Mc_{j} \). The sum in the denominator consists of the same plugins, but for every possible realization of \( X^{n} \).%
\footnote{As before, for continuous data, the sum is replaced by an integral.} %
Hence, the denominator can be interpreted as a measure of the model's collective goodness-of-fit or the model's fit capacity. Consequently, for every set of observations \( x^{n}_{\obs} \), the NML code outputs a number between zero and one that can be transformed into a non-negative number by taking the negative logarithm as%
\footnote{Quite deceivingly the minus sign actually makes this definition positive, as \( - \log (y) = \log (1/y) \geq 0 \) if \( 0 \leq y \leq 1 \).} %
\begin{align}
\label{descriptionLength}
- \log p_{\text{NML}}(x^{n}_{\obs} \, | \, \Mc_{j}) = - \log f_{j} \big (x^{n}_{\obs} \, | \, \hat{\theta}_{j}(x^{n}_{\obs}) \big ) + \underbrace{\log \sum f_{j}(x^{n} \, | \, \hat{\theta}_{j}(x^{n}))}_{\text{Model complexity}},
\end{align}
which is called the description length of model \( \Mc_{j} \). Within the MDL framework, the model with the shortest description length is the model that best describes the observed data \( x^{n}_{\obs} \).

The model complexity term is typically hard to compute, but \citet{rissanen1996fisher} showed that it can be well-approximated by the dimensionality and the geometrical complexity terms. That is,
\begin{align}
\nonumber
\text{FIA} = & - \log f_{j} \big (x^{n}_{\obs} \, | \, \hat{\theta}_{j}(x^{n}_{\obs}) \big ) + \frac{d_{j}}{2} \log \frac{n}{2 \pi} + \log \left ( \int_{\Theta} \sqrt{\text{det }  I_{\Mc_{j}}(\theta_{j})  } \, \der \theta_{j} \right ),
\end{align}
is an approximation of the description length of model \( \Mc_{j} \). The determinant is simply the absolute value when the number of free parameters \( d_{j} \) is equal to one. Furthermore, the integral in the geometrical complexity term coincides with the normalizing constant of the Jeffreys's prior, which represented the volume of the model. In other words, a model's fit capacity is proportional to its volume in model space as one would expect.

In sum, within the MDL philosophy, a model is selected if it yields the shortest description length, as this model uses the functional relationship \( f_{j} \) that best extracts the regularities from \( x^{n}_{\obs} \). As the description length is often hard to compute, we approximate it with FIA instead \citep{heck2014model}. To do so, we have to characterize (1) all possible outcomes of \( X \), (2) propose at least two models which will be pitted against each other, and (3) identify the model characteristics: the MLE \( \hat{\theta}_{j} \) corresponding to \( \Mc_{j} \), and its volume \( V_{\Mc_{j}} \). In the remainder of this section we show that FIA selects the model that is closest to the data with an additional penalty for model complexity.

\subsection{A new running example and the geometry of a random variable with $w=3$ outcomes}
To graphically illustrate the model selection procedure underlying MDL we introduce a random variable \( X \) that has \( w=3 \) number of potential outcomes. %
\begin{voorbeeld}[A psychological task with three outcomes]
\label{psychTask}
In the training phase of a source-memory task, the participant is presented with two lists of words on a computer screen. List \( \Lc \) is projected on the left-hand side and list \( \Rc \) is projected on the right-hand side. In the test phase, the participant is presented with two words, side by side, that can stem from either list, thus, \( ll, lr, rl, rr \). At each trial, the participant is asked to categorize these pairs as either:
\begin{itemize}
\item \( L \) meaning both words come from the left list, i.e., \( ``ll" \),
\item \( M \) meaning the words are mixed, i.e., \( ``lr" \) or \( ``rl" \),
\item \( R \) meaning both words come from the right list, i.e., \( ``rr" \).
\end{itemize}%
For simplicity we assume that the participant will be presented with \( n \) test pairs \( X^{n} \) of equal difficulty. \( \hfill \diamond \)
\end{voorbeeld}

\begin{figure}[h]
\centering
\includegraphics[width = 1 \textwidth]{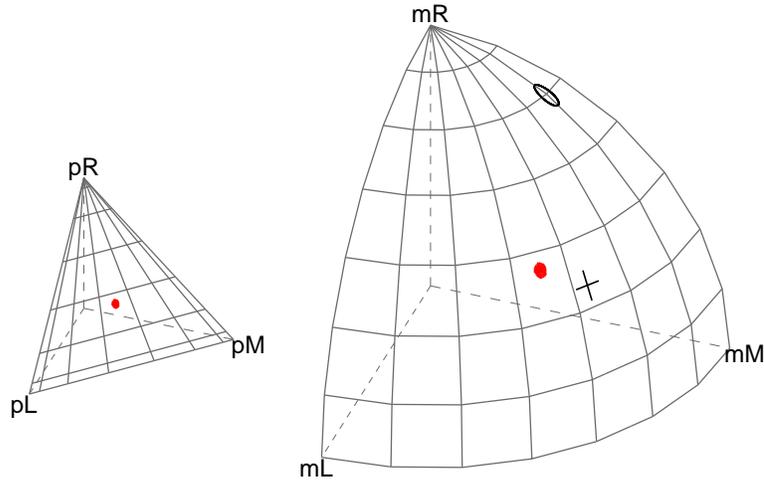}
\caption{Every point on the sphere corresponds to a pmf of a categorical distribution with \( w=3 \) categories. In particular, the (red) dot refers to the pmf \( p_{e}(x)=[1/3, 1/3, 1/3] \), the circle represents the pmf given by \( p(X)=[0.01, 0.18, 0.81] \), while the cross represents the pmf \( p(X)=[0.25, 0.5, 0.25] \).}
\label{L2a}
\end{figure}

For the graphical illustration of this new running example, we generalize the ideas presented in Section~{\ref{secCompleteSet2d}} from \( w=2 \) to \( w=3 \). Recall that a pmf of \( X \) with \( w \) number of outcomes can be written as a \( w \)-dimensional vector. For the task described above we know that a data generating pmf defines the three chances \( p(X)=[p(L), p(M), p(R)] \) with which \( X \) generates the outcomes \( [L, M, R] \) respectively.%
\footnote{As before we write \( p(X)=[p(L), p(M), p(R)] \) with a capital \( X \) to denote all the \( w \) number of chances simultaneously and we used the shorthand notation \( p(L)=p(X=L) \), \( p(M)=p(X=M) \) and \( p(R)=p(X=R) \).} %
As chances cannot be negative, (i) we require that \( 0 \leq p(x)=P(X=x) \) for every outcome \( x \) in \( \Xc \), and (ii) to explicitly convey that there are \( w=3 \) outcomes, and none more, these \( w=3 \) chances have to sum to one, that is, \( \sum_{x \in \Xc} p(x)=1 \). We call the largest set of functions that adhere to conditions (i) and (ii) the complete set of pmfs \( \Pc \). The three chances with which a pmf \( p(X) \) generates outcomes of \( X \) can be simultaneously represented in three-dimensional space with \( p(L)=P(X=L) \) on the left most axis, \( p(M)=P(X=M) \) on the right most axis and \( p(R)=P(X=R) \) on the vertical axis as shown in the left panel of \refFig{L2a}.%
\footnote{This is the three-dimensional generalization of \refFig{2dGeoA}.} %
In the most extreme case, we have the pmf \( p(X)=[1, 0, 0] \), \( p(X)=[0,1, 0] \) or \( p(X)=[0, 0,1] \), which correspond to the corners of the triangle indicated by \( pL, pM \) and \( pR \), respectively. These three extremes are linked by a triangular plane in the left panel of \refFig{L2a}. \emph{Any} pmf --and the true pmf \( p^{*}(X) \) in particular-- can be uniquely identified with a vector on the triangular plane and vice versa. For instance, a possible true pmf of \( X \) is \( p_{e}(X)=[1/3, 1/3, 1/3] \) (i.e., the outcomes \( L, M \) and \( R \) are generated with the same chance) depicted as a (red) dot on the simplex. 

This vector representation allows us to associate to each pmf of \( X \) the Euclidean norm. For instance, the representation in the left panel of \refFig{L2a} leads to an extreme pmf \( p(X)=[1, 0, 0] \) that is one unit long, while \( p_{e}(X)=[1/3, 1/3, 1/3] \) is only \( \sqrt{(1/3)^2+(1/3)^2+(1/3)^2} \approx 0.58 \) units away from the origin. As before, we can avoid this mismatch in lengths by considering the vectors \( m(X)=2 \sqrt{p(X)} \), instead. Any pmf that is identified as \( m(X) \) is now two units away from the origin. The model space \( \Mc \) is the collection of all transformed pmfs and represented as the surface of (the positive part of) the sphere in the right panel of \refFig{L2a}. By representing the set of all possible pmfs of \( X \) as \( m(X)=2\sqrt{p(X)} \), we adopted our intuitive notion of distance. As a result, the selection mechanism underlying MDL can be made intuitive by simply looking at the forthcoming plots.

\subsection{The individual-word and the only-mixed strategy}
\begin{figure}[h]
\centering
\includegraphics[width = 1 \textwidth]{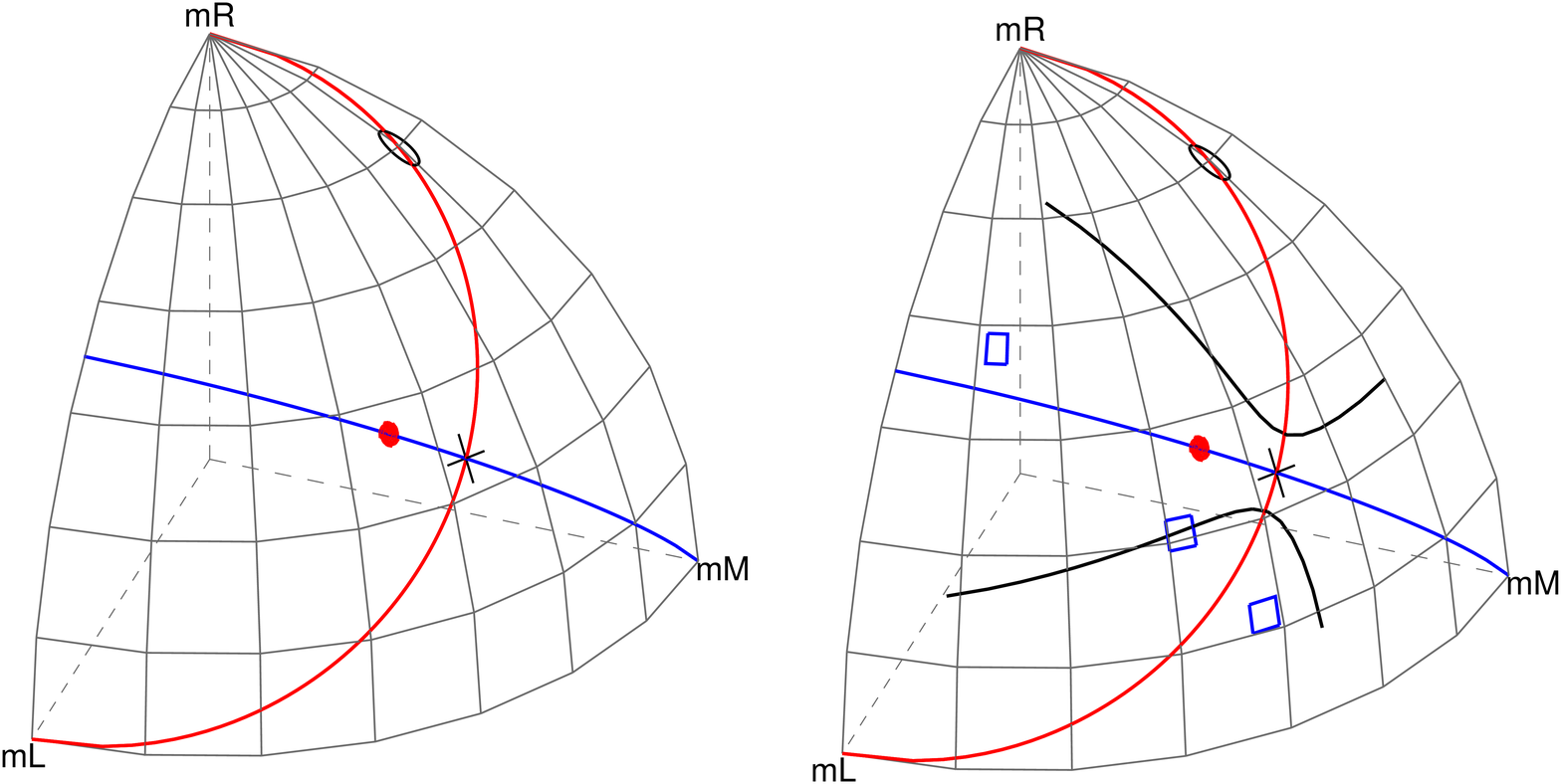}
\caption{Left panel: The set of pmfs that are defined by the individual-list strategy \( \Mc_{1} \) forms a curve that goes through both the cross and the circle, while the pmfs of the only-mixed strategy \( \Mc_{2} \) correspond to the curve that goes through both the cross and the dot. Right panel: The model selected by FIA can be thought of as the model closest to the empirical pmf with an additional penalty for model complexity. The selection between the individual-list and the only-mixed strategy by FIA based on \( n=30 \) trials is formalized by the additional curves --the only-mixed strategy is preferred over the individual-list strategy, when the observations yield an empirical pmf that lies between the two non-decision curves. The top, middle and bottom squares corresponding to the data sets \( x^{n}_{\obs, 1}, x^{n}_{\obs, 2} \) and \( x^{n}_{\obs, 3} \) in \refTab{mdlTable}, which are best suited to \( \Mc_{2} \), either, and \( \Mc_{1} \), respectively. The additional penalty is most noticeable at the cross, where the two models share a pmf. Observations with \( n=30 \) yielding an empirical pmf in this area are automatically assigned to the simpler model, i.e., the only-mixed strategy \( \Mc_{2} \).}
\label{modelSelection}
\end{figure}

To ease the exposition, we assume that both words presented to the participant come from the right list \( \Rc \), thus, \( ``rr" \) for the two models introduced below. %
\begin{figure}
\centering
\includegraphics[width = 1 \textwidth]{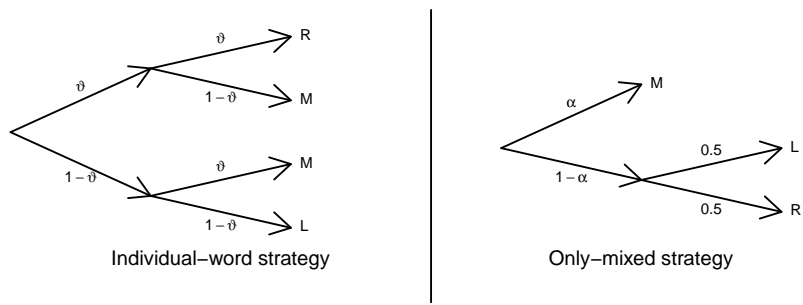}
\caption{Two MPT models that theorize how a participant chooses the outcomes \( L \), \( M \), or \( R \) in the source-memory task described in the main text. The left panel schematically describes the individual-word strategy, while the right model schematically describes the only-mixed strategy.}
\label{MPT}
\end{figure}%
As model \( \Mc_{1} \) we take the so-called individual-word strategy. Within this model \( \Mc_{1} \), the parameter is \( \theta_{1}=\vartheta \), which we interpret as the participant's ``right-list recognition ability''. With chance \( \vartheta \) the participant then correctly recognizes that the first word originates from the right list and repeats this procedure for the second word, after which the participant categorizes the word pair as \( L, M \), or \( R \), see the left panel of \refFig{MPT} for a schematic description of this strategy as a processing tree. Fixing the participant's ``right-list recognition ability'' \( \vartheta \) yields the following pmf
\begin{align}
\label{2BinomPdf}
f_{1}(X \, | \, \vartheta)=\big [ (1-\vartheta)^{2}, 2 \vartheta(1-\vartheta), \vartheta^{2} \big ].
\end{align}
For instance, when the participant's true ability is \( \vartheta^{*}=0.9 \), the three outcomes \( [L, M, R] \) are then generated with the following three chances \( f_{1}(X \, | \, 0.9)=[0.01, 0.18, 0.81] \), which is plotted as a circle in \refFig{L2a}. On the other hand, when \( \vartheta^{*}=0.5 \) the participant's generating pmf is then \( f_{1}(X \, | \, \vartheta=0.5)=[0.25, 0.5, 0.25] \), which is depicted as the cross in model space \( \Mc \). The set of pmfs so defined forms a curve that goes through both the cross and the circle, see the left panel of \refFig{modelSelection}. %

As a competing model \( \Mc_{2} \), we take the so-called only-mixed strategy. For the task described in \refEx{psychTask}, we might pose that participants from a certain clinical group are only capable of recognizing mixed word pairs and that they are unable to distinguish the pairs \( ``rr" \) from \( ``ll" \) resulting in a random guess between the responses \( L \) and \( R \), see the right panel of \refFig{MPT} for the processing tree. Within this model \( \Mc_{2} \) the parameter is \( \theta_{2}=\alpha \), which is interpreted as the participant's ``mixed-list differentiability skill'' and fixing it yields the following pmf 
\begin{align}
f_{2}(X \, | \, \alpha)= \big [(1-\alpha)/2, \alpha, (1-\alpha)/2 \big ]. %
\end{align}
For instance, when the participant's true differentiability is \( \alpha^{*}=1/3 \), the three outcomes \( [L, M, R] \) are then generated with the equal chances \( f_{2}(X \, | \, 1/3)=[1/3, 1/3, 1/3] \), which, as before, is plotted as the dot in \refFig{modelSelection}. On the other hand, when \( \alpha^{*}=0.5 \) the participant's generating pmf is then given by \( f_{2}(X \, | \, \alpha=0.5)=[0.25, 0.5, 0.25] \), i.e., the cross. The set of pmfs so defined forms a curve that goes through both the dot and the cross, see the left panel of \refFig{modelSelection}. 

The plots show that the models \( \Mc_{1} \) and \( \Mc_{2} \) are neither saturated nor nested, as the two models define proper subsets of \( \Mc \) and only overlap at the cross. Furthermore, the plots also show that \( \Mc_{1} \) and \( \Mc_{2} \) are both one-dimensional, as each model is represented as a line in model space. Hence, the dimensionality terms in all three information criteria are the same. Moreover, AIC and BIC will only discriminate these two models based on goodness-of-fit alone. This particular model comparison, thus, allows us to highlight the role Fisher information plays in the MDL model selection philosophy.

\subsection{Model characteristics}
\subsubsection{The maximum likelihood estimators}
\label{mleProjection}
For FIA we need to compute the goodness-of-fit terms, thus, we need to identify the MLEs for the parameters within each model. For the models at hand, the MLEs are %
%
\begin{align}
\hat{\theta}_{1}=\hat{\vartheta}=(Y_{M}+ 2 Y_{R})/(2n) \text{ for } \Mc_{1} \text{, and } \hat{\theta}_{2}=\hat{\alpha}=Y_{M}/n \text{ for } \Mc_{2},
\end{align}
where \( Y_{L}, Y_{M} \) and \( Y_{R}=n-Y_{L}-Y_{M} \) are the number of \( L, M \) and \( R \) responses in the data consisting of \( n \) trials. 

Estimation is a within model operation and it can be viewed as projecting the so-called \emph{empirical (i.e., observed) pmf} corresponding to the data onto the model. For iid data with \( w=3 \) outcomes the empirical pmf corresponding to \( x^{n}_{\obs} \) is defined as \( \hat{p}_{\obs}(X)= [y_{L}/n, y_{M}/n, y_{R} / n] \). Hence, the empirical pmf gives the relative occurrence of each outcome in the sample. For instance, the observations \( x^{n}_{\obs} \) consisting of \( [y_{L}=3, y_{M}=3, y_{R}=3] \) responses correspond to the observed pmf \( \hat{p}_{\obs}(X)=[1/3, 1/3, 1/3] \), i.e., the dot in \refFig{modelSelection}. Note that this observed pmf \( \hat{p}_{\obs}(X) \) does not reside on the curve of \( \Mc_{1} \). 

Nonetheless, when we use the MLE \( \hat{\vartheta} \) of \( \Mc_{1} \), we as researchers bestow the participant with a ``right-list recognition ability'' \( \vartheta \) and implicitly assume that she used the individual-word strategy to generate the observations. In other words, we only consider the pmfs on the curve of \( \Mc_{1} \) as viable explanations of how the participant generated her responses. For the data at hand, we have the estimate \( \hat{\vartheta}_{\obs}=0.5 \). If we were to generalize the observations \( x^{n}_{\obs} \) under \( \Mc_{1} \), we would then plug this estimate into the functional relationship \( f_{1} \) resulting in the predictive pmf \( f_{1}(X \, | \, \hat{\vartheta}_{\obs})=[0.25, 0.5, 0.25] \). Hence, even though the number of \( L, M \) and \( R \) responses were equal in the observations \( x^{n}_{\obs} \), under \( \Mc_{1} \) we expect that this participant will answer with twice as many \( M \) responses compared to the \( L \) and \( R \) responses in a next set of test items. Thus, for predictions, part of the data is ignored and considered as noise. 

Geometrically, the generalization \( f_{1}(X \, | \, \hat{\vartheta}_{\obs}) \) is a result of projecting the observed pmf \( \hat{p}_{\obs}(X) \), i.e., the dot, onto the cross that does reside on the curve of \( \Mc_{1} \).%
\footnote{This resulting pmf \( f_{1}(X \, | \, \hat{\vartheta}_{\obs}) \) is also known as the Kullback-Leibler projection of the empirical pmf \( \hat{p}_{\obs}(X) \) onto the model \( \Mc_{1} \). \citet{white1982maximum} used this projection to study the behavior of the MLE under model misspecification.} %
Observe that amongst all pmfs on \( \Mc_{1} \), the projected pmf is closest to the empirical pmf \( \hat{p}_{\obs}(X) \). Under \( \Mc_{1} \) the projected pmf \( f_{1}(X \, | \, \hat{\vartheta}_{\obs}) \), i.e., the cross, is perceived as structural, while any deviations from the curve of \( \Mc_{1} \) is labeled as noise. When generalizing the observations, we ignore noise. Hence, by estimating the parameter \( \vartheta \), we implicitly restrict our predictions to only those pmfs that are defined by \( \Mc_{1} \). Moreover, evaluating the prediction at \( x^{n}_{\obs} \) and, subsequently, taking the negative logarithm yields the goodness-of-fit term; in this case, \( - \log f_{1}(x^{n}_{\obs} \, | \, \hat{\vartheta}_{\obs}=0.5)=10.4 \). 

Which part of the data is perceived as structural or as noise depends on the model. For instance, when we use the MLE \( \hat{\alpha} \), we restrict our predictions to the pmfs of \( \Mc_{2} \). For the data at hand, we get \( \hat{\alpha}_{\obs}=1/3 \) and the plugin yields \( f_{2}(X \, | \, \hat{\alpha}_{\obs})=[1/3, 1/3, 1/3] \). Again, amongst all pmfs on \( \Mc_{2} \), the projected pmf is closest to the empirical pmf \( \hat{p}_{\obs}(X) \). In this case, the generalization under \( \Mc_{2} \) coincides with the observed pmf \( \hat{p}_{\obs}(X) \). Hence, under \( \Mc_{2} \) there is no noise, as the empirical pmf \( \hat{p}_{\obs}(X) \) was already on the model. Geometrically, this means that \( \Mc_{2} \) is closer to the empirical pmf than \( \Mc_{1} \), which results in a lower goodness-of-fit term \( - \log f_{2}(x^{n}_{\obs} \, | \, \hat{\alpha}_{\obs}=1/3)=9.9 \). 

This geometric interpretation allows us to make intuitive that data sets with the same goodness-of-fit terms will be as far from \( \Mc_{1} \) as from \( \Mc_{2} \). Equivalently, \( \Mc_{1} \) and \( \Mc_{2} \) identify the same amount of noise within \( x^{n}_{\obs} \), when the two models fit the observations equally well. For instance, \refFig{modelSelection} shows that observations \( x^{n}_{\obs} \) with an empirical pmf \( \hat{p}_{\obs}(X) = [0.25, 0.5, 0.25] \) are equally far from \( \Mc_{1} \) as from \( \Mc_{2} \). Note that the closest pmf on \( \Mc_{1} \) and \( \Mc_{2} \) are both equal to the empirical pmf, as \( f_{1}(X \, | \, \hat{\vartheta}_{\obs}=0.5)=\hat{p}_{\obs}(X)=f_{2}(X \, | \, \hat{\alpha}_{\obs}=1/2) \). As a result, the two goodness-of-fit terms will be equal to each other. 

In sum, goodness-of-fit measures a model's proximity to the observed data. Consequently, models that take up more volume in model space will be able to be closer to a larger number of data sets. In particular, when, say, \( \Mc_{3} \) is nested within \( \Mc_{4} \), this means that the distance between \( \hat{p}_{\obs}(X) \) and \( \Mc_{3} \) (noise) is at least the distance between \( \hat{p}_{\obs}(X) \) and \( \Mc_{4} \). Equivalently, for any data set, \( \Mc_{4} \) will automatically label more of the observations as structural. Models that excessively identify parts of the observations as structural are known to overfit the data. Overfitting has an adverse effect on generalizability, especially when \( n \) is small, as \( \hat{p}_{\obs}(X) \) is then dominated by sampling error. In effect, the more voluminous model will then use this sampling error, rather than the structure, for its predictions. To guard ourselves from overfitting, thus, bad generalizability, the information criteria AIC, BIC and FIA all penalize for model complexity. AIC and BIC only do this via the dimensionality terms, while FIA also take the models' volumes into account. 

\subsubsection{Geometrical complexity}
For both models the dimensionality term is given by \( \tfrac{1}{2} \log ( \tfrac{n}{2 \pi}) \). Recall that the geometrical complexity term is the logarithm of the model's volume, which for the individual-word and the only-mixed strategy are given by %
\begin{align}
V_{\Mc_{1}}& =\int_{0}^{1} \sqrt{I_{\Mc_{1}} (\theta)} \der \theta = \sqrt{2} \pi \text{ and } \\
V_{\Mc_{2}}& =\int_{0}^{1} \sqrt{I_{\Mc_{2}}(\alpha)} \der \alpha = \pi ,
\end{align}
%
respectively. Hence, the individual-word strategy is a more complex model, because it has a larger volume, thus, capacity to fit data compared to the only-mixed strategy. After taking logs, we see that the individual-word strategy incurs an additional penalty of \( 1/2 \log(2) \) compared to the only-mixed strategy. 

\subsection{Model selection based on the minimum description length principle}
With all model characteristics at hand, we only need observations to illustrate that MDL model selection boils down to selecting the model that is closest to the observations with an additional penalty for model complexity. %
\begin{table}[h]
\caption{The description lengths for three observations \( x^{n}_{\obs}=[y_{L}, y_{M}, y_{R}] \), where \( y_{L}, y_{M}, y_{R} \) are the number of observed responses \( L, M \) and \( R \) respectively.}
\centering
\label{mdlTable}
\begin{tabular}{llll}
  \hline
  \( x^{n}_{\obs}=[y_{L}, y_{M}, y_{R}] \) & \( \text{FIA}_{\Mc_{1}}(x^{n}_{\obs} )  \) & \( \text{FIA}_{\Mc_{2}}(x^{n}_{\obs})  \) & Preferred model \\ %
  \hline
  \( x^{n}_{\obs, 1}=[12, 1, 17] \) &  42 & 26 &   \( \Mc_{2} \) \\ %
  \( x^{n}_{\obs, 2}=[14, 10, 6] \) &  34 &  34 & tie \\ %
  \( x^{n}_{\obs, 3}=[12, 16, 2] \) &  29 & 32  & \( \Mc_{1} \) \\ %
   \hline
\end{tabular}
\end{table}
\refTab{mdlTable} shows three data sets \( x^{n}_{\obs, 1}, x^{n}_{\obs, 2}, x^{n}_{\obs, 3} \) with \( n=30 \) observations. The three associated empirical pmfs are plotted as the top, middle and lower rectangles in the right panel of \refFig{modelSelection}, respectively. \refTab{mdlTable} also shows the approximation of each model's description length using FIA. Note that the first observed pmf, the top rectangle in \refFig{modelSelection}, is closer to \( \Mc_{2} \) than to \( \Mc_{1} \), while the third empirical pmf, the lower rectangle, is closer to \( \Mc_{1} \). Of particular interest is the middle rectangle, which lies on an additional black curve that we refer to as a non-decision curve; observations that correspond to an empirical pmf that lies on this curve are described equally well by \( \Mc_{1} \) and \( \Mc_{2} \). For this specific comparison, we have the following decision rule: FIA selects \( \Mc_{2} \) as the preferred model whenever the observations correspond to an empirical pmf between the two non-decision curves, otherwise, FIA selects \( \Mc_{1} \). \refFig{modelSelection} shows that FIA, indeed, selects the model that is closest to the data except in the area where the two models overlap --observations %
consisting of \( n=30 \) trials with an empirical pmf near the cross are considered better described by the simpler model \( \Mc_{2} \). Hence, this yields an incorrect decision even when the empirical pmf is exactly equal to the true data generating pmf that is given by, say, \( f_{1}(X \, | \, \vartheta=0.51) \). This automatic preference for the simpler model, however, decreases as \( n \) increases. %
\begin{figure}[h]
\centering
\includegraphics[width = 1 \textwidth]{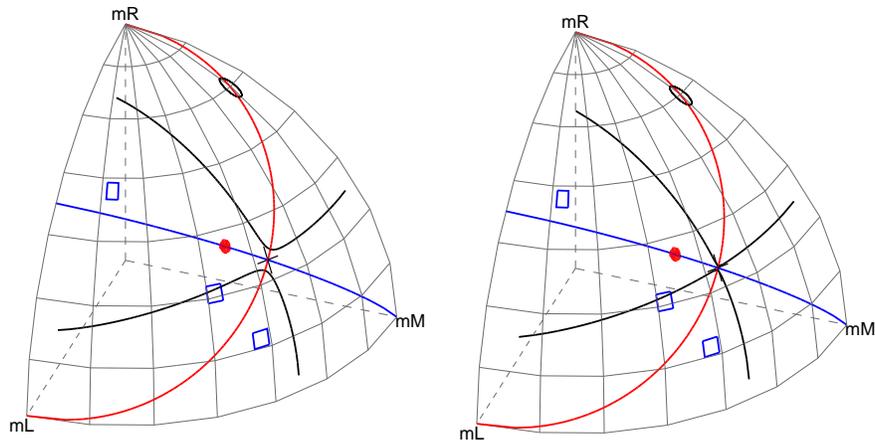}
\caption{For \( n \) large the additional penalty for model complexity becomes irrelevant. The plotted non-decision curves are based on \( n=120 \) and \( n =10{,}000 \) trials in the left and right panel respectively. In the right panel only the goodness-of-fit matters in the model comparison. The model selected is then the model that is closest to the observations.}
\label{modelSelection2}
\end{figure}%
The left and right panel of \refFig{modelSelection2} show the non-decision curves when \( n=120 \) and \( n \) (extremely) large, respectively. As a result of moving non-decision bounds, the data set \( x^{n}_{\obs, 4}=[56, 40, 24] \) that has the same observed pmf as \( x^{n}_{\obs, 2} \), i.e., the middle rectangle, will now be better described by model \( \Mc_{1} \). 

For (extremely) large \( n \), the additional penalty due to \( \Mc_{1} \) being more voluptuous than \( \Mc_{2} \) becomes irrelevant and the sphere is then separated into quadrants: observations corresponding to an empirical pmf in the top-left or bottom-right quadrant are better suited to the only-mixed strategy, while the top-right and bottom-left quadrants indicate a preference for the individual-word strategy \( \Mc_{1} \). Note that pmfs on the non-decision curves in the right panel of \refFig{modelSelection2} are as far apart from \( \Mc_{1} \) as from \( \Mc_{2} \), which agrees with our geometric interpretation of goodness-of-fit as a measure of the model's proximity to the data. This quadrant division is only based on the two models' goodness-of-fit terms and yields the same selection as one would get from BIC (e.g., \citealp{rissanen1996fisher}). %
For large \( n \), FIA, thus, selects the model that is closest to the empirical pmf. This behavior is desirable, because asymptotically the empirical pmf is not distinguishable from the true data generating pmf. As such, the model that is closest to the empirical pmf will then also be closest to the true pmf. Hence, FIA asymptotically selects the model that is closest to the true pmf. As a result, the projected pmf within the closest model is then expected to yield the best predictions amongst the competing models. 

\subsection{Fisher information and generalizability}
Model selection by MDL is sometimes perceived as a formalization of Occam's razor (e.g., \citealp{balasubramanian1996geometric, grunwald1998minimum}), a principle that states that the most parsimonious model should be chosen when the models under consideration fit the observed data equally well. This preference for the parsimonious model is based on the belief that the simpler model is better at predicting new (as yet unseen) data coming from the same source, as was shown by \citet{pitt2002toward} with simulated data.

To make intuitive why the more parsimonious model, on average, leads to better predictions, we assume, for simplicity, that the true data generating pmf is given by \( f(X \, | \, \theta^{*}) \), thus, the existence of a true parameter value \( \theta^{*} \). As the observations are expected to be contaminated with sampling error, we also expect an estimation error, i.e., a distance \( \der \theta \) between the maximum likelihood estimate \( \hat{\theta}_{\obs} \) and the true \( \theta^{*} \). Recall that in the construction of Jeffreys's prior Fisher information was used to convert displacement in model space to distances on parameter space. Conversely, Fisher information transforms the estimation error in parameter space to a generalization error in model space. Moreover, the larger the Fisher information at \( \theta^{*} \) is, the more it will expand the estimation error into a displacement between the prediction \( f(X \, | \, \hat{\theta}_{\obs}) \) and the true pmf \( f(X \, | \, \theta^{*}) \). Thus, a larger Fisher information at \( \theta^{*} \) will push the prediction further from the true pmf resulting in a bad generalization. Smaller models have, on average, a smaller Fisher information at \( \theta^{*} \) and will therefore lead to more stable predictions that are closer to the true data generating pmf. %
Note that the generalization scheme based on the MLE plugin \( f(X \, | \, \hat{\theta}_{\obs}) \) ignores the error at each generalization step. The Bayesian counterpart, on the other hand, does take these errors into account, see \citet{dawid2011posterior}, \citet{ly2017replication}, \citet{marsman2016four} and see \citet{erven2012catching}, \citet{grunwald2016fast}, \citet{van2014almost}, \citet{wagenmakers2006accumulative} for a prequential view of generalizability.

\section{Concluding Comments}
Fisher information is a central statistical concept that is of considerable relevance for mathematical psychologists. We illustrated the use of Fisher information in three different statistical paradigms: in the frequentist paradigm, Fisher information was used to construct hypothesis tests and confidence intervals; in the Bayesian paradigm, Fisher information was used to specify a default, parameterization-invariant prior distribution; lastly, in the paradigm of information theory, data compression, and minimum description length, Fisher information was used to measure model complexity. Note that these three paradigms highlight three uses of the functional relationship \( f \) between potential observations \( x^{n} \) and the parameters \( \theta \). Firstly, in the frequentist setting, the second argument was fixed at a supposedly known parameter value \( \theta_{0} \) or \( \hat{\theta}_{\obs} \) resulting in a probability mass function, a function of the potential outcomes \( f ( \cdot \, \, | \, \theta_{0} ) \). Secondly, in the Bayesian setting, the first argument was fixed at the observed data resulting in a likelihood function, a function of the parameters \( f ( x_{\obs} \, | \, \cdot) \). Lastly, in the information geometric setting both arguments were free to vary, i.e., \( f ( \cdot \, | \, \cdot ) \) and plugged in by the observed data and the maximum likelihood estimate. 

To ease the exposition we only considered Fisher information of one-dimensional parameters. The generalization of the concepts introduced here to vector valued \( \theta \) can be found in the appendix. A complete treatment of all the uses of Fisher information throughout statistics would require a book (e.g., \citealp{frieden2004science}) rather than a tutorial article. Due to the vastness of the subject, the present account is by no means comprehensive. Our goal was to use concrete examples to provide more insight about Fisher information, something that may benefit psychologists who propose, develop, and compare mathematical models for psychological processes. Other uses of Fisher information are in the detection of model misspecification (\citealp{golden1995making}; \citealp{golden2000statistical}; \citealp{waldorp2005wald}; \citealp{waldorp2009robust}; \citealp{waldorp2011effective}; \citealp{white1982maximum}), in the reconciliation of frequentist and Bayesian estimation methods through the Bernstein-von Mises theorem (\citealp{bickel2012semiparametric}; \citealp{rivoirard2012bernstein}; \citealp{van1998asymptotic}; \citealp{yang2000asymptotics}), in statistical decision theory (e.g., \citealp{berger1985statistical}; \citealp{hajek1972local}; \citealp{korostelev2011mathematical}; \citealp{ray2016minimax}; \citealp{wald1949statistical}), in the specification of objective priors for more complex models (e.g., \citealp{ghosal1997non}; \citealp{grazian2015jeffreys}; \citealp{kleijn2013criteria}), and computational statistics and generalized MCMC sampling in particular (e.g., \citealp{banterle2015accelerating}; \citealp{girolami2011riemann}; \citealp{grazian2014approximate}; \citealp{gronau2017tutorial}).

In sum, Fisher information is a key concept in statistical modeling. We hope to have provided an accessible and concrete tutorial article that explains the concept and some of its uses for applications that are of particular interest to mathematical psychologists.


\clearpage
\bibliographystyle{imsart-nameyear}
\bibliography{Alexander}

\appendix
\section{Generalization to Vector-Valued Parameters: The Fisher Information Matrix}
\label{appendixFiMatrix}
Let \( X \) be a random variable, \( \vec{\theta}=(\theta_{1}, \ldots, \theta_{d}) \) a vector of parameters, and \( f \) a functional relationship that relates \( \vec{\theta} \) to the potential outcomes \( x \) of \( X \). As before, it is assumed that by fixing \( \vec{\theta} \) in \( f \) we get the pmf \( p_{\vec{\theta}}(x)=f(x \, | \, \vec{\theta}) \), which is a function of \( x \). The pmf \( p_{\vec{\theta}}(x) \) fully determines the chances with which \( X \) takes on the events in the outcome space \( \Xc \). The Fisher information of the vector \( \vec{\theta} \in \R^{d} \) is a positive semidefinite symmetric matrix of dimension \( d \times d \) with the entry at the \( i \)th row and \( j \)th column given by %
\begin{align}
I_{X}(\vec{\theta})_{i, j} = & \Cov \Big ( \dot{l}(X \, | \, \vec{\theta}), \dot{l}^{T}(X \, | \, \vec{\theta}) \Big )_{i, j}, \\
= & \begin{cases} 
\sum_{x \in \Xc} \Big ( \tfrac{\partial }{\partial \theta_{i}} l(x \, | \, \vec{\theta}), \tfrac{\partial }{\partial \theta_{j}} l(x \, | \, \vec{\theta}) \Big ) p_{\vec{\theta}}(x) & \text{if } X \text{ is discrete,} \\
\int_{x \in \Xc} \Big ( \tfrac{\partial }{\partial \theta_{i}} l(x \, | \, \vec{\theta}), \tfrac{\partial }{\partial \theta_{j}} l(x \, | \, \vec{\theta}) \Big ) p_{\vec{\theta}}(x) \der x & \text{if } X \text{ is continuous.}
\end{cases}
\end{align}
where \( l(x \, | \, \vec{\theta})=\log f (x \, | \, \vec{\theta}) \) is the log-likelihood function, \( \tfrac{ \partial}{ \partial \theta_{i}} l(x \, | \, \vec{\theta}) \) is the score function, that is, the partial derivative with respect to the \( i \)th component of the vector \( \vec{\theta} \) and the dot is short-hand notation for the vector of the partial derivatives with respect to \( \theta=(\theta_{1}, \ldots, \theta_{d}) \). Thus, \( \dot{l}(x \, | \, \vec{\theta}) \) is a \( d \times 1 \) column vector of score functions, while \( \dot{l}^{T}(x \, | \, \vec{\theta}) \) is a \( 1 \times d \) row vector of score functions at the outcome \( x \). The partial derivative is evaluated at \( \vec{\theta} \), the same \( \vec{\theta} \) that is used in the pmf \( p_{\vec{\theta}}(x) \) for the weighting. In Appendix~{\ref{appendixReg}} it is shown that the score functions are expected to be zero, which explains why \( I_{X}(\vec{\theta}) \) is a covariance matrix. 

Under mild regularity conditions the \( i,j \)th entry of the Fisher information matrix can be equivalently calculated via the negative expectation of the second order partial derivates, that is,
\begin{align}
\label{FIM}
I_{X}(\vec{\theta})_{i, j} = & - E \Big ( \tfrac{\partial^{2}}{\partial \theta_{i} \partial \theta_{j}} l(X \, | \, \vec{\theta}) \Big ), \\
= & %
\begin{cases} 
- \sum_{x \in \Xc} \tfrac{\partial^{2} }{\partial \theta_{i} \partial \theta_{j}} \log f(x \, | \, \vec{\theta}) p_{\vec{\theta}}(x) & \text{if } X \text{ is discrete,} \\
- \int_{x \in \Xc} \tfrac{\partial^{2} }{\partial \theta_{i} \partial \theta_{j}} \log f(x \, | \, \vec{\theta}) p_{\vec{\theta}}(x) \der x & \text{if } X \text{ is continuous.}
\end{cases}
\end{align}%
Note that the sum (thus, integral in the continuous case) is with respect to the outcomes \( x \) of \( X \). 

\begin{voorbeeld}[Fisher information for normally distributed random variables]
When \( X \) is normally distributed, i.e., \( X \sim \Nc(\mu, \sigma^{2}) \), it has the following probability density function (pdf) 
\begin{align}
f(x \, | \, \vec{\theta}) = \frac{1}{\sqrt{2 \pi} \sigma } \exp \big ( - \frac{1}{2 \sigma^{2}} ( x - \mu)^{2} \big ),
\end{align}
where the parameters are collected into the vector \( \vec{\theta} ={\mu \choose \sigma} \), with \( \mu \in \R \) and \( \sigma > 0 \). The score vector at a specific \( \theta= { \mu \choose \sigma} \) is the following vector of functions of \( x \) %
\begin{align}
\dot{l}(x \, | \, \vec{\theta}) = %
 \begin{pmatrix} \tfrac{\partial} { \partial \mu } l(x \, | \, \vec{\theta}) \\ \tfrac{\partial} { \partial \sigma } l(x \, | \, \vec{\theta}) \end{pmatrix} = 
 \begin{pmatrix} \frac{x- \mu}{\sigma^{2}} \\  \frac{(x-\mu)^{2}}{\sigma^{3}} -\frac{1}{\sigma} \end{pmatrix} .
\end{align}
The unit Fisher information matrix \( I_{X}(\vec{\theta}) \) is a \( 2 \times 2 \) symmetric positive semidefinite matrix, consisting of expectations of partial derivatives. Equivalently, \( I_{X}(\vec{\theta}) \) can be calculated using the second order partials derivatives %
\begin{align}
\label{FINorm}
I_{X}(\vec{\theta}) = - E 
\begin{pmatrix}
\frac{\partial^{2}}{\partial \mu \partial \mu} \log f(x \, | \, \mu, \sigma^{2}) & \frac{\partial^{2}}{\partial \mu \partial \sigma } \log f(x \, | \, \mu, \sigma) \\
\frac{\partial^{2}}{\partial \sigma \partial \mu} \log f(x \, | \, \mu, \sigma) & \frac{\partial^{2}}{\partial \sigma \partial \sigma} \log f(x \, | \, \mu, \sigma)
\end{pmatrix}
 =
 \begin{pmatrix}
\frac{1}{\sigma^{2}} & 0 \\
0 & \frac{2}{\sigma^{2}}
\end{pmatrix}.
\end{align}%
The off-diagonal elements are in general not zero. If the \( i,j \)th entry is zero we say that \( \theta_{i} \) and \( \theta_{j} \) are orthogonal to each other, see Appendix~{\ref{appendixOrtho}} below. \( \hfill \diamond \)
\end{voorbeeld}%
For iid trials \( X^{n} = (X_{1}, \ldots, X_{n}) \) with \( X \sim p_{\theta}(x) \), the Fisher information matrix for \( X^{n} \) is given by \( I_{X^{n}}(\vec{\theta}) = n I_{X}(\vec{\theta}) \). Thus, for vector-valued parameters \( \vec{\theta} \) the Fisher information matrix remains additive. 

In the remainder of the text, we simply use \( \theta \) for both one-dimensional and vector-valued parameters. Similarly, depending on the context it should be clear whether \( I_{X}(\theta) \) is a number or a matrix. 

\section{Frequentist Statistics based on Asymptotic Normality}
The construction of the hypothesis tests and confidence intervals in the frequentist section were all based on the MLE being asymptotically normal.

\subsection{Asymptotic normality of the MLE for vector-valued parameters}
\label{appendixFreq}
For so-called regular parametric models, see Appendix~{\ref{appendixReg}}, the MLE for vector-valued parameters \( \theta \) converges in distribution to a multivariate normal distribution, that is, %
\begin{align}
\label{mleCLTM0}
\sqrt{n} ( \hat{\theta} - \theta^{*}) \overset{D}{\rightarrow} \Nc_{d} \Big ( 0, I^{-1}_{X}(\theta^{*}) \Big ) , \text{ as } n \rightarrow \infty,
\end{align}
where \( \Nc_{d} \) is a \( d \)-dimensional multivariate normal distribution, and \( I_{X}^{-1}(\theta^{*}) \) the inverse Fisher information matrix at the true value \( \theta^{*} \). For \( n \) large enough, we can, thus, approximate the sampling distribution of the ``error'' of the MLE by a normal distribution, thus, 
\begin{align}
\label{mleCLTM}
(\hat{\theta} - \theta^{*}) \overset{D}{\approx} \Nc_{d} \Big ( 0 , \tfrac{1}{n} I^{-1}_{X}(\theta^{*}) \Big ) \text{, we repeat, approximately.}
\end{align}
In practice, we fix \( n \) and replace the true sampling distribution by this normal distribution. Hence, we incur an approximation error that is only negligible whenever \( n \) is large enough. What constitutes \( n \) large enough depends on the true data generating pmf \( p^{*}(x) \) that is unknown in practice. In other words, the hypothesis tests and confidence intervals given in the main text based on the replacement of the true sampling distribution by this normal distribution might not be appropriate. In particular, this means that a hypothesis tests at a significance level of 5\% based on the asymptotic normal distribution, instead of the true sampling distribution, might actually yield a type 1 error rate of, say, 42\%. Similarly, as a result of the approximation error, a 95\%-confidence interval might only encapsulate the true parameter in, say, 20\% of the time that we repeat the experiment. 

\subsection{Asymptotic normality of the MLE and the central limit theorem}
Asymptotic normality of the MLE can be thought of as a refinement of the central limit theorem. The (Lindeberg-L\'{e}vy) CLT is a general statement about the sampling distribution of the sample mean estimator \( \bar{X}=\tfrac{1}{n} \sum_{i=1}^{n} X_{i} \) based on iid trials of \( X \) with common population mean \( \theta=E(X) \) and variance \( \Var(X) < \infty \). More specifically, the CLT states that, with a proper scaling, the sample mean \( \bar{X} \) centred around the true \( \theta^{*} \) will converge in distribution to a normal distribution, that is, \( \sqrt{n} (\bar{X} - \theta^{*} )  \overset{D}{\rightarrow} \Nc \big ( 0 , \Var(X) \big ) \). In practice, we replace the true sampling distribution by this normal distribution at fixed \( n \) and hope that \( n \) is large enough. Hence, for fixed \( n \) we then suppose that the ``error'' is distributed as \( (\bar{X} - \theta^{*}) \overset{D}{\approx} \Nc ( 0, \tfrac{1}{n} \Var(X)) \) and we ignore the approximation error. %
In particular, when we know that the population variance is \( \Var(X)=1 \), we then know that we require an experiment with \( n=100 \) samples for \( \bar{X} \) to generate estimates within \( 0.196 \) distance from \( \theta \) with approximately 95\% chance, that is, \( P( | \bar{X} - \theta | \leq 0.196) \approx 0.95 \).%
\footnote{As before, chance refers to the relative frequency, that is, when we repeat the experiment \( k=200 \) times, each with \( n=100 \), we get \( k \) number of estimates and approximately 95\% of these \( k \) number of estimates are then expected to be within \( 0.196 \) distance away from the true population mean \( \theta^{*} \).} %
This calculation was based on our knowledge of the normal distribution \( \Nc(0, 0.01) \), which has its 97.5\% quantile at 0.196. In the examples below we re-use this calculation by matching the asymptotic variances to \( 0.01 \).%
\footnote{Technically, an asymptotic variance is free of \( n \), but we mean the approximate variance at finite \( n \). For the CLT this means \( \tfrac{1}{n} \sigma^{2} \).} %
The 95\% statement only holds approximately, because we do not know whether \( n=100 \) is large enough for the CLT to hold, i.e., this probability could be well below 23\%. Note that the CLT holds under very general conditions; the population mean and variance both need to exist, i.e., be finite. The distributional form of \( X \) is irrelevant for the statement of the CLT.

On the other hand, to even compute the MLE we not only require that the population quantities to exists and be finite, but we also need to know the functional relationship \( f \) that relates these parameters to the outcomes of \( X \). When we assume more (and nature adheres to these additional conditions), we know more, and are then able to give stronger statements. We give three examples.

\begin{voorbeeld}[Asymptotic normality of the MLE vs the CLT: The Gaussian distribution]
\label{gaussEx}
If \( X \) has a Gaussian (normal) distribution, i.e., \( X \sim \Nc(\theta, \sigma^{2}) \), with \( \sigma^{2} \) known, then the MLE is the sample mean and the unit Fisher information is \( I_{X}(\theta) = 1/ \sigma^{2} \). Asymptotic normality of the MLE leads to the same statement as the CLT, that is, \( \sqrt{n}(\hat{\theta} - \theta^{*}) \overset{D}{\rightarrow} \Nc ( 0,  \sigma^{2}) \). Hence, asymptotically we do not gain anything by going from the CLT to asymptotic normality of the MLE. The additional knowledge of \( f(x  \, | \, \theta) \) being normal does, however, allow us to come to the rare conclusion that the normal approximation holds exactly for every finite \( n \), thus, \( (\hat{\theta} - \theta^{*} )\overset{D}{=} \Nc ( 0, \tfrac{1}{n} \sigma^{2}) \). In all other cases, whenever \( X \not \sim \Nc(\theta, \sigma^{2}) \), we always have an approximation.%
\footnote{This is a direct result of Cram\'{e}r's theorem that states that whenever \( X \) is independent of \( Y \) and \( Z=X+Y \) with \( Z \) a normal distribution, then \( X \) and \( Y \) themselves are necessarily normally distributed.} %
Thus, whenever \( \sigma^{2}=1 \) and \( n=100 \) we know that \( P(| \hat{\theta} - \theta^{*}| \leq 0.196) = 0.95 \) holds exactly. \( \hfill \diamond \)
\end{voorbeeld}

\begin{voorbeeld}[Asymptotic normality of the MLE vs the CLT: The Laplace distribution]
\label{laplaceEx}
If \( X \) has a Laplace distribution with scale \( b \), i.e., \( X \sim \text{Laplace}(\theta, b) \), then its population mean and variance are \( \theta=E(X) \) and \( 2 b^{2} = \Var(X) \), respectively. 

In this case, the MLE is the sample median \( \hat{M} \) and the unit Fisher information is \( I_{X}(\theta) = 1/b^{2} \). Asymptotic normality of the MLE implies that we can approximate the sampling distribution by the normal distribution, that is, \( ( \hat{\theta} - \theta^{*} ) \overset{D}{\approx} \Nc(0, \tfrac{1}{n} b^{2} ) \), when \( n \) is large enough. Given that the population variance is \( \Var(X)=1 \), we know that \( b=1/\sqrt{2} \), yielding a variance of \( \tfrac{1}{2 n} \) in our normal approximation to the sampling distribution. Matching this variance to \( 0.01 \) shows that we now require only \( n=50 \) samples for the estimator to generate estimates within 0.196 distance away from the true value \( \theta^{*} \) with 95\% chance. As before, the validity of this statement only holds approximately, i.e., whenever the normal approximation to the sampling distribution of the MLE at \( n=50 \) is not too bad.

Hence, the additional knowledge of \( f(x \, | \, \theta) \) being Laplace allows us to use an estimator, i.e., the MLE, that has a lower asymptotic variance. Exploiting this knowledge allowed us to design an experiment with twice as few participants. \( \hfill \diamond \)
\end{voorbeeld}

\begin{voorbeeld}[Asymptotic normality of the MLE vs the CLT: The Cauchy distribution]
\label{cauchyEx}
If \( X \) has a Cauchy distribution centred around \( \theta \) with scale \( 1 \), i.e., \( X \sim \text{Cauchy}(\theta, 1) \), then \( X \) does not have a finite population variance, nor a finite population mean. As such, the CLT cannot be used. Even worse, \citet{fisher1922mathematical} showed that the sample mean as an estimator for \( \theta \) is in this case useless, as the sampling distribution of the sample mean is a Cauchy distribution that does not depend on \( n \), namely, \( \bar{X} \sim \text{Cauchy}(\theta, 1) \). As such, using the first observation alone to estimate \( \theta \) is as good as combining the information of \( n=100 \) samples in the sample mean estimator. Hence, after seeing the first observation no additional information about \( \theta \) is gained using the sample mean \( \bar{X} \), not even if we increase \( n \). 

The sample median estimator \( \hat{M} \) performs better. Again, \citet{fisher1922mathematical} already knew that for \( n \) large enough that \( (\hat{M}  - \theta^{*}) \overset{D}{\approx} \Nc(0, \tfrac{1}{n}  \tfrac{\pi^{2}}{2} ) \). The MLE is even better, but unfortunately, in this case, it cannot be given as an explicit function of the data.%
\footnote{Given observations \( x^{n}_{\obs} \) the maximum likelihood estimate \( \hat{\theta}_{\obs} \) is the number for which the score function \( \dot{l}(x^{n}_{\obs}  \, | \, \theta) = \sum_{i=1}^{n} \tfrac{2 (x_{\obs, i} - \theta)}{ 1 + (x_{\obs, i} - \theta)^{2}} \) is zero. This optimization cannot be solved analytically and there are \( 2 n \) solutions to this equation.} %
The Fisher information can be given explicitly, namely, \( I_{X}(\theta)=1/2 \). Asymptotic normality of the MLE implies that \( ( \hat{\theta} - \theta^{*} ) \overset{D}{\approx} \Nc(0, \tfrac{1}{n} 2) \), when \( n \) is large enough. Matching the variances in the approximation based on the normal distribution to \( 0.01 \) shows that we require \( n=25 \pi^{2} \approx 247 \) for the sample median and \( n=200 \) samples for the MLE to generate estimates within 0.196 distance away from the true value of value \( \theta^{*} \) with approximate 95\% chance. \( \hfill \diamond \)
\end{voorbeeld}

\subsection{Efficiency of the MLE: The H\'{a}jek-LeCam convolution theorem and the Cram\'{e}r-Fr\'{e}chet-Rao information lower bound} 
The previous examples showed that the MLE is an estimator that leads to a smaller sample size requirement, because it is the estimator with the lower asymptotic variance. This lower asymptotic variance is a result of the MLE making explicit use of the functional relationship between the samples \( x^{n}_{\obs} \) and the target \( \theta \) in the population. Given any such \( f \), one might wonder whether the MLE is the estimator with the \emph{lowest possible} asymptotic variance. The answer is affirmative, whenever we restrict ourselves to the broad class of so-called regular estimators. 

A \emph{regular estimator} \( T_{n}=t_{n}(X_{n}) \) is a function of the data that has a limiting distribution that does not change too much, whenever we change the parameters in the neighborhood of the true value \( \theta^{*} \), see \citet[][{p. 115}]{van1998asymptotic} for a precise definition. The H\'{a}jek-LeCam convolution theorem characterizes the aforementioned limiting distribution as a convolution, i.e., a sum of, the independent statistics \( \Delta_{\theta^{*}} \) and \( Z_{\theta^{*}} \). That is, for any regular estimator \( T_{n} \) and every possible true value \( \theta^{*} \) we have %
\begin{align}
\sqrt{n} ( T_{n} - \theta^{*}) \overset{D}{\rightarrow} \Delta_{\theta^{*}} + Z_{\theta^{*}}, \text{ as } n \rightarrow \infty, 
\end{align}
where \( Z_{\theta^{*}} \sim \Nc \big (0, I^{-1}_{X}  (\theta^{*}) \big ) \) and where \( \Delta_{\theta^{*}} \) has an arbitrary distribution. By independence, the variance of the asymptotic distribution is simply the sum of the variances. As the variance of \( \Delta_{\theta^{*}} \) cannot be negative, we know that the asymptotic variance of any regular estimator \( T_{n} \) is bounded from below, that is, \( \Var(\Delta_{\theta^{*}}) + I_{X}^{-1}(\theta^{*}) \geq I^{-1}_{X}(\theta^{*}) \). %

The MLE is a regular estimator with \( \Delta_{\theta^{*}} \) equal to the fixed true value \( \theta^{*} \), thus, \( \Var(\Delta_{\theta^{*}}) =0 \). As such, the MLE has an asymptotic variance \( I^{-1}_{X}(\theta^{*}) \) that is equal to the lower bound given above. Hence, amongst the broad class of regular estimators, the MLE performs best. This result was already foreshadowed by \citet{fisher1922mathematical}, though it took another 50 years before this statement was made mathematically rigorous \citep{hajek1970characterization, inagaki1970limiting, lecam1970assumptions, van2002statistical, yang1999conversation}, see also \citet{ghosh1985efficiency} for a beautiful review. 

We stress that the normal approximation to the true sampling distribution only holds when \( n \) is large enough. In practice, \( n \) is relatively small and the replacement of the true sampling distribution by the normal approximation can, thus, lead to confidence intervals and hypothesis tests that perform poorly \citep{brown2001interval}. This can be very detrimental, especially, when we are dealing with hard decisions such as the rejection or non-rejection of a hypothesis. %

A simpler version of the H\'{a}jek-LeCam convolution theorem is known as the Cram\'{e}r-Fr\'{e}chet-Rao information lower bound (\citealp{cramer1946methods}; \citealp{frechet1943extension}; \citealp{rao1945information}), which also holds for finite \( n \). This theorem states that the variance of an unbiased estimator \( T_{n} \) cannot be lower than the inverse Fisher information, that is, \( n \Var(T_{n}) \geq I_{X}^{-1}(\theta^{*}) \). We call an estimator \( T_{n}=t(X^{n}) \) \emph{unbiased} if for every possible true value \( \theta^{*} \) and at each fixed \( n \), its expectation is equal to the true value, that is, \( E(T_{n}) = \theta^{*} \). Hence, this lower bound shows that Fisher information is not only a concept that is useful for large samples. 

Unfortunately, the class of unbiased estimators is rather restrictive (in general, it does not include the MLE) and the lower bound cannot be attained whenever the parameter is of more than one dimensions \citep{wijsman1973attainment}. Consequently, for vector-valued parameters \( \theta \), this information lower bound does not inform us, whether we should stop our search for a better estimator. 

The H\'{a}jek-LeCam convolution theorem implies that for \( n \) large enough the MLE \( \hat{\theta} \) is the best performing statistic. For the MLE to be superior, however, the data do need to be generated as specified by the functional relationship \( f \). In reality, we do not know whether the data are indeed generated as specified by \( f \), which is why we should also try to empirically test such an assumption. For instance, we might believe that the data are normally distributed, while in fact they were generated according to a Cauchy distribution. This incorrect assumption implies that we should use the sample mean, but Example~{\ref{cauchyEx}} showed the futility of such estimator. Model misspecification, in addition to hard decisions based on the normal approximation, might be the main culprit of the crisis of replicability. Hence, more research on the detection of model misspecification is desirable and expected (e.g., \citealp{grunwald2016safe}; \citealp{grunwald2014inconsistency}; \citealp{ommen2016robust}). 

\section{Bayesian use of the Fisher-Rao Metric: The Jeffreys's Prior}
\label{appendixBayes}
We make intuitive that the Jeffreys's prior is a uniform prior on the model \( \Mc_{\Theta} \), i.e., 
\begin{align}
P(m^{*} \in J_{m}) = \frac{1}{V} \int_{J_{m}} 1 \der m_{\theta}(X) = \int_{\theta_{a}}^{\theta_{b}} \sqrt{I_{X}(\theta)} \der \theta, 
\end{align}
where \( J_{m} = \big ( m_{\theta_{a}}(X), m_{\theta_{b}}(X) \big ) \) is an interval of pmfs in model space. To do so, we explain why the differential \( \der m_{\theta}(X) \), a displacement in model space, is converted into \( \sqrt{I_{X}(\theta)} \der \theta \) in parameter space. The elaboration below boils down to an explanation of arc length computations using integration by substitution. 

\subsection{Tangent vectors}
First note that we swapped the area of integration by substituting the interval \( J_{m}= \big ( m_{\theta_{a}}(X), m_{\theta_{b}}(X) \big ) \) consisting of pmfs in function space \( \Mc_{\Theta} \) by the interval \( (\theta_{a}, \theta_{b}) \) in parameter space. This is made possible by the parameter functional \( \nu \) with domain \( \Mc_{\Theta} \) and range \( \Theta \) that uniquely assigns to any (transformed) pmf \( m_{a}(X) \in \Mc_{\Theta} \) a parameter value \( \theta_{a} \in \Theta \). In this case, we have \( \theta_{a}=\nu(m_{a}(X))=( \tfrac{1}{2} m_{a}(1))^{2} \). Uniqueness of the assignment implies that the resulting parameter values \( \theta_{a} \) and \( \theta_{b} \) in \( \Theta \) differ from each other whenever \( m_{a}(X) \) and \( m_{b}(X) \) in \( \Mc_{\Theta} \) differ from each other. %
\begin{figure}[h]
   \centering
\begin{tabular}{cc}
    \begin{minipage}{.4 \textwidth}
    \centering
    \includegraphics[width= \linewidth]{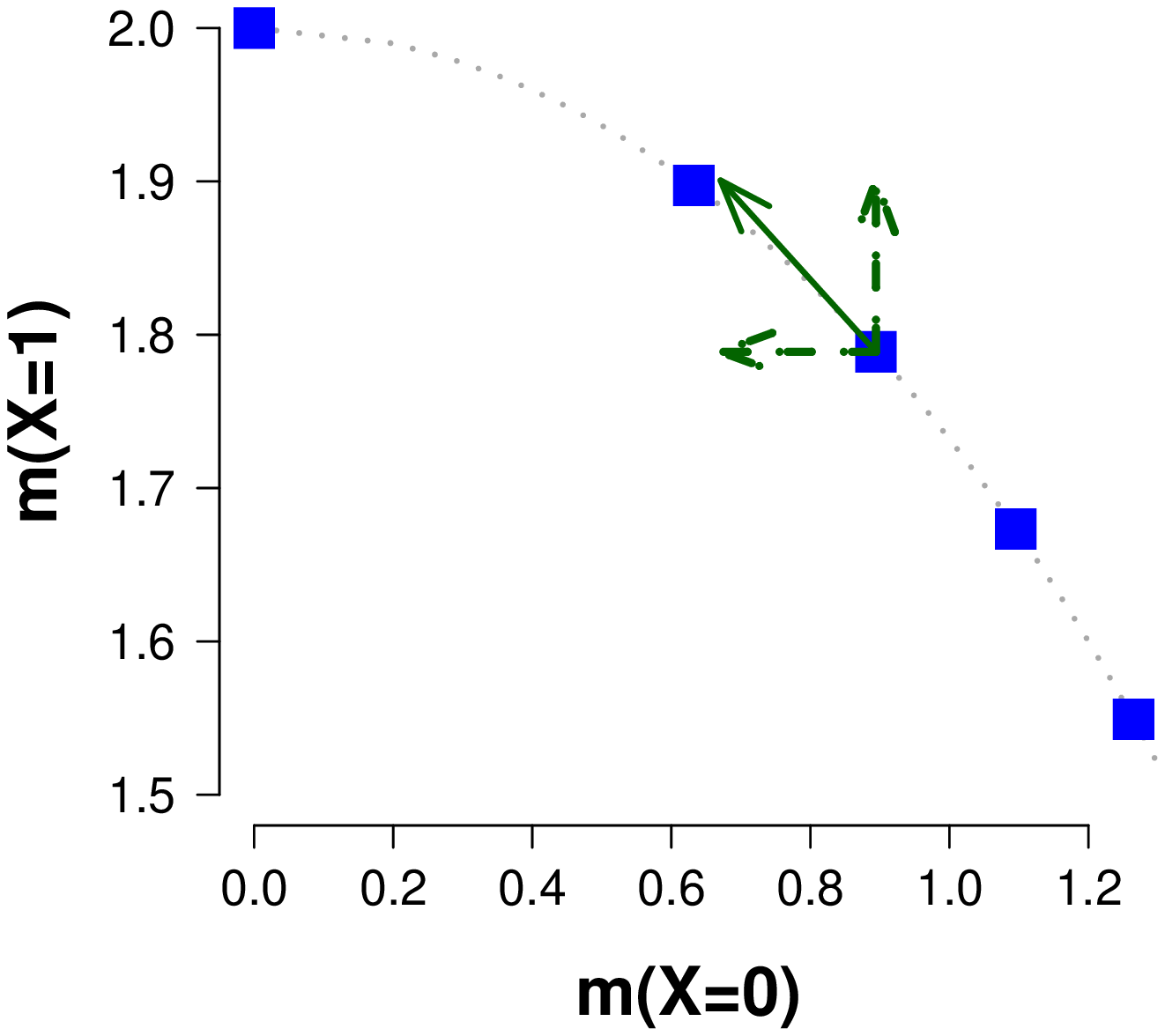}
    \end{minipage}  & %
    \begin{minipage}{.4 \textwidth}
    \centering
    \includegraphics[width= \linewidth]{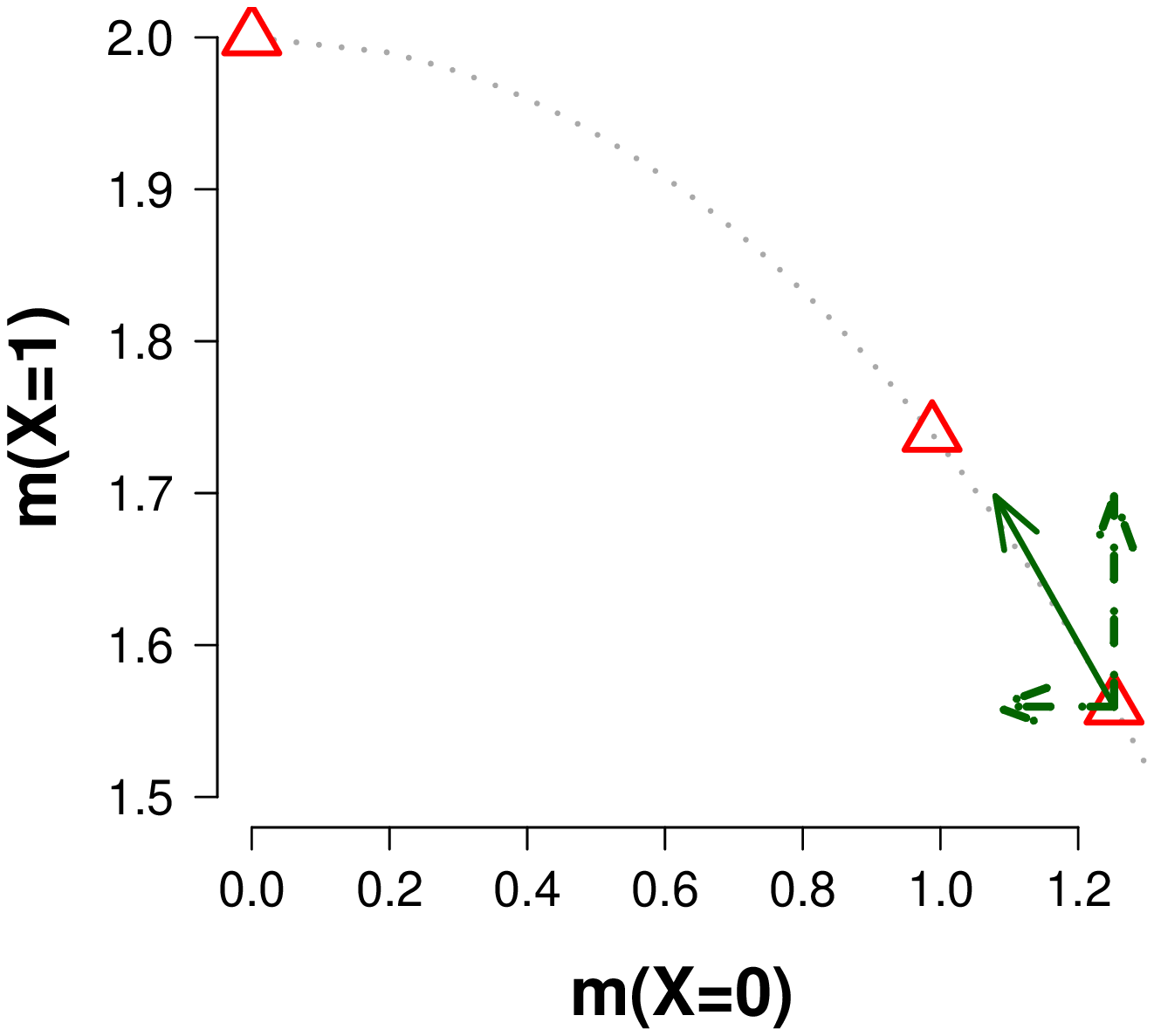}
    \end{minipage}  %
\end{tabular}    
\caption{The full arrow represents the simultaneous displacement in model space based on the Taylor approximation \refEq{l2Error} in terms of \( \theta \) at \( m_{\theta_{a}}(X) \), where \( \theta_{a}=0.8 \) (left panel) and in terms of \( \phi \) at \( m_{\phi_{a}}(X) \) where \( \phi_{a}=0.6 \pi \) (right panel). The dotted line represents a part of the Bernoulli model and note that the full arrow is tangent to the model.}
\label{2dGeoC}
\end{figure}  
For example, the map \( \nu : \Mc_{\Theta} \rightarrow \Theta \) implies that in the left panel of \refFig{2dGeoC} the third square from the left with coordinates \( m_{a}(X)=[0.89, 1.79] \) can be labeled by \( \theta_{a}=0.8 \approx ( \tfrac{1}{2} (1.79) )^{2} \), while the second square from the left with coordinates \( m_{b}(X)=[0.63, 1.90] \) can be labeled by \( \theta_{b}=0.9 \approx ( \tfrac{1}{2} (1.90) )^{2} \). 

To calculate the arc length of the curve \( J_{m} \) consisting of functions in \( \Mc_{\Theta} \), we first approximate \( J_{m} \) by a finite sum of tangent vectors, i.e., straight lines. The approximation of the arc length is the sum of the length of these straight lines. The associated approximation error goes to zero, when we increase the number of tangent vectors and change the sum into an integral sign, as in the usual definition of an integral. First we discuss tangent vectors. 

In the left panel in \refFig{2dGeoC}, we depicted the tangent vector at \( m_{\theta_{a}}(X) \) as the full arrow. This full arrow is constructed from its components: one broken arrow that is parallel to the horizontal axis associated with the outcome \( x=0 \), and one broken arrow that is parallel to the vertical axis associated with the outcome \( x=1 \). The arrows parallel to the axes are derived by first fixing \( X=x \) followed by a Taylor expansion of the parameterization \( \theta \mapsto m_{\theta}(x) \) at \( \theta_{a} \). The Taylor expansion is derived by differentiating with respect to \( \theta \) at \( \theta_{a} \) yielding the following ``linear'' function of the distance \( \der \theta = | \theta_{b} - \theta_{a} | \) in parameter space, 
\begin{align}
\label{l2Error0}
\der m_{\theta_{a}}(x) = m_{\theta_{b}}(x) - m_{\theta_{a}}(x) = \underbrace{\frac{\der m_{\theta_{a}}(x)}{\der \theta}}_{A_{\theta_{a}}(x)} \, \der \theta + \underbrace{o( \der \theta )}_{B_{\theta_{a}}(x)},
\end{align}
where the slope, a function of \( x \), \( A_{\theta_{a}}(x) \) at \( m_{\theta_{a}}(x) \) in the direction of \( x \) is given by %
\begin{align}
\label{l2Error}
A_{\theta_{a}}(x) =\frac{\der m_{\theta_{a}}(x) }{\der \theta} =\tfrac{1}{2} \big \{ \underbrace{\tfrac{\der}{ \der \theta} \log f(x \, | \, \theta_{a})}_{\text{score function}} \big \} m_{\theta_{a}}(x) ,
\end{align}
and with an ``intercept'' \( B_{\theta_{a}}(x)=o ( \der \theta) \) that goes fast to zero whenever \( \der \theta \rightarrow 0 \). Thus, for \( \der \theta \) small, the intercept \( B_{\theta_{a}}(x) \) is practically zero. Hence, we approximate the displacement between \( m_{\theta_{a}}(x) \) and \( m_{\theta_{b}}(x) \) by a straight line.

\begin{voorbeeld}[Tangent vectors]
\label{tangentVectors}
In the right panel of \refFig{2dGeoC} the right most triangle is given by \( m_{\phi_{a}}(X)=[1.25, 1.56] \), while the triangle in the middle refers to \( m_{\phi_{b}}(X)=[0.99, 1.74] \). Using the functional \( \tilde{\nu} \), i.e., the inverse of the parameterization, \( \phi \mapsto 2 \sqrt{f (x \, | \, \phi)} \), where \( f(x \, | \, \phi)= \big ( \frac{1}{2} + \frac{1}{2} \big ( \tfrac{\phi}{\pi} \big )^{3} \big )^{x} \big ( \frac{1}{2} - \frac{1}{2} \big ( \tfrac{\phi}{\pi} \big )^{3} \big )^{1-x} \), we find that these two pmfs correspond to \( \phi_{a}=0.6 \pi \) and \( \phi_{b}=0.8 \pi \).

The tangent vector at \( m_{\phi_{a}}(X) \) is constructed from its components. For the horizontal displacement, we fill in \( x=0 \) in \( \log f(x \, | \, \phi) \) followed by the derivation with respect to \( \phi \) at \( \phi_{a} \) and a multiplication by \( m_{\phi_{a}}(x) \) resulting in %
\begin{align}
\frac{\der m_{\phi{a}}(0)}{\der \phi} \der \phi = & \tfrac{1}{2} \Big \{ \tfrac{\der}{ \der \phi} \log f(0 \, | \, \phi_{a}) \Big \} m_{\phi_{a}}(0) \, \der \phi, \\
= & - {3 \phi_{a}^{2} \over \sqrt{2 \pi^{3} (\pi^{3} + \phi_{a}^{3}})} \, \der \phi ,
\end{align}
where \( \der \phi= | \phi_{b} - \phi_{a} | \) is the distance in parameter space \( \Phi \). The minus sign indicates that the displacement along the horizontal axis is from right to left. Filling in \( \der \phi = | \phi_{b}-\phi_{a} | =0.2 \pi \) and \( \phi_{a}=0.6 \pi \) yields a horizontal displacement of \( 0.17 \) at \( m_{\phi_{a}} (0 ) \) from right to left in model space. Similarly, the vertical displacement in terms of \( \phi \) is calculated by first filling in \( x=1 \) and leads to 
\begin{align}
 \frac{\der m_{\phi_{a}}(1)}{\der \phi} \, \der \phi = & \tfrac{1}{2} \Big \{ \tfrac{\der}{\der \phi} \log f(1 \, | \, \phi_{a}) \Big \} m_{\phi_{a}}(1) \, \der \phi, \\
= & {3 \phi_{a}^{2} \over \sqrt{2 \pi^{3} (\pi^{3} - \phi_{a}^{3}})} \, \der \phi.
\end{align}
By filling in \( \der \phi =0.2 \) and \( \phi_{a}=0.6 \pi \), we see that a change of \( \der \phi= 0.2 \pi \) at \( \phi_{a} = 0.6 \pi \) in the parameter space corresponds to a vertical displacement of \( 0.14 \) at \( m_{\phi_{a}} (1) \) from bottom to top in model space. Note that the axes in \refFig{2dGeoC} are scaled differently. 

The combined displacement \( \frac{\der m_{\phi_{a}}( X)}{\der \phi} \der \phi \) at \( m_{\phi_{a}}(X) \) is the sum of the two broken arrows and plotted as a full arrow in the right panel of \refFig{2dGeoC}. \( \hfill \diamond \)
\end{voorbeeld}
The length of the tangent vector \( \frac{\der m_{\theta_{a}}(X)}{ \der \theta} \) at the vector \( m_{\theta_{a}}(X) \) is calculated by taking the root of the sum of its squared component, the natural measure of distance we adopted above and this yields %
\begin{align}
\label{arcLength}
\Big \| \tfrac{\der m_{\theta_{a}}(X)}{\der \theta} \der  \theta \Big \|_{2} & = \sqrt{\sum_{x \in \Xc} \Big ( \tfrac{\der m_{\theta_{a}}(x) }{\der \theta} \Big)^{2} (\der  \theta)^{2} }, \\
&= \sqrt{\sum_{x \in \Xc} \Big ( \tfrac{\der }{\der \theta} \log f (x \, | \, \theta_{a}) \Big )^{2} p_{\theta_{a}} (x)} \der  \theta 
& = \sqrt{I_{X}(\theta_{a})} \der  \theta .
\end{align}
The second equality follows from the definition of \( {\der m_{\theta_{a}}(X) \over \der \theta} \), i.e., \refEq{l2Error}, and the last equality is due to the definition of Fisher information. 

\begin{voorbeeld}[name=Length of the tangent vectors, continue=tangentVectors]
The length of the tangent vector in the right panel of \refFig{2dGeoC} can be calculated as the root of the sums of squares of its components, that is, \( \| \tfrac{\der m_{\phi_{a}}(X)}{\der \phi} \der  \phi \|_{2} = \sqrt{(-0.14)^2+0.17^2} = 0.22 \). Alternatively, we can first calculate the square root of the Fisher information at \( \phi_{a}=0.6 \pi \), i.e., 
\begin{align}
\sqrt{I(\phi_{a})} = \frac{3 \phi_{a}^{2}}{\sqrt{ \pi^{6} - \phi^{6}}} = 0.35,
\end{align}
and a multiplication with \( \der \phi = 0.2 \pi \) results in \( \| \tfrac{\der m_{\phi_{a}}(X)}{\der \phi} \|_{2} \der  \phi = 0.22 \). \( \hfill \diamond \)
\end{voorbeeld}

More generally, to approximate the length between pmfs \( m_{\theta_{a}}(X) \) and \( m_{\theta_{b}}(X) \), we first identify \( \nu(m_{\theta_{a}}(X))=\theta_{a} \) and multiply this with the distance \( \der \theta = | \theta_{a} - \nu(m_{\theta_{b}}(X))| \) in parameter space, i.e., %
\begin{align}
\label{Jacobian}
\der m_{\theta} (X) = \Big \| {\der m_{\theta}(X) \over \der \theta} \Big \|_{2} \, \der \theta = \sqrt{I_{X}(\theta)} \, \der \theta.
\end{align}
In other words, the root of the Fisher information converts a small distance \( \der \theta \) at \( \theta_{a} \) to a displacement in model space at \( m_{\theta_{a}}(X) \). 

\subsection{The Fisher-Rao metric}
By virtue of the parameter functional \( \nu \), we send an interval of pmfs \( J_{m}=\Big ( m_{\theta_{a}}(X) , m_{\theta_{b}}(X) \Big ) \) in the function space \( \Mc_{\Theta} \) to the interval \( (\theta_{a}, \theta_{b}) \) in the parameter space \( \Theta \). In addition, with the conversion of \( \der m_{\theta}(X) = \sqrt{I_{X}(\theta)} \, \der \theta \) we integrate by substitution, that is, %
\begin{align}
P \Big ( m^{*}(X) \in J_{m} \Big )= \frac{1}{V} \int_{m_{\theta_{a}}(X)}^{m_{\theta_{b}}(X)} 1 \der m_{\theta}(X) = \frac{1}{V} \int_{\theta_{a}}^{\theta_{b}} \sqrt{I_{X}(\theta)} \der \theta.
\end{align}
In particular, choosing \( J_{\theta}=\Mc_{\Theta} \) yields the normalizing constant \( V = \int_{0}^{1} \sqrt{I_{X}(\theta)} \der \theta \). The interpretation of \( V \) as being the total length of \( \Mc_{\Theta} \) is due to the use of \( \der m_{\theta}(X) \) as the metric, a measure of distance, in model space. To honour Calyampudi Radhakrishna Rao's \citeyearpar{rao1945information} contribution to the theory, this metric is also known as the Fisher-Rao metric (e.g., \citealp{amari1987differentialBook}; \citealp{atkinson1981rao}; \citealp{burbea1984informative}; \citealp{burbea1982entropy, burbea1984differential}; \citealp{dawid1977further}; \citealp{efron1975defining}; \citealp{kass2011geometrical}). 

\subsection{Fisher-Rao metric for vector-valued parameters}
\label{fisherRaoMetric}
\subsubsection{The parameter functional \( \nu: \Pc \rightarrow B\) and the categorical distribution}
For random variables with \( w \) number of outcomes, the largest set of pmfs \( \Pc \) is the collection of functions \( p \) on \( \Xc \) such that (i) \( 0 \leq p(x)=P(X=x) \) for every outcome \( x \) in \( \Xc \), and (ii) to explicitly convey that there are \( w \) outcomes, and none more, these \( w \) chances have to sum to one, that is, \( \sum_{x \in \Xc} p(x)=1 \). The complete set of pmfs \( \Pc \) can be parameterized using the functional \( \nu \) that assigns to each \( w \)-dimensional pmf \( p(X) \) a parameter \( \beta \in \R^{w-1} \). 

For instance, given a pmf \( p(X)=[p(L), p(M), p(R)] \) we typically use the functional \( \nu: \Pc \rightarrow \R^{2} \) that takes the first two coordinates, that is, \( \nu(p(X)) = \beta = {\beta_{1} \choose \beta_{2}} \), where \( \beta_{1}=p(L) \) and \( \beta_{2}=p(M) \). The range of this functional \( \nu \) is the parameter space \( B=[0, 1] \times [0, \beta_{1}] \). Conversely, the inverse of the functional \( \nu \) is the parameterization \( \beta \mapsto p_{\beta}(X)=[ \beta_{1}, \beta_{2}, 1- \beta_{1} - \beta_{2}] \), where (i') \( 0 \leq \beta_{1}, \beta_{2} \) and (ii') \( \beta_{1}+\beta_{2} \leq 1 \). The restrictions (i') and (ii') imply that the parameterization has domain \( B \) and the largest set of pmfs \( \Pc \) as its range. By virtue of the functional \( \nu \) and its inverse, that is, the parameterization \( \beta \mapsto p_{\beta}(X) \), we conclude that the parameter space \( B \) and the complete set of pmfs \( \Pc \) are isomorphic. This means that each pmf \( p(X) \in \Pc \) can be uniquely identified with a parameter \( \beta \in B \) and vice versa. The inverse of \( \nu \) implies that the parameters \( \beta \in B \) are functionally related to the potential outcomes \( x \) of \( X \) as %
\begin{align}
\label{categoricalPmf}
f(x \, | \, \beta) = \beta_{1}^{x_{L}} \beta_{2}^{x_{M}} (1- \beta_{1} - \beta_{2})^{x_{R}},
\end{align}
where \( x_{L}, x_{M} \) and \( x_{R} \) are the number of \( L, M \) and \( R \) responses in one trial --we either have \( x=[x_{L}, x_{M}, x_{R}]=[1, 0,0] \), \( x=[0, 1, 0] \), or \( x=[0, 0, 1] \). The model \( f(x \, | \, \beta) \) can be regarded as the generalization of the Bernoulli model to \( w=3 \) categories. In effect, the parameters \( \beta_{1} \) and \( \beta_{2} \) can be interpreted as a participant's propensity of choosing \( L \) and \( M \), respectively. If \( X^{n} \) consists of \( n \) iid categorical random variables with the outcomes \( [L, M, R] \), the joint pmf of \( X^{n} \) is then%
\begin{align}
\label{trinomial}
f(x^{n} \, | \, \beta)= \beta_{1}^{y_{L}} \beta_{2}^{y_{M}} (1 - \beta_{1} - \beta_{2})^{y_{R}}, 
\end{align}
where \( y_{L}, y_{M} \) and \( y_{R}=n-y_{L}-y_{M} \) are the number of \( L, M \) and \( R \) responses in \( n \) trials. As before, the representation of the pmfs as the vectors \( m_{\beta}(X)=[2 \sqrt{\beta_{1}}, 2 \sqrt{ \beta_{2}}, 2 \sqrt{1-\beta_{1} - \beta_{2}}] \) form the surface of (the positive part of) the sphere of radius two, thus, \( \Mc = \Mc_{B} \), see \refFig{fiOrthogonal}. The extreme pmfs indicated by \( mL, mM \) and \( mR \) in the figure are indexed by the parameter values \( \beta=(1, 0) \), \( \beta=(0, 1) \) and \( \beta=(0, 0) \), respectively. 

\subsubsection{The stick-breaking parameterization of the categorical distribution}
Alternatively, we could also have used a ``stick-breaking'' parameter functional \( \tilde{\nu} \) that sends each pmf in \( \Pc \) to the vector of parameters \( \tilde{\nu}(p(X)) = {\gamma_{1} \choose \gamma_{2}} \), where \( \gamma_{1}=p_{L} \) and \( \gamma_{2}= p_{M}/(1- p_{L}) \).%
\footnote{This only works if \( p_{L} < 1 \). When \( p(x_{1})=1 \), we simply set \( \gamma_{2}=0 \), thus, \( \gamma=(1, 0) \).} %
Again the parameter \( \gamma={\gamma_{1} \choose \gamma_{2}} \) is only a label, but this time the range of \( \tilde{\nu} \) is the parameter space \( \Gamma=[0,1] \times [0,1] \). The functional relationship \( f \) associated to \( \gamma \) is given by 
\begin{align}
f(x \, | \, \gamma) = \gamma_{1}^{x_{L}} \big ( (1- \gamma_{1}) \gamma_{2} \big )^{x_{M}} \big ( (1- \gamma_{1}) (1- \gamma_{2}) \big )^{x_{R}}.
\end{align}
For each \( \gamma \) we can transform the pmf into the vector 
\begin{align}
m_{\gamma}(X) = [2 \sqrt{\gamma_{1}}, 2 \sqrt{(1- \gamma_{1}) \gamma_{2}} , 2 \sqrt{(1- \gamma_{1})(1-\gamma_{2})}] ,
\end{align}
and write \( \Mc_{\Gamma} \) for the collection of vectors so defined. As before, this collection coincides with the full model, i.e., \( \Mc_{\Gamma}=\Mc \). In other words, by virtue of the functional \( \tilde{\nu} \) and its inverse \( \gamma \mapsto p_{\gamma}(x) =f(x \, | \, \gamma) \) we conclude that the parameter space \( \Gamma \) and the complete set of pmfs \( \Mc \) are isomorphic. Because \( \Mc=\Mc_{B} \) this means that we also have an isomorphism between the parameter space \( B \) and \( \Gamma \) via \( \Mc \), even though \( B \) is a strict subset of \( \Gamma \). Note that this equivalence goes via parameterization \( \beta \mapsto m_{\beta}(X)  \) and the functional \( \tilde{\nu} \). 

\subsubsection{Multidimensional Jeffreys's prior via the Fisher information matrix and orthogonal parameters}
\label{appendixOrtho}
The multidimensional Jeffreys's prior is parameterization-invariant and has as normalization constant \( V = \int \sqrt{\det  I_{X}(\theta) } \der \theta \), where \( \det  I_{X}(\theta)  \) is the determinant of the Fisher information matrix. 

In the previous subsection we argued that the categorical distribution in terms of \( \beta \) or parameterized with \( \gamma \) are equivalent to each other, that is, \( \Mc_{B} = \Mc=\Mc_{\Gamma} \). However, these two parameterizations describe the model space \( \Mc \) quite differently. In this subsection we use the Fisher information to show that the parameterization in terms of \( \gamma \) is sometimes preferred over \( \beta \). 

The complete model \( \Mc \) is easier described by \( \gamma \), because the parameters are orthogonal. We say that two parameters are \emph{orthogonal to each other} whenever the corresponding off-diagonal entries in the Fisher information matrix are zero. The Fisher information matrices in terms of \( \beta \) and \( \gamma \) are 
\begin{align}
I_{X}(\beta) = \frac{1}{1- \beta_{1} - \beta_{2}} \begin{pmatrix}
1- \beta_{2} & 1 \\
1 & 1- \beta_{1}
\end{pmatrix}
\, \text{ and } \, 
I_{X}(\gamma) = \begin{pmatrix}
\frac{1}{\gamma_{1}(1-\gamma_{1})} & 0 \\
0 & \frac{1-\gamma_{1}}{\gamma_{2}(1-\gamma_{2})}
\end{pmatrix},
\end{align}
respectively. %
\begin{figure}[h]
\centering
\includegraphics[width = 1 \textwidth]{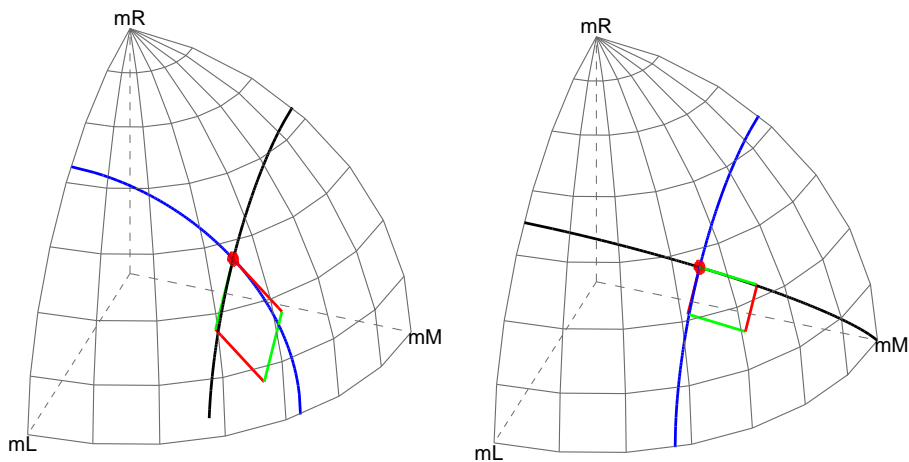}
\caption{When the off-diagonal entries are zero, the tangent vectors are orthogonal. Left panel: The tangent vectors at \( p_{\beta^{*}}(X)=[1/3,1/3,1/3] \) span a diamond with an area given by \( \sqrt{\det I(\beta^{*})} \der  \beta \). The black curve is the submodel with \( \beta_{2}=1/3 \) fixed and \( \beta_{1} \) free to vary and yields a green tangent vector. The blue curve is the submodel with \( \beta_{1}=1/3 \) fixed and \( \beta_{2} \) free to vary. Right panel: The tangent vectors at the same pmf in terms of \( \gamma \), thus, \( p_{\gamma^{*}}(X) \), span a rectangle with an area given by \( \sqrt{ \det I(\gamma^{*})} \der  \gamma \). The black curve is the submodel with \( \gamma_{2}=1/2 \) fixed and \( \gamma_{1} \) free to vary and yields a green tangent vector. The blue curve is the submodel with \( \gamma_{1}=1/3 \) fixed and \( \gamma_{2} \) free to vary.}
\label{fiOrthogonal}
\end{figure}%
The left panel of \refFig{fiOrthogonal} shows the tangent vectors  at \( p_{\beta^{*}}(X)=[1/3, 1/3, 1/3] \) in model space, where \( \beta^{*}=(1/3, 1/3) \). The green tangent vector corresponds to \( \frac{ \partial m_{\beta^{*}}(X)}{ \partial \beta_{1}} \), thus, with \( \beta_{2}=1/3 \) fixed and \( \beta_{1} \) free to vary, while the red tangent vector corresponds to \( \frac{ \partial m_{\beta^{*}}(X)}{ \partial \beta_{2}} \), thus, with \( \beta_{1}=1/3 \) and \( \beta_{2} \) free to vary. The area of the diamond spanned by these two tangent vectors is \( \sqrt{ \det I(\beta^{*})} \der \beta_{1} \der \beta_{2} \), where we have taken \( \der \beta_{1}=0.1 \) and \( \der \beta_{2}=0.1 \). 

The right panel of \refFig{fiOrthogonal} shows the tangent vectors at the same point \( p_{\gamma^{*}}(X)=[1/3, 1/3, 1/3] \), where \( \gamma^{*}=(1/3, 1/2) \). The green tangent vector corresponds to \( \frac{ \partial m_{\gamma^{*}}(X)}{ \partial \gamma_{1}} \), thus, with \( \gamma_{2}=1/2 \) fixed and \( \gamma_{1} \) free to vary, while the red tangent vector corresponds to \( \frac{ \partial m_{\gamma^{*}}(X)}{ \partial \gamma_{2}} \), thus, with \( \gamma_{1}=1/3 \) and \( \gamma_{2} \) free to vary. By glancing over the plots, we see that the two tangent vectors are indeed orthogonal. The area of the rectangle spanned by these these two tangent vectors is \( \sqrt{ \det I(\gamma^{*})} \der \gamma_{1} \der \gamma_{2} \), where we have taken \( \der \gamma_{1} = \der \gamma_{2}=0.1 \). 

There are now two ways to calculate the normalizing constant of the Jeffreys's prior, the area, more generally volume, of the model \( \Mc \). In terms of \( \beta \) this leads to 
\begin{align}
V = \int_{0}^{1} \left ( \int_{0}^{\beta_{1}} \frac{1}{1 - \beta_{1} - \beta_{2}} \sqrt{ \beta_{1} \beta_{2} - \beta_{1} - \beta_{2} } \der \beta_{2} \right ) \der \beta_{1}.
\end{align}
Observe that the inner integral depends on the value of \( \beta_{1} \) from the outer integral. This coupling is reflected by the non-zero off-diagonal term of the Fisher information matrix \( I_{X}(\beta) \) corresponding to \( \beta_{1} \) and \( \beta_{2} \). On the other hand, orthogonality implies that the two parameters can be treated independently of each other. That is, knowing and fixing \( \gamma_{1} \) and changing \( \gamma_{2} \) will not affect \( m_{\gamma}(X) \) via \( \gamma_{1} \). This means that the double integral decouples %
\begin{align}
V = \int_{0}^{1} \left ( \int_{0}^{1} \frac{1}{\sqrt{\gamma_{1} \gamma_{2} ( 1 - \gamma_{2})}} \der \gamma_{1} \right ) \der \gamma_{2} = \int_{0}^{1} \frac{1}{\sqrt{\gamma_{1}}} \der \gamma_{1} \int_{0}^{1} \frac{1}{\sqrt{\gamma_{2} ( 1 - \gamma_{2})}} \der \gamma_{2} = 2 \pi.
\end{align}
Using standard geometry we verify that this is indeed the area of \( \Mc \), as an eighth of the surface area of a sphere of radius two is given by \( \tfrac{1}{8} 4 \pi 2^{2} = 2 \pi \). 

Orthogonality is relevant in Bayesian analysis, as it provides an argument to choose a prior on a vector-valued parameter that factorizes (e.g., \citealp{berger1998bayes}; \citealp{huzurbazar1950probability, huzurbazar1956sufficient}; \citealp{jeffreys1961theory}; \citealp{kass1992approximate}; \citealp{ly2016harold, ly2016evaluation}), see also \citet{cox1987parameter, mitchell1962sufficient}.

By taking a random variable \( X \) with \( w=3 \) outcomes, we were able to visualize the geometry of model space. For more general \( X \) these plots get more complicated and perhaps even impossible to draw. Nonetheless, the ideas conveyed here extend, even to continuous \( X \), whenever the model adheres to the regularity conditions given in Appendix~{\ref{appendixReg}}. 

\section{MDL: Coding Theoretical Background}
\label{appendixMdl}
\subsection{Coding theory, code length and log-loss}
A coding system translates words, i.e., outcomes of a random variable \( X \), into code words with code lengths that behave like a pmf. Code lengths can be measured with a logarithm, which motivates the adoption of log-loss, defined below, as \emph{the} decision criterion within the MDL paradigm. The coding theoretical terminologies introduced here are illustrated using the random variable \( X \) with \( w=3 \) potential outcomes. 

\subsubsection{Kraft-McMillan inequality: From code lengths of a specific coding system to a pmf}
For the source-memory task we encoded the outcomes as \( L, M \) and \( R \), but when we communicate a participant's responses \( x^{n}_{\obs} \) to a collaborator over the internet, we have to encode the observations \( x^{n}_{\obs} \) as zeroes and ones. For instance, we might use a coding system \( \tilde{C} \) with code words \( \tilde{C}(X=L)=00 \), \( \tilde{C}(X=M)=01 \) and \( \tilde{C}(X=R)=10 \). This coding system \( \tilde{C} \) will transform any set of responses \( x^{n}_{\obs} \) into a code string \( \tilde{C}(x^{n}_{\obs}) \) consisting of \( 2 n \) bits. Alternatively, we can use a coding system \( C \) with code words \( C(X=L)=10 \), \( C(X=M)=0 \) and \( C(X=R)=11 \), instead. Depending on the actual observations \( x^{n}_{\obs} \), this coding system outputs code strings \( C(x^{n}_{\obs}) \) with varying code lengths that range from \( n \) to \( 2 n \) bits. For example, if a participant responded with \( x^{n}_{\obs}=(M, R, M, L, L, M, M, M) \) in \( n=8 \) trials, the coding system \( C \) would then output the 11-bit long code string \( C(x^{n}_{\obs})=01101010000 \). In contrast, the first coding system \( \tilde{C} \) will always output a \( 16 \)-bit long code string when \( n=8 \). Shorter code strings are desirable as they will lead to a smaller load on the communication network and they are less likely to be intercepted by ``competing'' researchers.

Note that the shorter code length \( C(x^{n}_{\obs})=01101010000 \) of 11-bits is a result of having code words of unequal lengths. The fact that one of the code word is shorter does not interfere with the decoding, since no code word is a prefix of another code word. As such, we refer to \( C \) as a prefix (free) coding system. This implies that the 11-bit long code string \( C(x^{n}_{\obs}) \) is self-punctuated and that it can be uniquely deciphered by simply reading the code string from left to right resulting in the retrieval of \( x^{n}_{\obs} \). Note that the code lengths of \( C \) inherit the randomness of the data. In particular, the coding system \( C \) produces a shorter code string with high chance, if the participant generates the outcome \( M \) with high chance. In the extreme case, the coding system \( C \) produces the 8-bits long code string \( C(x^{n})=00000000 \) with 100\% (respectively, 0\%) chance, if the participant generates the outcome \( M \) with 100\% (respectively, 0\%) chance. 
More generally, Kraft and McMillan (\citealp{kraft1949device}; \citealp{mcmillan1956two}) showed that \emph{any} uniquely decipherable (prefix) coding system from the outcome space \( \Xc \) with \( w \) outcomes to an alphabet with \( D \) elements must satisfy the inequality 
\begin{align}
\sum_{i=1}^{w} D^{-l_{i}} \leq 1,
\end{align}
where \( l_{i} \) is the code length of the outcome \( w \). In our example, we have taken \( D=2 \) and code length of 2, 1 and 2 bits for the response \( L, M \) and \( R \) respectively. Indeed, \( 2^{-2}+2^{-1}+2^{-2}=1 \). Hence, code lengths behave like the logarithm (with base \( D \)) of a pmf. 

\subsubsection{Shannon-Fano algorithm: From a pmf to a coding system with specific code lengths}
Given a data generating pmf \( p^{*}(X) \), we can use the so-called \emph{Shannon-Fano algorithm} (e.g., \citealp[Ch.~{5}]{cover2006elements}) to construct a prefix coding system \( C^{*} \). The idea behind this algorithm is to give the outcome \( x \) that is generated with the highest chance the shortest code length. To do so, we encode the outcome \( x \) as a code word \( C^{*}(x) \) that consists of \( - \log_{2} p^{*}(x) \) bits.%
\footnote{When we use the logarithm with base two, \( \log_{2}(y) \), we get the code length in bits, while the natural logarithm, \( \log(y) \), yields the code length in nats. Any result in terms of the natural logarithm can be equivalently described in terms of the logarithm with base two, as \( \log(y) = \log(2) \log_{2}(y) \).} %

For instance, when a participant generates the outcomes \( [L, M, R] \) according to the chances \( p^{*}(X)=[0.25, 0.5, 0.25] \) the Shannon-Fano algorithm prescribes that we should encode the outcome \( L \) with \( - \log_{2} (0.25) = 2 \), \( M \) with \( - \log_{2}(0.5)=1 \) and \( R \) with 2 bits; the coding system \( C \) given above.%
\footnote{Due to rounding, the Shannon-Fano algorithm actually produces code words \( C(x) \) that are at most one bit larger than the ideal code length \( - \log_{2} p^{*}(x) \). We avoid further discussions on rounding. Moreover, in the following we consider the natural logarithm instead.} %
The Shannon-Fano algorithm works similarly for any other given pmf \( p_{\beta}(X) \). Hence, the Kraft-McMillan inequality and its inverse, i.e., the Shannon-Fano algorithm imply that pmfs and uniquely decipherable coding systems are equivalent to each other. As such we have an additional interpretation of a pmf. To distinguish the different uses, we write \( f(X \, | \, \beta) \) when we view the pmf as a coding system, while we retain the notation \( p_{\beta}(X) \) when we view the pmf as a data generating device. In the remainder of this section we will not explicitly construct any other coding system, as the coding system itself is irrelevant for the discussion at hand --only the code lengths matter.

\subsubsection{Entropy, cross entropy, log-loss}
With the true data generating pmf \( p^{*}(X) \) at hand, thus, also the true coding system \( f(X \, | \, \beta^{*}) \), we can calculate the (population) average code length per trial
\begin{align}
\label{shannonEntropy}
H(p^{*}(X))=H \Big ( p^{*}(X) \| \, f(X \, | \, \beta^{*}) \Big ) = \sum_{x \in \Xc} - \log f(x \, | \, \beta^{*}) p^{*}(x).
\end{align}
Whenever we use the logarithm with base 2, we refer to this quantity \( H(p^{*}(X)) \) as the Shannon entropy.%
\footnote{Shannon denoted this quantity with an \( H \) to refer to the capital Greek letter for eta. It seems that John von Neumann convinced Claude Shannon to call this quantity entropy rather than information \citep{tribus1971energy}.} %
If the true pmf is \( p^{*}(X)=[0.25, 0.5, 0.25] \) we have an average code length of \( 1.5 \) bits per trail whenever we use the true coding system \( f(X \, | \, \beta^{*}) \). Thus, we expect to use 12 bits to encode observations consisting of \( n=8 \) trials. 

As coding theorists, we have no control over the true data generating pmf \( p^{*}(X) \), but we can choose the coding system \( f(X \, | \, \beta) \) to encode the observations. The (population) average code length per trial is given by
\begin{align}
\label{crossEntropy}
H(p^{*}(X) \| \, \beta ) = H \Big (p^{*}(X) \| \, f(X \, | \, \beta) \Big ) = \sum_{x \in \Xc} - \log f(x \, | \, \beta) p^{*}(x).
\end{align}
The quantity \( H (p^{*}(X) \| \,\beta) \) is also known as the cross entropy from the true pmf \( p^{*}(X) \) to the postulated \( f(X \, | \, \beta) \).%
\footnote{Observe that the entropy \( H(p^{*}(X)) \) is the just the cross entropy from the true \( p^{*}(X) \) to the true coding system \( f(X \, | \, \beta^{*}) \).} %
For instance, when we use the pmf \( f(X \, | \, \beta)=[0.01, 0.18, 0.81] \) to encode data that are generated according to \( p^{*}(X)=[0.25, 0.5, 0.25] \), we will use \( 2.97 \) bits on average per trial. Clearly, this is much more than the \( 1.5 \) bits per trial that we get from using the true coding system \( f(X \, | \, \beta^{*}) \). 

More generally, \citet{shannon1948mathematical} showed that the cross entropy can never be smaller than the entropy, i.e., \( H(p^{*}(X)) \leq H (p^{*}(X) \| \, \beta) \). In other words, we always get a larger average code length, whenever we use the wrong coding system \( f(X \, | \, \beta) \). To see why this holds, we decompose the cross entropy as a sum of the entropy and the Kullback-Leibler divergence,%
\footnote{The KL-divergence is also known as the relative entropy.} %
and show that the latter cannot be negative. This decomposition follows from the definition of cross entropy and a subsequent addition and subtraction of the entropy resulting in
\begin{align}
H(p^{*}(X) \| \, \beta) = H(p^{*}(X)) + \underbrace{\sum_{x \in \Xc} \big ( \log \frac{p^{*}(x)}{ f(x \, | \, \beta^{*})} \big ) p^{*}(x) }_{D(p^{*}(X) \| \beta)},
\end{align}
where \( D(p^{*}(X) \| \beta) \) defines the \emph{Kullback-Leibler divergence} from the true pmf \( p^{*}(X) \) to the postulated coding system \( f(X \, | \, \beta) \). Using the so-called \emph{Jensen's inequality} it can be shown that the KL-divergence is non-negative and that it is only zero whenever \( f(X \, | \, \beta) =p^{*}(X) \). Thus, the cross entropy can never be smaller than the entropy. Consequently, to minimize the load on the communication network, we have to minimize the cross entropy with respect to the parameter \( \beta \). Unfortunately, however, we cannot do this in practice, because the cross entropy is a population quantity based on the unknown true pmf \( p^{*}(X) \). Instead, we do the next best thing by replacing the true \( p^{*}(X) \) in \refEq{crossEntropy} by the empirical pmf that gives the relative occurrences of the outcomes in the sample rather than in the population. Hence, for any postulated \( f( X \, | \, \beta) \), with \( \beta \) fixed, we approximate the population average defined in \refEq{crossEntropy} by the sample average %
\begin{align}
\label{logLoss}
H(x^{n}_{\obs} \| \,\beta)=H \Big (\hat{p}_{\obs}(X) \| \, f(X \, | \, \beta) \Big ) = \sum_{i=1}^{n} - \log f(x_{\obs, i} \, | \, \beta) = - \log f(x^{n}_{\obs} \, | \, \beta).
\end{align}
We call the quantity \( H(x^{n}_{\obs} \| \, \beta) \) the log-loss from the observed data \( x^{n}_{\obs} \), i.e., the empirical pmf \( \hat{p}_{\obs}(X) \), to the coding system \( f(X \, | \, \beta) \). 

\subsection{Data compression and statistical inference}
The entropy inequality \( H(p^{*}(X)) \leq H( p^{*}(X) \| \beta) \) implies that the coding theorist's goal of finding the coding system \( f(X \, | \, \beta) \) with the shortest average code length is in fact equivalent to the statistical goal of finding the true data generating process \( p^{*}(X) \). The coding theorist's best guess is the coding system \( f(X \, | \, \beta) \) that minimizes the log-loss from \( x^{n}_{\obs} \) to the model \( \Mc_{B} \). Note that minimizing the negative log-likelihood is the same as maximizing the likelihood. Hence, the log-loss is minimized by the coding system associated with the MLE, thus, the predictive pmf \( f(X \, | \, \hat{\beta}_{\obs}) \). Furthermore, the cross entropy decomposition shows that minimization of the log-loss is equivalent to minimization of the KL-divergence from the observations \( x^{n}_{\obs} \) to the model \( \Mc_{B} \). The advantage of having the optimization problem formulated in terms of KL-divergence is that it has a known lower bound, namely, zero. Moreover, whenever the KL-divergence from \( x^{n}_{\obs} \) to the code \( f(X \, | \, \hat{\beta}_{\obs}) \) is larger than zero, we then know that the empirical pmf associated to the observations does not reside on the model. In particular, Section~{\ref{mleProjection}} showed that the MLE plugin, \( f(X \, | \, \hat{\beta}_{\obs}) \) is the pmf on the model that is closest to the data. This geometric interpretation is due to the fact that we retrieve the Fisher-Rao metric, when we take the second derivative of the KL-divergence with respect to \( \beta \) \citep{kullback1951information}. This connection between the KL-divergence and Fisher information is exploited in \citet{ghosal1997non} to generalize the Jeffreys's prior to nonparametric models, see also \citet{van2014renyi} for the relationship between KL-divergence and the broader class of divergence measures developed by \citet{renyi1961measures}, see also \citet{campbell1965coding}.

\section{Regularity conditions}
\label{appendixReg}
A more mathematically rigorous exposition of the subject would have had this section as the starting point, rather than the last section of the appendix. %
%
%
The regularity conditions given below can be seen as a summary, and guidelines for model builders. If we as scientists construct models such that these conditions are met, we can then use the results presented in the main text. We first give a more general notion of statistical models, then state the regularity conditions followed by a brief discussion on these conditions. 

The goal of statistical inference is to find the true probability measure \( P^{*} \) that governs the chances with which \( X \) takes on its events. A model \( \Ps_{\Theta} \) defines a subset of \( \Ps \), the largest collection of all possible probability measures. We as model builders choose \( \Ps_{\Theta} \) and perceive each probability measure \( P \) within \( \Ps_{\Theta} \) as a possible explanation of how the events of \( X \) were or will be generated. When \( P^{*} \in \Ps_{\Theta} \) we have a well-specified model and when \( P^{*} \nin \Ps_{\Theta} \), we say that the model is misspecified. 

By taking \( \Ps_{\Theta} \) to be equal to the largest possible collection \( \Ps \), we will not be misspecified. Unfortunately, this choice is not helpful as the complete set is hard to track and leads to uninterpretable inferences. Instead, we typically construct the candidate set \( \Ps_{\Theta} \) using a parameterization that sends a label \( \theta \in \Theta \) to a probability measure \( P_{\theta} \). For instance, we might take the label \( \theta = {\mu \choose \sigma^{2}} \) from the parameter space \( \Theta = \R \times (0, \infty) \) and interpret these two numbers as the population mean and variance of a normal probability \( P_{\theta} \). This distributional choice is typical in psychology, because it allows for very tractable inference with parameters that are generally overinterpreted. Unfortunately, the normal distribution comes with rather stringent assumptions resulting in a high risk of misspecification. More specifically, the normal distribution is far too ideal, as it supposes that the population is nicely symmetrically centred at its population mean and outliers are practically not expected due to its tail behavior. 

Statistical modeling is concerned with the intelligent construction of the candidate set \( \Ps_{\Theta} \) such that it encapsulates the true probability measure \( P^{*} \). In other words, the restriction of \( \Ps \) to \( \Ps_{\Theta} \) in a meaningful manner. Consequently, the goal of statistical inference is to give an informed guess \( \tilde{P} \) within \( \Ps_{\Theta} \) for \( P^{*} \) based on the data. This guess should give us insights to how the data \emph{were} generated and how yet unseen data \emph{will be generated}. Hence, the goal is not to find the parameters as they are mere labels. Of course parameters can be helpful, but they should not be the goal of inference. 

Note that our general description of a model as a candidate set \( \Ps_{\Theta} \) does not involve any structure --thus, the members of \( \Ps_{\Theta} \) do not need to be related to each other in any sense. We use the parameterization to transfer the structure of our labels \( \Theta \) to a structure on \( \Ps_{\Theta} \). To do so, we require that \( \Theta \) is a nice open subset of \( \R^{d} \). Furthermore, we require that each label defines a member \( P_{\theta} \) of \( \Ps_{\Theta} \) unambiguously. This means that if \( \theta^{*} \) and \( \theta \) differ from each other that the resulting pair of probability measure \( P_{\theta^{*}} \) and \( P_{\theta} \) also differ from each other. Equivalently, we call a parameterization identifiable whenever \( \theta^{*}=\theta \) leads to \( P_{\theta^{*}}=P_{\theta} \). Conversely, identifiability implies that when we know everything about \( P_{\theta} \), we can then also use the inverse of the parameterization to pinpoint the unique \( \theta \) that corresponds to \( P_{\theta} \). We write \( \nu: \Pc_{\Theta} \rightarrow \Theta \) for the functional that attaches to each probability measure \( P \) a label \( \theta \). For instance, \( \nu \) could be defined on the family of normal distribution such that \( P \mapsto \nu(P) = { E_{P}(X) \choose \Var_{P}(X)} = {\mu \choose \sigma^{2}} \). In this case we have \( \nu(\Ps_{\Theta}) = \Theta \) and, therefore, a one-to-one correspondence between the probability measures \( P_{\theta} \in \Ps_{\Theta} \) and the parameters \( \theta \in \Theta \). 

By virtue of the parameterization and its inverse \( \nu \), we can now transfer additional structure from \( \Theta \) to \( \Ps_{\Theta} \). We assume that each probability measure \( P_{\theta} \) that is defined on the events of \( X \) can be identified with a probability density function (pdf) \( p_{\theta}(x) \) that is defined on the outcomes of \( X \). For this assumption, we require that the set \( \Ps_{\Theta} \) is dominated by a so-called countably additive measure \( \lambda \). When \( X \) is continuous, we usually take for \( \lambda \) the Lebesgue measure that assigns to each interval of the form \( (a, b) \) a length of \( b - a \). Domination allows us to express the probability of \( X \) falling in the range \( (a, b) \) under \( P_{\theta} \) by the ``area under the curve of \( p_{\theta}(x) \)'', that is, \( P_{\theta} \big (X \in (a, b) \big ) = \int_{a}^{b} p_{\theta}(x) \der x \). For discrete variables \( X \) taking values in \( \Xc =\{ x_{1}, x_{2}, x_{3}, \ldots \} \), we take \( \lambda \) to be the counting measure. Consequently, the probability of observing the event \( X \in A \) where \( A =\{ a=x_{1}, x_{2}, \ldots, b=x_{k} \} \) is calculated by summing the pmf at each outcome, that is, \( P_{\theta}(X \in A) = \sum_{x=a}^{x=b} p_{\theta}(x) \). Thus, we represent \( \Ps_{\Theta} \) as the set \( \Pc_{\Theta} = \{ p_{\theta}(x) \, : \, \theta \in \Theta, P_{\theta}(x) = \int_{-\infty}^{x} p_{\theta}(y) \der y \text{ for all x} \in \Xc \} \) in function space. With this representation of \( \Pc_{\Theta} \) in function space, the parameterization is now essentially the functional relationship \( f \) that pushes each \( \theta \) in \( \Theta \) to a pdf \( p_{\theta}(x) \). If we choose \( f \) to be regular, we can then also transfer additional topological structure from \( \Theta \) to \( \Pc_{\Theta} \). 

\begin{definitie}[Regular parametric model]
We call the model \( \Pc_{\Theta} \) a \emph{regular parametric model}, if the parameterization \( \theta \mapsto p_{\theta}(x) = f(x \, | \, \theta) \), that is, the functional relationship \( f \), satisfies the following conditions %
\begin{enumerate}[label=(\roman*)]
\item its domain \( \Theta \) is an open subset of \( \R^{d} \),
\item at each possible true value \( \theta^{*} \in \Theta \), the spherical representation \( \theta \mapsto m_{\theta}(x) = 2 \sqrt{p_{\theta}(x)}= 2 \sqrt{f(x \, | \, \theta)} \) is so-called Fr\'{e}chet differentiable in \( L_{2}(\lambda) \). The tangent function, i.e., the ``derivative'' in function space, at \( m_{\theta^{*}}(x) \) is then given by 
\begin{align}
\frac{ \der m_{\theta}(x) }{\der \theta} \der \theta = \tfrac{1}{2} (\theta - \theta^{*})^{T} \dot{l}(x \, | \, \theta^{*}) m_{\theta^{*}}(x),
\end{align}
where \( \dot{l}(x \, | \, \theta^{*}) \) is a \( d \)-dimensional vector of score functions in \( L_{2}(P_{\theta^{*}}) \),
\item the Fisher information matrix \( I_{X}(\theta) \) is non-singular,
\item the map \( \theta \mapsto \dot{l}(x \, | \, \theta) m_{\theta}(x) \) is continuous from \( \Theta \) to \( L_{2}^{d}(\lambda) \). %
\end{enumerate}
Note that (ii) allows us to generalize the geometrical concepts discussed in Appendix~{\ref{fisherRaoMetric}} to more general random variables \( X \). \( \hfill \diamond \)
\end{definitie}

We provide some intuition. Condition (i) implies that \( \Theta \) inherits the topological structure of \( \R^{d} \). In particular, we have an inner product on \( \R^{d} \) that allows us to project vectors onto each other, a norm that allows us to measure the length of a vector, and the Euclidean metric that allows us to measure the distance between two vectors by taking the square root of the sums of squares, that is, \( \| \theta^{*} - \theta \|_{2} = \sqrt{\sum_{i=1}^{d} \big (\theta^{*}_{i} - \theta_{i} \big )^{2}} \). For \( d = 1 \) this norm is just the absolute value, which is why we previously denoted this as \( | \theta^{*} - \theta | \). 

Condition (ii) implies that the measurement of distances in \( \R^{d} \) generalizes to the measurement of distance in function space \( L_{2}(\lambda) \). Intuitively, we perceive functions as vectors and say that a function \( h \) is a member of \( L_{2}(\lambda) \), if it has a finite norm (length), i.e., \( \| h(x) \|_{L_{2}(\lambda)} < \infty \), meaning
\begin{align}
\| h(x) \|_{L_{2}(\lambda)} & =%
\begin{cases}
 \sqrt{ \int_{\Xc} [h(x)]^{2} \der x } & \text{ if } X \text{ takes on outcomes on } \R, \\
\sqrt{\sum_{x \in \Xc} [h(x)]^{2} } & \text{ if } X \text{ is discrete.}
 \end{cases}
\end{align}
As visualized in the main text, by considering \( \Mc_{\Theta} = \{ m_{\theta}(x) = \sqrt{p_{\theta}(x)} \, | \, p_{\theta} \in \Pc_{\theta} \} \) we relate \( \Theta \) to a subset of the sphere with radius two in the function space \( L_{2}(\lambda) \). In particular, Section~{\ref{fiInMdl}} showed that whenever the parameter is one-dimensional, thus, a line, that the resulting collection \( \Mc_{\Theta} \) also defines a line in model space. Similarly, Appendix~{\ref{fisherRaoMetric}} showed that whenever the parameter space is a subset of \( [0, 1] \times [0, 1] \) that the resulting \( \Mc_{\Theta} \) also forms a plain. 

Fr\'{e}chet differentiability at \( \theta^{*} \) is formalized as 
\begin{align}
\frac{ \| m_{\theta}(x) - m_{\theta^{*}}(x) - \tfrac{1}{2} (\theta - \theta^{*})^{T} \dot{l}(x \, | \, \theta^{*}) m_{\theta^{*}}(x) \|_{L_{2}(\lambda)}}{ \| \theta - \theta^{*} \|_{2}} \rightarrow 0 .
\end{align}
This implies that the linearization term \( \tfrac{1}{2} (\theta - \theta^{*})^{T} \dot{l}(x \, | \, \theta^{*}) m_{\theta^{*}}(x) \) is a good approximation to the ``error'' \( m_{\theta}(x) - m_{\theta^{*}}(x) \) in the model \( \Mc_{\Theta} \), whenever \( \theta \) is close to \( \theta^{*} \) given that the score functions \( \dot{l}(x \, | \, \theta^{*}) \) do not blow up. More specifically, this means that each component of \( \dot{l}(x  \, | \, \theta^{*}) \) has a finite norm. We say that the component \( \tfrac{\partial } { \partial \theta_{i}} l(x \, | \, \theta^{*}) \) is in \( L_{2}(P_{\theta^{*}}) \), if \( \| \tfrac{\partial } { \partial \theta_{i}} l(x \, | \, \theta^{*}) \|_{L_{2}(P_{\theta^{*}})} < \infty \), meaning 
\begin{align}
\Big \| \tfrac{\partial } { \partial \theta_{i}} l(x \, | \, \theta^{*}) \Big \|_{L_{2}(P_{\theta^{*}})} =
 \begin{cases}
\sqrt{\int_{x \in \Xc} \big ( \tfrac{\partial} { \partial \theta_{i}} l(x \, | \, \theta^{*}) \big )^{2} p_{\theta^{*}}(x) \der x} & \text{if } X \text{ is continuous,} \\
\sqrt{\sum_{x \in \Xc} \big ( \tfrac{\partial} { \partial \theta_{i}} l(x \, | \, \theta^{*}) \big )^{2} p_{\theta^{*}}(x) } & \text{if } X \text{ is discrete.} \\
\end{cases}
\end{align}
This condition is visualized in \refFig{2dGeoC} and \refFig{fiOrthogonal} by tangent vectors with finite lengths. Under \( P_{\theta^{*}} \), each component \( i=1, \ldots, d \) of the tangent vector is expected to be zero, that is, %
\begin{align}
 \begin{cases}
\int_{x \in \Xc} \tfrac{\partial} { \partial \theta_{i}} l(x \, | \, \theta^{*}) p_{\theta^{*}}(x) = 0 & \text{if } X \text{ is continuous,} \\
\sum_{x \in \Xc} \tfrac{\partial} { \partial \theta_{i}} l(x \, | \, \theta^{*}) p_{\theta^{*}}(x) =0 & \text{if } X \text{ is discrete.} \\
\end{cases}
\end{align}
This condition follows from the chain rule applied to the logarithm and an exchange of the order of integration with respect to \( x \), and derivation with respect to \( \theta_{i} \), as
\begin{align}
\int_{x \in \Xc} \tfrac{ \partial }{ \partial \theta_{i}} l(x \, | \, \theta^{*}) p_{\theta^{*}}(x) \der x = \int_{x \in \Xc} \tfrac{ \partial }{ \partial \theta_{i}} p_{\theta^{*}}(x) \der x = \tfrac{ \partial } { \partial \theta _{i} } \int_{x \in \Xc} p_{\theta^{*}}(x) \der x = \tfrac{ \partial } { \partial \theta _{i} } 1 = 0.
\end{align}
Note that if \( \int \tfrac{\partial}{\partial \theta_{i}} p_{\theta^{*}}(x) \der x > 0 \), then a small change at \( \theta^{*} \) will lead to a function \( p_{\theta^{*}+\der  \theta}(x) \) that does not integrate to one and, therefore, not a pdf. 

Condition (iii) implies that the model does not collapse to a lower dimension. For instance, when the parameter space is a plain the resulting model \( \Mc_{\Theta} \) cannot be line. Lastly, condition (iv) implies that the tangent functions change smoothly as we move from \( m_{\theta^{*}}(x) \) to \( m_{\theta}(x) \) on the sphere in \( L_{2}(\lambda) \), where \( \theta \) is a parameter value in the neighborhood of \( \theta^{*} \). 

The following conditions are stronger, thus, less general, but avoid Fr\'{e}chet differentiability and are typically easier to check.

\begin{lemma}
Let \( \Theta \subset \R^{d} \) be open. At each possible true value \( \theta^{*} \in \Theta \), we assume that \( p_{\theta}(x) \) is continuously differentiable in \( \theta \) for \( \lambda \)-almost all \( x \)  with tangent vector \( \dot{p}_{\theta^{*}}(x) \). We define the score function at \( x \) as %
\begin{align}
\label{newScore}
 \dot{l}(x  \, | \, \theta^{*}) = \frac{\dot{p}_{\theta^{*}}(x)}{  p_{\theta^{*}}(x) } 1_{[p_{\theta^{*}} > 0]}(x),
 \end{align}
 where \( 1_{[p_{\theta^{*}} > 0]}(x) \) is the indicator function 
\begin{align}
1_{[p_{\theta^{*}} > 0]}(x) = \begin{cases}
1 & \text{ for all } x \text{ such that } p_{\theta^{*}}(x) > 0, \\
0 & \text{ otherwise} .
 \end{cases}
\end{align}
The parameterization \( \theta \mapsto P_{\theta} \) is regular, if the norm of the score vector \refEq{newScore} is finite in quadratic mean, that is, \( \dot{l}(X \, | \, \theta^{*})  \in L_{2}(P_{\theta^{*}}) \), and if the corresponding Fisher information matrix based on the score functions \refEq{newScore} is non-singular and continuous in \( \theta \). \( \hfill \diamond \)
\end{lemma}
There are many better sources than the current manuscript on this topic that are mathematically much more rigorous and better written. For instance, \citet{bickel1993efficient} give a proof of the lemma above and many more beautiful, but sometimes rather (agonizingly) technically challenging, results. For a more accessible, but no less elegant, exposition of the theory we highly recommend \citet{van1998asymptotic}. 
\end{document}